\documentclass[a4paper,10pt]{article}
\usepackage[utf8]{inputenc}

\usepackage{amscd,  graphicx, color, mathrsfs}
\usepackage{lineno}
\usepackage{dsfont}

\usepackage{amsmath, amsthm, amssymb, amsfonts, enumerate}
\usepackage[colorlinks=true,linkcolor=blue,urlcolor=blue]{hyperref}

\usepackage{amsfonts,latexsym,amstext}
\usepackage[active]{srcltx}
\usepackage{xcolor}
\usepackage{amsmath,amsfonts,calrsfs, amssymb,color}
\usepackage{amsthm}
\usepackage{amssymb}

\newtheorem{theorem}{Theorem}[section]
\newtheorem{lemma}[theorem]{Lemma}
\newtheorem{proposition}[theorem]{Proposition}
\newtheorem{definition}[theorem]{Definition}
\newtheorem{example}[theorem]{Example}
\newtheorem{hypothesis}[theorem]{Hypothesis}

\newtheorem{remark}[theorem]{Remark}
\newtheorem{corollary}[theorem]{Corollary}

\numberwithin{equation}{section}

\def\qed{{\hfill\hbox{\enspace${ \square}$}} \smallskip}
\def\sqr#1#2{{\vcenter{\vbox{\hrule height .#2pt \hbox{\vrule
 width .#2pt height#1pt \kern#1pt \vrule
width .#2pt} \hrule height .#2pt}}}}
\def\square{\mathchoice\sqr54\sqr54\sqr{4.1}3\sqr{3.5}3}

\def\ds{\begin{displaystyle}}
\def\eds{\end{displaystyle}}
\def\dis{\displaystyle }
\def\<{\langle }
\def\>{\rangle }

\def\R{\mathbb R}
\def\N{\mathbb N}

\def\E{\mathbb E}
\def\P{\mathbb P}
\def\Q{\mathbb Q}

\def\calf{{\cal F}}
\def\calg{{\cal G}}

\def\calm{{\cal M}}

\def\call{{\cal L}}
\def\cals{{\cal S}}

\allowdisplaybreaks

\title{A nonlinear Bismut-Elworthy formula for HJB equations with quadratic Hamiltonian in Banach spaces}

\usepackage{authblk}

\author[1]{Davide ADDONA\thanks{davide.addona@unimib.it}
} 
\author[1]{Elena BANDINI\thanks{elena.bandini@unimib.it}} 
\author[1]{Federica MASIERO\thanks{federica.masiero@unimib.it}} 
\affil[1]{Dipartimento di Matematica e Applicazioni, Universit\`a di Milano-Bicocca, Milano, Italy}

\date{}

\begin{document}
\maketitle

\begin{abstract}
We consider a Backward Stochastic Differential Equation (BSDE for short) in a Markovian framework for the pair of processes $(Y,Z)$, with generator with quadratic growth with respect to $Z$. The forward equation is an evolution equation in an abstract Banach space. 
 We prove an analogue of the Bismut-Elworty formula when the diffusion operator has a pseudo-inverse not necessarily bounded and  when the generator has quadratic growth with respect to $Z$. In particular, our model covers the case of the heat equation in space dimension greater than or equal to 2.  We apply these results to solve semilinear Kolmogorov equations for the unknown $v$, with nonlinear term with quadratic growth with respect to $\nabla v$  and  final condition only bounded and continuous, and to solve stochastic optimal control problems with quadratic growth.
\end{abstract}


\medskip


 \section{Introduction}\label{Sec_introduction}

In this paper  we  deal with  Markovian BSDEs  whose  generator has quadratic growth with respect to $Z$, and we  generalize to this framework  the Bismut-Elworthy type formula introduced in  \cite{futeBismut}, where the Lispchitz case was studied.
More precisely, our BSDE is related to a forward stochastic differential equation of the form
\begin{equation}
\left\{
\begin{array}
[c]{l}%
dX^{t,x}_\tau  =AX^{t,x}_\tau d\tau+F( X^{t,x}_\tau)d\tau+(-A)^{-\alpha}dW_\tau
,\text{ \ \ \ }\tau\in\left[  t,T\right],  \\
X^{t,x}_t =x\in E,
\end{array}
\right.  \label{forwardintro}%
\end{equation}
where $E$ is a Banach space  which is continuously and densely embedded in a real and separable Hilbert space $H$. The operator $A$ is the generator of a contraction analytic semigroup in $H$, 
 which turns out to be strongly continuous or analytic in $E$, and  $\left\lbrace W_\tau,\,\tau\geq 0\right\rbrace$
is a cylindrical Wiener process in $H$.
We assume that the stochastic convolution%
\[
w_{A}(\tau) =\int_{0}^{\tau}e^{( \tau-s) A}(-A)^{-\alpha}dW_{s}
\]
is well defined as a Gaussian process in $H$, and that it admits an $E$-continuous version.

The presence of the diffusion operator $(-A)^{-\alpha }$ in \eqref{forwardintro} allows us to deal with stochastic heat equations in $2$ and $3$ space dimensions, while stochastic heat equations in one space dimension can be considered without any regularization of the white noise, 
that is in the case with $\alpha=0$. Moreover, we consider dissipative maps $F$ in  \eqref{forwardintro} in order to have more generality in the structure of the equation. Notice that, under this latter assumption,  $F$ is well defined only on the Banach space $E$, 
while  it is not even defined on the whole Hilbert space $H$; this is a natural situation arising in many evolution equations, see e.g. \cite{dpz92} and \cite{dpzErgodicity}.

The solution of equation (\ref{forwardintro}) will be denoted by $X$, or also by $X^{t,x}$, to stress the dependence on the initial conditions, and the transition semigroup related to $X^{t,x}$ will be denoted by
\[
 P_{t,\tau}[\phi](x):=\E\phi(X_\tau^{t,x}),\quad \phi \in B_b(E).
\]
At least formally, the generator of $P_{t,\tau}$ is the second order differential operator
\[
 (\call f)(x)=\frac{1}{2}(Tr ((-A)^{-\alpha}(-A^*)^{-\alpha} \nabla^2 f)(x)+\<Ax,\nabla f(x)\>+\<F(x),\nabla f(x)\>.
\]
This is the link with the solution, in mild sense, of the semilinear Kolmogorov equation in $E$ (see e.g. \cite{DP3}):
\begin{equation}
\left\{
\begin{array}
[c]{l}%
\frac{\partial v}{\partial t}(t,x)=-\call v\left(  t,x\right)
+\psi\left( t,x,v(t,x),\nabla v(t,x) (-A)^{-\alpha}  \right)  ,\text{ \ \ \ \ }t\in\left[  0,T\right]
,\text{ }x\in E,\\
v(T,x)=\phi\left(  x\right).
\end{array}
\right.  \label{KolmoIntro}%
\end{equation}
We recall that by mild solution of equation (\ref{KolmoIntro}) we mean a bounded and continuous function $v:[0,T]\times H \rightarrow H$, once G\^ateaux differentiable with respect to $x$, and
satisfying the integral equality
\begin{equation} 
v(t,x)=P_{t,T}\left[  \phi\right]  \left(  x\right)  +\int_{t}^{T}%
P_{t,s}\left[  \psi(s,\cdot, v(s,\cdot), \nabla v\left(  s,\cdot\right)(-A)^{-\alpha}) 
\right]  \left(  x\right)  ds, 
\text{\ \ }t\in\left[  0,T\right],\text{ }x\in E.\label{solmildkolmointro}%
\end{equation}
Second order differential equations are a widely studied topic in the literature, see e.g. \cite{DP3}.
In the case of $\psi$ only locally Lipschitz continuous, we cite \cite{Go2}, \cite{MR}, \cite{Mas1} and also \cite{Mas-spa}, where in
particular the quadratic case is studied with datum $\phi$ only continuous. We also mention the monograph \cite{Ce}, where semilinear Kolmogorov equations related to forward equations of reaction diffusion type more general than the one considered here are studied, but requiring Lipschitz continuity of the final datum.

We will consider equation \eqref{KolmoIntro} under the assumptions that the final datum $\phi$ is bounded and continuous,
and that $\psi$ has quadratic growth with respect to the derivative $\nabla v{(-A)^{-\alpha}}$.
In order to prove existence and uniqueness of a mild solution of the form \eqref{solmildkolmointro} for the Kolmogorov equation (\ref{KolmoIntro}), we aim at   representing this mild solution  in terms of a Markovian BSDE of the form
\begin{equation}\label{bsdeIntro}
    \left\{\begin{array}{l}
 dY_\tau=-\psi(\tau,X_\tau,Y_\tau,Z_\tau)\;d\tau+Z_\tau\;dW_\tau,
  \\\dis
  Y_T=\phi(X_T).
\end{array}\right.
\end{equation}
We recall that, in order to solve partial differential equations by means of BSDEs, one of the crucial tasks is the identification of $Z $ with the derivative of $Y$ taken in the directions  of the diffusion operator. In this regard, we refer to the seminal paper \cite{PaPe92} for the finite dimensional case, and to \cite{fute} for the infinite dimensional extension in Hilbert spaces: in both papers the driver $\psi $ is Lipschitz continuous  in $Y$ and in $Z$, and  $\psi$ and $\phi$ are differentiable. We also mention \cite{MasBan}, where  an extension to the Banach space case is studied  with the same assumptions of Lipschitz continuity and differentiability on the data.

In the present paper we do not make differentiability assumptions on the  coefficients: thank to a  variant of the nonlinear Bismut-Elworthy formula for BSDEs introduced in \cite{futeBismut}, 
 we are still able to prove that the solution of the BSDE (\ref{bsdeIntro}) gives the mild solution of the Kolmogorov equation (\ref{KolmoIntro}). 
Bismut-Elworthy formulas for the transition semigroup of equations of type \eqref{forwardintro} with invertible diffusion operator are a classical topic in the literature, see e.g.  \cite{DP3}. In \cite{Ce}  the case of an operator like the one in (\ref{forwardintro}), with pseudo-inverse which is not necessarily bounded,  is also considered.
According to these classical Bismut formulas, for every $0\leq t<\tau\leq T,\,x\in H$, $h\in H$, and for every bounded and continuous real function $f$ defined on $H$, one has
 \begin{equation}\label{Bismutintro}
 \<\nabla_x P_{t,\tau}[f](x), h\>=\E f\left(X_\tau^{t,x}\right)U^{h,t,x}_\tau,
\end{equation}
where ($G(r,X^{t,x}_r)$ being the general diffusion operator)
\begin{equation*}
 U^{h,t,x}_\tau:=\dfrac{1}{\tau-t}\int_t^\tau \<G^{-1}(r,X^{t,x}_r)\nabla_xX_r^{t,x}h,dW_r \>.
\end{equation*}
\newline In \cite{futeBismut}  a nonlinear Bismut-Elworthy formula for the process $Y$ solution of the BSDE (\ref{bsdeIntro})  is proved when $\psi $ is Lipschitz continuous with respect to $Z$ and the process $X$ takes its values in a Hilbert space $H$. According to this formula, for $0\leq t<\tau\leq T,\,x\in H$, for every direction $h\in H$,
\begin{equation}\label{Bismut-nonlin-intro}
 \E\left[ \nabla_x\,Y^{t,x}_\tau h \right]=
\E\Big[\int_\tau^T\psi\left(r,X_r^{t,x},Y_r^{t,x},Z_r^{t,x}\right)U^{h,t,x}_r\,dr\Big]
+\E\left[  \phi(X_T^{t,x})U^{h,t,x}_T\right].
\end{equation}
Formula (\ref{Bismut-nonlin-intro}) is used in  \cite{futeBismut} to solve a semilinear Kolmogorov equation of the form of (\ref{KolmoIntro}). When the Hamiltonian function $\psi$ is Lipschitz continuous with respect to the derivative of $ v$, semilinear Kolmogorov equations of the type of \eqref{KolmoIntro} can be solved also by using the estimates coming from  the classical Bismut formulas (\ref{Bismutintro}) and by  a fixed point argument, see e.g. \cite{Ce}, \cite{DP3},
\cite{Go1}. In the quadratic case this procedure does not work anymore:
for this reason,  nonlinear versions of Bismut-Elworthy formulas,  that give an alternative way to solve equations like  \eqref{KolmoIntro}, are particularly interesting in such a framework.
In \cite{Mas-spa},
a nonlinear version of the Bismut-Elworthy formula has been  provided   and  has been applied to semilinear Kolmogorov equations of the type of (\ref{KolmoIntro}), with quadratic hamiltonian, and 
in a Hilbert space.

In the present paper, we generalize \eqref{Bismut-nonlin-intro} to the Banach space framework, and to the case of diffusion operator $(-A)^{-\alpha}$ that has unbounded pseudo-inverse operator.
In this context, the   nonlinear Bismut formula (\ref{Bismut-nonlin-intro}) has  its own independent interest, and  moreover it allows to solve the Kolmogorov equation with Hamiltonian function quadratic with respect to $\nabla v(-A)^{-\alpha}$. 
We first  provide an analogous of the nonlinear Bismut formula given in \cite{futeBismut} in the case of   Banach space framework and  Lispchitz continuous generator.
Then, 
we prove a nonlinear Bismut formula in the quadratic case when $\psi$ and $\phi$ are  differentiable. 
To this end,  denoted by $(Y^{t,x},Z^{t,x})$ a solution to the Markovian BSDE (\ref{bsdeIntro}) and assuming that $\phi$ and $\psi$ are differentiable, the two  main ingredients  are 
the identification
\begin{equation}\label{identif-Z-intro}
  Z^{t,x}_t=\nabla_x\,Y^{t,x}_t (-A)^{-\alpha},\quad t\in[0,T], \,x\in E, 
\end{equation} 
and an a priori estimate on $Z^{t,x}$ 
of the form ($C$ being a constant depending on $t,\;T,\;A,\;F,\;\Vert\phi\Vert_\infty$)
\begin{equation}\label{stimabismutintro}
\vert Z^{t,x}_t\vert_H\leq C (T-t)^{-1/2},
\end{equation}
which is obtained with techniques similar to the ones used in \cite{Mas-spa}, see also \cite{BaoDeHu} and \cite{Ri}. Both 
\eqref{identif-Z-intro} and \eqref{stimabismutintro} are new in the Banach space framework and in the case of quadratic generator with respect to $z$.
 Finally, differentiability assumptions are removed
by an approximation procedure, obtained by suitably generalizing the one introduced in \cite{PZ}.

Our results can be applied to a stochastic optimal control problem consisting in  minimizing a cost functional of the form
\begin{equation}\label{contIntro}
J\left(  t,x,u\right) 
=\Big[\mathbb{E}\int_{t}^{T}
l\left(s,X^{u}_s,u_s\right)ds+\mathbb{E}\phi\left(X^{u}_T\right)\Big]
\end{equation}
over all the admissible controls $u$  taking values in $H$ and   not necessarily bounded.
Here $l$ has quadratic growth with respect to $u$, and $X^{u}$
is the solution of  the controlled state equation 
\begin{equation}
\left\{
\begin{array}
[c]{l}%
dX^{u}_\tau =AX^{u}_\tau d\tau +F(X^u_\tau)d\tau+ Qu_\tau  d\tau+(-A)^{-\alpha}dW_\tau ,\text{ \ \ \ }\tau\in\left[  t,T\right] \\
X^{u}_t  =x, 
\end{array}
\right.  \label{sdecontrolintro}%
\end{equation}
with $Q=I$ or $Q=(-A)^{-\alpha}$.
The aim of this latter part of the work is to characterize the value function as the solution of the associated Hamilton Jacobi Bellman (HJB in the following) equation, and to provide a feedback law for optimal controls.
  If $Q=(-A)^{-\alpha}$, namely when the controls affect the system only through the noise (the so called structure condition holds true),   the optimal control problem \eqref{contIntro} can be  completely solved, see Theorem \ref{thm_rel_fond_opt_cont_str}. 
When $Q=I$,  the optimal control problem can be  completely solved 
by restricting ourselves to the class of more regular controls taking values in $D((-A)^{-\alpha})$, see Theorem \ref{thm_rel_fond_opt_cont_spec_cost}.
 In the general case of $Q=I$ and $H$-valued controls, we are able to provide an  ``$\varepsilon$-optimal solution" of the problem in the sense that the value function can be approximated by a sequence of functions which are solutions of approximating HJB equations, and we can obtain an  $\varepsilon$-optimal control in feedback form, see Theorem \ref{T:convJu_n}.

The paper is organized as follows: in Section \ref{Sec:2} we fix the notations and we give the results on the forward process. In Section \ref{Sec_forwbac} we introduce the forward backward system: here  the main results are 
the identification \eqref{identif-Z-intro} of $Z_t^{t,x}$ with $\nabla Y^{t,x}_t(-A)^{-\alpha}$, which is new in the case of $\psi$ quadratic with respect to $z$ and in the Banach space framework, and 
the a  priori estimate \eqref{stimabismutintro} on $Z$ not involving derivatives of the coefficients of the BSDE. 
In Section \ref{sez-Bismut-lip} we give the nonlinear Bismut formula (\ref{Bismut-nonlin-intro}) in the Banach space $E$ and with $\psi$ Lipschitz continuous with respect to $z$, then in Section \ref{sez-Bismut-quad} we extend  formula (\ref{Bismut-nonlin-intro}) to the case of $\psi $ quadratic with respect to $z$. In both  Sections \ref{sez-Bismut-lip} and \ref{sez-Bismut-quad}, the Bismut formula is applied to solve the corresponding semilinear Kolmogorov equation (\ref{KolmoIntro}). Finally in Section \ref{sez-appl-contr} we apply the previous results to solve the stochastic optimal control problem \eqref{contIntro}.

 \section{Notations and preliminary results on the forward process}\label{Sec:2}

We assume that $E$  is a real and separable Banach space which admits a Schauder basis, and that $E$ is continuously and
densely embedded in a real and separable Hilbert space $H$.
$E$ and $H$ are respectively endowed with the norms $|\cdot|_E$ and $|\cdot|_H$.
We fix a complete probability space  $(\Omega, \mathcal F, \P)$ endowed with a filtration $\{\mathcal F_t, \,\,t \geq 0\}$ satisfying the usual conditions.

We list below some notations that are used in the paper.  Let $K$ be a given Banach space endowed with the norm $|\cdot|_K$. For any $p,q \in [1,\infty)$ and any $t \in [0,\,T]$, we set 
\begin{itemize}
\item $L^p(0,T;K)$ the space of  
$K$-valued measurable 
functions defined on $[0,\,T]$, normed by
\begin{align*}
\|f\|_{L^p(0,T;K)}:=\Big(\int_0^T|f_s|^p_K ds\Big)^{1/p}.
\end{align*}
\item $L^q(\Omega;L^p(0,T;K))$ the space of adapted processes $(u_s)_{s \in [0,\,T]}$, defined on $[0,\,T]$ and  with values in $K$,  normed by
\begin{align*}
\|u\|_{L^q(\Omega;L^p(0,T;K))}:=\Big[\mathbb E\Big (\int_0^T|u_s|^p_K ds\Big)^{q/p}\Big]^{1/q}.
\end{align*} 
\item ${\cals^p}((t,T];K)$ (resp. ${\cals^p}([t,T];K)$)  the space of all  adapted processes $(X_s)_{s\in [t,T]}$, continuous on $(t,T]$ (resp. on $[t,T]$) and   with values in $K$,   normed by 
\begin{align*}
\|X\|_{{\cals^p([t,T];K)}}=\|X\|_{{\cals^p((t,T];K)}}:=\mathbb E\Big[\sup_{s\in[t,T]}|X_s|_K^p\Big]^{1/p}.
\end{align*}
If $K = \R$ we simply write ${\cals^p([t,T])}$. 
\item ${\calm^p}([t,T];K)$  the space of all predictable processes $(Z_s)_{s\in [t,T]}$ with values in $K$ normed by
\begin{align*}
\|Z\|_{\calm^p([t,T];K)}:=\mathbb E\Big[\Big(\int_t^T|Z_s|_K ^2ds\Big)^{p/2}\Big]^{1/p}.
\end{align*}
If $K = \R$ we simply write ${\calm^p([t,T])}$. 
\end{itemize}

\noindent We  denote by $L(E,K)$ the space of all bounded linear operators from $E$ to $K$, endowed with the usual operator norm. $E^\ast$  denotes the dual space of $E$, and  $\langle \cdot, \cdot\rangle_{E\times E^*}$ denotes the duality between $E$ and $E^\ast$. 

\noindent  We say that a function $f : E \rightarrow K$ belongs to the class $\mathcal G^1(E,K)$ if $f$ is continuous and G\^ateaux differentiable on $E$ and if the gradient $\nabla f : E \rightarrow L(E,K)$
 is strongly continuous. If $K=\R$ we  simply write $\mathcal G^1(E)$.
We  say  that $f : [0,\,T] \times E \rightarrow \R$ is in $\calg^{0,1}\left(  \left[  0,T\right]  \times E\right)$ if $f$ is continuous and G\^ateaux differentiable with respect to every $x \in E$ and  the gradient $\nabla f : [0,\,T] \times E \rightarrow L(E,\R)$ is strongly continuous. For more details on this classes of G\^ateaux differentiable functions see  \cite[Section 2.2]{fute}.

\subsection{The forward equation}

We are given the Markov process $X$ in $E$ (also denoted $X^{t,x}$ to stress the dependence
on the initial conditions)  solution to the equation 
\begin{equation}
\label{forward}
\left\{
\begin{array}{ll}
dX_\tau^{t,x}=AX_\tau^{t,x}d\tau+F(X_\tau^{t,x})d\tau+(-A)^{-\alpha} dW_\tau, & \tau\in[t,T], \vspace{1mm} \\
X_t^{t,x}=x\in E,
\end{array}
\right.
\end{equation}
where $(W_\tau)_{\tau\in [0,T]}$ is a cylindrical Wiener process with values in $H$, see e.g. \cite{dpz92} for details on cylindrical Wiener processes in infinite dimensions. From now on $\{\mathcal{F}_{\tau}, \tau \geq 0\}$ will be the natural filtration generated by
the Wiener process and augmented in the usual way.

\noindent We assume the following on the coefficients of equation \eqref{forward}.
\begin{hypothesis}
 \label{ip_forward}
\begin{enumerate}
\item $A$ is a linear operator which generates a contraction analytic semigroup $(e^{t A})_{t\geq0}$ on the Hilbert space $H$ and there exist $c,\omega>0$ such that $|e^{tA}h|_H\leq ce^{-\omega t}|h|_H$ for any $h\in H$ and any $t\geq 0$. Further, the restriction of $A$ to $E$ generates  a contraction $C_0$ (or analytic) semigroup on $E$. 
\item The stochastic convolution
\begin{align*}
w^A(s,t):=\int_s^te^{(t-u)A}(-A)^{-\alpha}dW_u, \quad 0\leq s<t\leq T,
\end{align*}
admits an $E$-continuous version, and, for any $p\geq 2$,
$\mathbb E[\sup_{t\in[0,T]}|w^A(t)|_E^p]<+\infty$
(when $s=0$ we write $w^A(t)$ instead of $w^A(0,t)$).
\item $F:D(F)\subset H\rightarrow H$ is a measurable and dissipative map, and $E\subseteq D(F)$.
\item The restriction $F_E$ of $F$ to $E$  is a map from $E$ to $E$ which is measurable and dissipative
 (where no confusion is possible, we simply write $F$ instead of $F_E$).
 $ F \in \calg^1(H,H)$ and $F_E$ is Fr\'echet differentiable. Further, there exist $a,c,\gamma>0$, $m\in \N$ and for any $z\in E$ an element $z^*\in\partial |z|_E$, such that, for any $x \in E$, $h\in H$, 
 \begin{align*}
 |F_E(x)|_{E} & \leq c(1+|x|_E^{2m+1}), \\
\|\nabla F(x)\|_{\mathcal L(E)} & \leq c(1+|x|_E^{2m}),  \\
\langle F(x+z)-F(x),z^*\rangle_{E\times E^*} & \leq -a|z|_E^{2m+1}+c(1+|x|^\gamma_E),  \\
|\nabla F(x) h|_H & \leq c\left(1+|x|_E^{2m}\right)|h|_H.
\end{align*}
\item $\alpha\in(0,1/2)$.
\end{enumerate}
\end{hypothesis}
By Hypothesis \ref{ip_forward}-1. and the Kuratowski
theorem, see e.g. \cite{Par}, Chapter I, Theorem 3.9,  it follows that $E$ is a
Borel set in $H$.


\begin{remark}\label{R:Hyp_1}
Since by Hypothesis  \ref{ip_forward}-3.-4. $F$ is differentiable and dissipative, 
 we get
\begin{align*}
1\geq |z-\alpha DF(x)z|_E, \quad x,z\in E, \quad |z|_E=1, \quad \alpha>0.
\end{align*}
In particular, from the Hahn-Banach theorem,  there exists $z^*\in \partial|z|_E$ such that $|z-\alpha DF(x)z|_E=\langle z-\alpha DF(x)z,z^*\rangle_{E\times E^*}$, and therefore
$\langle DF(x)z,z^*\rangle _{E\times E^*}\leq 0$. 
Further, from \cite[Appendix D]{dpz92} we have 
\begin{align}
\label{der_-}
D_-|x|_Ey=\min\{\langle y,x^*\rangle_{E\times E^*}: x^*\in\partial |x|_E \}.
\end{align}
\end{remark}

\begin{remark}
\label{R2:hyp_1}
 Since $A$ generates a contraction semigroup on $E$, then $A$ is dissipative, and for any $x\in D(A)$ we have
$\langle Ax,x^*\rangle_{E\times E^*}\leq 0$,  $x^*\in \partial |x|_E$,
see Example D.8 in \cite{dpz92}.
\end{remark}
We now give an example of spaces $E$ and $H$ and of operator $A$ satisfying Hypothesis \ref{ip_forward}-1.-2.
\begin{example}
Let $d,n\in\N$ with $d\leq3$,  $\mathcal O\subset \R^d$ be an open bounded set,  $H:=L^2(\mathcal O;\R^n)$ and $E:=C(\overline{\mathcal O};\R^n)$. Further, let $A$ be the realization in $H$ of the operator 
\begin{align*}
\mathcal A-(\rho+1)I=(\Delta-(\rho+1)I,\ldots,\Delta-(\rho+1)I),
\end{align*}
with boundary conditions $\mathcal Bu=0$, where $\mathcal B=(\mathcal B_1,\ldots,\mathcal B_n)$ and
\begin{align*}
\mathcal B_k=Id, \quad \textrm{or }\,\,\,\,\mathcal B_k:=\sum_{i=1}^d\nu_i(\xi)\frac{\partial}{\partial\xi_i}, \quad \xi\in\partial\mathcal O, \ k=1,\ldots,n,
\end{align*} 
where $\nu_i$ is the normal vector to the boundary of $\overline{\mathcal O}$. As shown for example in \cite{LU95}, $A$ satisfies Hypothesis \ref{ip_forward}-$1.$ Moreover, \cite[Lemma 6.1.2]{Ce} with $Q=(-A)^{-\alpha}$ shows that Hypothesis \ref{ip_forward}-$2.$ is satisfied  with this choice of $H$, $E$ and $A$.
\end{example}
In the following proposition we collect important results on the solution of the forward equation \eqref{forward}.
 We recall that, given $x\in E$ and $t\in[0,T]$, a mild solution to \eqref{forward} is an adapted process $X^{t,x}:[0,T]\times \Omega\rightarrow E$ which satisfies
\begin{align}
\label{forward_mild_sol}
X_\tau^{t,x}=e^{(\tau-t)A}x+\int_t^\tau e^{(\tau-s)A}F(X_s^{t,x})ds+\int_t^\tau e^{(\tau-s)A}dW_s, \quad  \tau\in[t,T], \ \P\textup{-a.s.}
\end{align}
\begin{proposition}\label{P:unificata}
Let Hypothesis \ref{ip_forward} hold true. Then the following hold. 
\begin{itemize}
\item[(i)] For any $x\in E$,  $t\in [0,T]$, the problem \eqref{forward} admits a unique mild solution $X^{t,x}\in \mathcal S^p((t,T];E)$, for any $p\geq1$. If  $A$ generates a strongly continuous semigroup on $E$, then  the process $X^{t,x}$ is also continuous up to $t$. Moreover, 
there exists a positive constant $c$ 
 such that,  
 for any $\tau\in[t,T]$, 
\begin{align}
|X_\tau^{t,x}|_E&\leq e^{c\tau}|x|_E+h(t,\tau),\quad \P\textup{-a.s.},\label{forward_estimate}
\end{align}
where
\begin{align*}
& h(t,\tau):=ce^{c(\tau-t)}\int_t^{\tau}\left(1+|w^A(t,s)|_E^{2m+1}\right) ds+\sup_{s\in[t,\tau]}|w^A(t,s)|_E. 
\end{align*}
\item[(ii)] For any $x\in E$, $t\in [0,T]$, the mild solution $X^{t,x}$ to \eqref{forward} is G\^ateaux differentiable as a map from $E$ to $\mathcal S^p([t,T];E)$, 
and
\begin{align}
&\sup_{x\in E, \tau\in[t,T]}|\nabla_x X_\tau^{t,x} z|_E\leq |z|_E, \quad  z\in E, \quad \P\textup{-a.s.} \label{estimate_der}
\end{align}
Moreover, $X^{t,x}$ is G\^ateaux differentiable as a map from $E$ to $\mathcal S^p([0,T];H)$, and
\begin{align}
\label{estimate_der_H}
\sup_{x\in E, \tau\in[t,T]}|\nabla_x X_\tau^{t,x} h|_H\leq |h|_H, \quad h \in H, \quad  \P\textup{-a.s.} 
\end{align}
\item[(iii)] For any $x\in E$, $t\in[0,T]$ and  $\tau\in[t,T]$, \begin{align}
\label{eq_mild_der_H}
\nabla_xX^{t,x}_\tau h=e^{(\tau-t)A}h+\int_t^\tau e^{(\tau-s)A}\nabla F(X_s^{t,x})\nabla_xX^{t,x}_s h \, ds, \quad h \in H, \quad \P\textup{-a.s.} 
\end{align}
\end{itemize}
\end{proposition}
\proof
Item $(i)$ can be proved arguing as in \cite[Theorem 7.13]{dpz92}.

\noindent The first part of $(ii)$ and inequality \eqref{estimate_der} follow from \cite[Propositions 3.10 \& 3.13]{MasBan}. We claim that 
\begin{align}\label{claim2}
\sup_{x\in E, \tau\in[t,T]}|\nabla_x X_\tau^{t,x} z|_H\leq |z|_H, \,\,\,\,z \in E\quad \P\textup{-a.s.}
\end{align}
 If the claim is true, since $E$ is densely embedded into $H$, by approximation we immediately deduce \eqref{estimate_der_H} for any $h\in H$.
In order to prove \eqref{claim2}, we consider $z\in E$ and the approximating processes
$G^n_\tau z:=nR(n,A)\nabla_x X^{t,x}_\tau z$, $n\in\N$, where $R(n,A):=(nI-A)^{-1}$. Then, $G_\tau^n z$ is a strict solution to
\begin{align*}
\frac{d}{dt}G_\tau^{n}z=AG_\tau^{n}z+\nabla F(X_\tau^{t,x})G_\tau^{n}z, \quad \tau\in(t,T], \quad G_t^{n}z=n R(n,A)z.
\end{align*}
The dissipativity of $F$ and $A$ implies
$\frac{d}{d\tau}|G^n_\tau z|_H^2\leq 0$, 
which gives $|G_\tau^nz|_H\leq |nR(n,A) z|_H$. Letting $n\rightarrow+\infty$ we get \eqref{claim2}.

It remains  to prove $(iii)$. To this end,  we recall that (see e.g. \cite{MasBan}), for any $x,z\in E$, the process $\nabla_xX^{t,x}_\tau z$ is a mild solution to 
\begin{align}
\label{first_variation}
\left\{
\begin{array}{ll}
d \zeta_{\tau}=A\zeta_{ \tau} d\tau+\nabla F(X_{\tau}^{t,x})\zeta_{\tau}, & \tau\in [t,T], \vspace{1mm}\\
\zeta_t=z\in E,
\end{array}
\right.
\end{align}and therefore
\begin{align}
\label{mild_sol_der_E}
\nabla_x X^{t,x}_\tau z=e^{(\tau-t)A}z+\int_t^\tau e^{(\tau-s)A}\nabla F(X_s^{t,x})\nabla_xX^{t,x}_s z \, ds, \quad \P\textup{-a.s.}
\end{align}
Let $h\in H$ and let $(h_n)\subset E$ be an approximating sequence of $h$ in $H$. If we replace $h_n$ to $z$ in \eqref{mild_sol_der_E}, from $(ii)$ we deduce that the left-hand side of \eqref{mild_sol_der_E} and the first term in the right-hand side of \eqref{mild_sol_der_E} converge respectively to $\nabla_xX^{t,x}_\tau h$ and to $e^{(\tau-t)A}h$, as $n \rightarrow+\infty$. As far as the integral in the right-hand side of \eqref{mild_sol_der_E} is considered, with $z$ replaced by $h_n$, again from $(ii)$ we infer that
\begin{align*}
e^{(\tau-s)A}\nabla F(X_s^{t,x})\nabla_xX^{t,x}_sh_n\rightarrow e^{(\tau-s)A}\nabla F(X_s^{t,x})\nabla_xX^{t,x}_sh, \quad \P\textup{-a.s.},
\end{align*}
as $n\rightarrow+\infty$. Thanks to Hypothesis \ref{forward}-$4.$, estimate \eqref{forward_estimate} and \eqref{estimate_der_H}, we can apply the dominated convergence theorem and therefore
\begin{align*}
\int_t^\tau e^{(\tau-s)A}\nabla F(X_s^{t,x})\nabla_xX^{t,x}_sh_n\rightarrow \int_t^\tau e^{(\tau-s)A}\nabla F(X_s^{t,x})\nabla_xX^{t,x}_sh, \quad \P\textup{-a.s.},
\end{align*}
as $n\rightarrow+\infty$, which gives \eqref{eq_mild_der_H}.
\endproof

Now we show that for any $x,z\in E$ and  any $t\in[0,T]$, the process $\nabla_x X^{t,x}z$ belongs to $D(-A)^{1/2}$ a.e. in $(t,T)$ and $\P$-a.s., and satisfies useful estimates.
\begin{proposition}\label{P1.9} Let Hypothesis \ref{ip_forward} holds true,  and let $x\in E$,  $z\in H$ and $t\in[0,T]$. Then,  $\nabla_xX^{t,x}z\in D((-A)^{1/2})$, a.e. in $(t,T)$ and $\P$-a.s., and there exists a positive constant $C$ such that, for any $\varepsilon\in [0,1/2]$,
\begin{align}
\label{weight_estimate_derivative}
\int_t^\tau|(-A)^{\varepsilon}\nabla_x X_s^{t,x}z|^2_Hds \leq C(\tau-t)^{1-2\varepsilon}|z|^2_H, \,\,\,\,\tau \in [t,\,T],\quad \P\textup{-a.s.}
\end{align}
\end{proposition}
\proof Let $x \in E$. We  prove \eqref{weight_estimate_derivative}  for $t=0$, 
 the case  $t\in[0,T]$ can be proved by   analogous computations.
  
\noindent We first assume that   $z\in E$. Let $\nabla_x X_t^x z$ be a strict solution to \eqref{first_variation}, otherwise we can approximate it by smooth processes, as in  the proof of  Proposition \ref{P:unificata}, item (ii). The dissipativity of $F$ in $H$ gives
\begin{align*}
\frac d{ds}|\nabla_x X_s^xz|_H^2=\langle A\nabla_xX_s^x z,\nabla_xX_s^x z\rangle_H
+ \langle \nabla F(X_s^x)\nabla_xX_s^x z, \nabla_xX_s^x z\rangle_H
\leq \langle A\nabla_xX_s^x z,\nabla_xX_s^x z\rangle_H,
\end{align*}
for any $s\in[0,T]$. Integrating between $0$ and $\tau\in [0,T]$ 
we get
\begin{align*}
|\nabla_x X_\tau^x z|_H^2+\int_0^\tau\langle -A\nabla_xX_s^x z,\nabla_xX_s^x z\rangle_H ds\leq |z|^2_H.
\end{align*}
Since $\langle -A\nabla_xX_s^x z,\nabla_xX_s^x z\rangle_H=|(-A)^{1/2}\nabla_x X_s^xz|_H^2$ for any $s\in[0,t]$, from \eqref{estimate_der_H} we deduce that $\nabla_x X_\tau^xz\in D((-A)^{1/2})$ for any $\tau\in [0,T]$. Thus  \eqref{weight_estimate_derivative} holds for $\varepsilon=1/2$,  $t=0$ and any $z \in E$.

\noindent Let us now consider $\varepsilon\in[0,1/2)$. From  interpolation estimates (see e.g. \cite[Section 2.2]{LU95})
\begin{align}\label{inter_est}
|(-A)^\varepsilon e^{tA}x|_H\leq C_\varepsilon|x|^{1-2\varepsilon}_H|(-A)^{1/2}x|^{2\varepsilon}_H, \quad x\in D((-A)^{1/2}).
\end{align}
By replacing  $x$ by $\nabla_x X_s^x z$ in \eqref{inter_est},  we get
\begin{align*}
\int_0^\tau|(-A)^\varepsilon \nabla_x X_s^x z|^2_H ds
\leq & C_\varepsilon^2 \int_0^\tau |(-A)^{1/2} \nabla_x X_s^x z|^{4\varepsilon}_H| \nabla_x X_s^x z|^{2-4\varepsilon}_H ds \\
\leq & C_\varepsilon^2 \Big(\int_0^\tau |(-A)^{1/2} \nabla_x X_s^x z|^2_H ds\Big)^{2\varepsilon}\Big(\int_0^\tau| \nabla_x X_s^x z|^{2}_H ds\Big)^{1-2\varepsilon} 
\leq C \tau^{1-2\varepsilon}|z|_H^2,
\end{align*}
with $C:=\sup_{\varepsilon\in(0,1/2)}C_\varepsilon^2$. We conclude that  \eqref{weight_estimate_derivative} holds for $t=0$, $\varepsilon\in [0,\,1/2]$ and any $z\in E$.

Let us now  consider $z\in H$, and let $(z_n)\subset E$ be an approximating sequence of $z$ in $H$. Then, from \eqref{estimate_der_H}, for any $\tau\in[0,T]$ we get 
\begin{equation}\label{convX}
\nabla_x X^{x}_\tau z_n\rightarrow \nabla_x X^x_\tau z \quad  \P\text{-a.s. in } H, \text{ as } n\rightarrow+\infty.
\end{equation}
Since \eqref{weight_estimate_derivative} holds for any $z\in E$, it follows that $((-A)^{1/2}\nabla_xX^x z_n)$ is a Cauchy sequence in $\calm^2([0,T];H)$, and therefore there exists a process $\xi\in \calm^2([0,T];H)$ such that $(-A)^{1/2}\nabla_xX^x z_n\rightarrow \xi$ in  $\calm^2([0,T];H)$. Since $(-A)^{-1/2}$ is a bounded operator on $H$, it follows that 
\begin{align*}
\nabla_x X^x z_n=(-A)^{-1/2}(-A)^{1/2}\nabla_xX^x z_n\rightarrow (-A)^{-1/2}\xi,
\end{align*}
in  $\calm^2([0,T];H)$. Therefore, also by (\ref{convX}), $(-A)^{-1/2}\xi=\nabla_x X^xz$ a.e. in $(0,T)$ and $\P$-a.s.,
 which means that $\nabla_x X^xz\in D((-A)^{1/2})$ a.e. in $(0,T)$ and $\P$-a.s.,
   and $(-A)^{1/2}\nabla_x X^xz=\xi$ a.e. in $(0,T)$ and $\P$-a.s.
     In particular, we get
$$
\int_0^t|(-A)^{1/2}\nabla_x X_s^{t,x}z|^2_H ds \leq C|z|^2_H.
$$
Again, by applying interpolation estimates we see that \eqref{weight_estimate_derivative} holds for  $\varepsilon\in[0,1/2]$,  $t=0$ and any $z\in H$.
\endproof

We end this section by  giving pointwise estimates of $(-A)^\alpha\nabla_x X_\tau^{t,x}z$. In particular, we improve the result of Proposition \ref{P1.9}, by obtaining that $\nabla_xX_\tau^{t,x}$ belongs to $D((-A)^\alpha)$ for any $\tau\in[t,T]$, $\P$-a.s.
\begin{proposition}\label{P:pointwise_est} Let Hypothesis \ref{ip_forward} holds true and let $x\in E$,  $z\in H$ and  $t\in[0,T]$. Then, for any $x\in E$ and  $z\in H$, 
\begin{align}
\label{sper_est_der}
\E \Big[\sup_{\tau \in [t,\,T]} |(-A)^\alpha \nabla_x X_\tau^{t,x}z|_{H}\Big]
\leq & C|z|_H \left((\tau-t)^{-\alpha}+(\tau-t)^{1-\alpha}\left(|x|_E^{2m+1}+C_T\right)\right), 
\end{align}
and if in addition $z\in D((-A)^\alpha)$, then \eqref{pointwise_est_der2} gives 
\begin{align}
\label{sper_est_der2}
\E \Big[\sup_{\tau \in [t,\,T]}|(-A)^\alpha \nabla_x X_\tau^{t,x}z|_{H}\Big]
\leq & C\left(|(-A)^\alpha z|_H +(\tau-t)^{1-\alpha}\left(|x|_E^{2m+1}+C_T\right)|z|_H \right),
\end{align}	
where
$C_T:=\E[\sup_{\tau\in[0,T]}|w^A(\tau)|_E^{2m+1}]$.
\end{proposition}
\proof
We prove estimate \eqref{pointwise_est_der}, then  \eqref{pointwise_est_der2} follows from analogous arguments.
 Fix $x\in E$, $z\in H$ and let us consider $t=0$. 
We recall that $A$ generates an analytic semigroup on $H$ and therefore $e^{t A}x$ belongs to $D((-A)^k)$ for any $k\in \N$ and any $h\in H$, and $|(-A)^\beta e^{tA}h|_H\leq C_\beta t^{-\beta}|h|_H$ for any $\beta\geq0$ and some positive constant $C_\beta$. This means that $\nabla_x X_\tau^{x}z\in D((-A)^\alpha)$ for any $\tau\in[0,T]$ and, recalling \eqref{eq_mild_der_H}, 
\begin{align*}
(-A)^\alpha \nabla_x X_\tau^{x}z=(-A)^\alpha e^{\tau A}z+\int_0^\tau (-A)^\alpha e^{(\tau-s)A}\nabla_xF(X_s^{x})\nabla_x X_s^{x}z dx, \quad \P\textup{-a.s.}
\end{align*}
From Hypothesis \ref{forward}-4. and \eqref{estimate_der_H} we deduce that
\begin{align*}
|(-A)^\alpha \nabla_x X_\tau^{x}z|_H
\leq & C_\alpha \tau^{-\alpha}|z|_H+c\,C_\alpha \int_0^\tau (\tau-s)^{-\alpha}|\nabla F(X_s^x)\nabla_x X_s^x z|_H ds \\
\leq & C|z|_H \Big(\tau^{-\alpha}+\tau^{1-\alpha}\Big(|x|_E^{2m+1}+\sup_{\tau\in[0,T]}|w^A(\tau)|_E^{2m+1}\Big)\Big),
\quad \P\textup{-a.s.},
\end{align*}
for some positive constant $C$ independent of $x,z$.
Then, for any $\tau \in (t,T]$, 
\begin{align}
\label{pointwise_est_der}
|(-A)^\alpha \nabla_x X_\tau^{t,x}z|_{H}
\leq & C|z|_H \Big((\tau-t)^{-\alpha}+(\tau-t)^{1-\alpha}\Big(|x|_E^{2m+1}+\sup_{\tau\in[t,T]}|w^A(\tau)|_E^{2m+1}\Big)\Big),\,\, \P\textup{-a.s.},
\end{align}
for some positive constant $C$ independent of $x,z,t$. Further, if $z\in D((-A)^\alpha)$, then
\begin{align}
\label{pointwise_est_der2}
|(-A)^\alpha \nabla_x X_\tau^{t,x}z|_{H}
\leq & C\Big(|(-A)^\alpha z|_H +(\tau-t)^{1-\alpha}\Big(|x|_E^{2m+1}+\sup_{\tau\in[t,T]}|w^A(\tau)|_E^{2m+1}\Big)|z|_H \Big),
\,\, \P\textup{-a.s.}
\end{align}
Taking the expectation in \eqref{pointwise_est_der} and \eqref{pointwise_est_der2} we get respectively  (\ref{sper_est_der}) and (\ref{sper_est_der2}).
\endproof

\section{The forward-backward system}\label{Sec_forwbac}
 We consider the following 
forward-backward system of stochastic differential equations (FBSDE for short) for the unknown $(X,Y,Z)$ (also denoted by
$(X^{t,x},Y^{t,x},Z^{t,x})$ to stress the dependence on the initial
conditions $t$ and $x$): for given $t\in [0,T]$ and $x\in E$,
\begin{equation}\label{fbsde}
    \left\{\begin{array}{ll}
\dis dX_\tau =
AX_\tau d\tau+ F(X_\tau)d\tau+(-A)^{-\alpha}dW_\tau,& \tau\in
[t,T],
\\\dis
X_t=x,
\\\dis
 dY_\tau=-\psi(\tau,X_\tau,Y_\tau,Z_\tau)\;d\tau+Z_\tau\;dW_\tau, & \tau\in[t,T],
  \\\dis
  Y_T=\phi(X_T).
\end{array}\right.
\end{equation}
The second equation is of backward type for the unknown $(Y,Z)$
and depends on the Markov process $X$.
Under suitable assumptions on the coefficients $\psi$
(the so-called generator of the BSDE)
and $\phi$ 
we  look for a solution consisting of a pair of  processes $(Y,Z) \in \cals^2([t,T])\times\calm^2([t,T];H)$.
More precisely,   we will assume that $\psi$ is Lipschitz continuous with respect to $y$ and locally Lipschitz continuous and with quadratic growth with respect to $z$, as stated below.
\begin{hypothesis}
\label{ip-psiphi} The functions  $\mathbb{\phi}:E\rightarrow\mathbb{R}$  and $\psi:[0,T]\times E\times \R\times H
\rightarrow\mathbb{R}$  in \eqref{fbsde} satisfy the following.
\begin{itemize}
\item[(i)] $\mathbb{\phi}$ is continuous,  and there exists a nonnegative constant $K_\phi$ such that 
$\vert \phi(x)\vert\leq K_\phi$
for every  $x\in E$.
\item[(ii)]	$\psi$ 
is measurable and,  for every fixed $t\in[0,T]$, the map $\psi(t,\cdot,\cdot,\cdot):E\times\R\times H\rightarrow\R$
is continuous.
 Moreover, there exist nonnegative constants $L_\psi$ and $K_\psi$
such that
\begin{align*}
&\vert \psi(t,x_1,y_1,z_1)- \psi(t,x_2,y_2,z_2)\vert\leq L_\psi\left(
\vert x_1-x_2\vert_E+
\vert y_1-y_2\vert+\vert z_1-z_2\vert_{H}
(1+\vert z_1\vert_{H}+\vert z_2\vert_{H})\right),\\
&\vert \psi(t,x,0,0)\vert\leq K_\psi,
\end{align*}
for every $t\in[0,T]$, 
$x_1, x_2 \in E$, 
$y_1,y_2\in\R$ and $z_1,z_2\in H$.
\end{itemize}
\end{hypothesis}

\begin{theorem}\label{teo-ex-bsdequadr} 
Assume that Hypotheses \ref{ip_forward} and  \ref{ip-psiphi} hold
true, and for any   $(t, x) \in [0, T ] \times E$, let $(X^{t,x},Y ^{t,x}, Z^{ t,x})$ be a solution to the FBSDE \eqref{fbsde}.
Then,  there exists a unique solution of the Markovian BSDE in \eqref{fbsde} such that 
\[
 \Vert Y^{t,x}\Vert _{\cals^2([t,T])}+\Vert Z^{t,x}\Vert_{\calm^2([t,T];H)}\leq C,
\]
where $C$ is a constant that may depend on $T,\, A,\,F,
\, K_\psi,\,L_\psi,\,K_\phi.$.
Moreover,  setting $v(t,x):= Y^{t,x}_t$, 
\begin{equation}\label{Id_Y_Markov_prop}
Y^{t,x}_s = v(s,X^{t,x}_s), \quad \P\textup{-a.s.},  \,\,s \in [t,\,T],
\end{equation}
and there  exists a Borel  function  $u: [t,\,T] \times E \rightarrow H$ such that 
\begin{equation}\label{Id_Z_Markov_prop}
 Z^{t,x}_s = u(s,X^{t,x}_s), \quad \P\textup{-a.s.},\,\, \textup{a.e.} \,\,s \in [t,\,T].
\end{equation} 
\end{theorem}
\proof The first part of the  result substantially follows from \cite{Kob}. 
Identities \eqref{Id_Y_Markov_prop}-\eqref{Id_Z_Markov_prop} are a consequence of the Markov property  of $X$,
see for instance Theorem 4.1 in \cite{ElKaPengQuenez} 
or  the proof of Theorem 5.1 in \cite{FT05}.
\endproof

We recall some further estimates for the solution $(Y,Z)$ of the forward-backward system (\ref{fbsde}).
In particular, 
$Z\in\calm^p([t,T];H)$, for any $p\geq 1$.
The corresponding proof can be found e.g. in \cite{Mas-spa}. 
\begin{proposition}\label{teo-stimap-bsdequadr}
Assume that Hypotheses \ref{ip_forward} and  
 \ref{ip-psiphi}
hold true,
and for any   $(t, x) \in [0, T ] \times E$, let $(X^{t,x},Y ^{t,x}, Z^{ t,x})$ be a solution to the FBSDE \eqref{fbsde}.
Then, for all $p\geq 1$, 
\[
 \Vert Y^{t,x}\Vert _{\cals^p([t,T])}+\Vert Z^{t,x}\Vert_{\calm^p([t,T];H)}\leq C,
\]
where $C$ is a constant that may depend on $T, \,A,\,F,\, K_\psi,\,L_\psi,\,K_\phi$.
\end{proposition}
At this point, we aim at  proving a stability result for the BSDE when the final datum and the generator are approximated by sequences of Fr\'echet differentiable functions $(\phi_n)_{n\geq 1},\, (\psi_{\ell})_{\ell\geq 1}$, converging pointwise respectively to $\phi$ and $\psi$, and such that, for all $t\in[0,T]$, $x,x_1,x_2\in E$, $y_1, y_2 \in \R$, $z_1, z_2 \in H$, 
\begin{align}
&\vert \phi_n(x)\vert\leq K_\phi, \quad
\vert \psi_\ell(t,x,0,0)\vert\leq K_\psi,
\label{stime_appr_1}\\
&|\psi_{\ell}(t,x_1,y_1,z_1)-\psi_{\ell}(t,x_2,y_2,z_2)|
\leq L_\psi(|x_1-x_2|_E+|y_1-y_2|+|z_1-z_2|_H(1+|z_1|_H+|z_2|_H)).\label{stime_appr_2}
\end{align}
To provide such approximations we extend the result in {\cite{PZ} valid for Hilbert
spaces: by using Schauder basis, the approximation performed in that paper can
be achieved also in Banach spaces, along the lines of what is done in \cite{MasBanEJP}. We start by introducing the following objects.
\begin{definition}
\label{appr_psil_phin}
\begin{itemize}
	\item [i)] Denote by $(e_n)_{n\geq 1}$ the normalized Schauder basis in $E$ and by $(h_n)$ an orthonormal basis of $H$.
For any $n\in\N$, we define the projections $Q_n:H\rightarrow \R^n$ and $P_n:E\rightarrow  \R^n$ as follows: 
$$Q_nz:=(z_1,\ldots,z_n), \quad P_nx:=(x_1,\ldots,x_n),
$$
 for any $z\in H$ and $x\in E$ with
$z=\sum_{n=1}^\infty z_nh_n$ and $x=\sum_{n=1}^\infty x_ne_n$, $z_n,x_n\in \R$.
\item [ii)] We consider  nonnegative smooth kernels $\vartheta\in C^\infty_c(\R)$ and $\rho_m\in C^\infty_c(\R^m)$,  $m\in\N$, such that
$${\rm  supp}\ \!(\vartheta)\subseteq \{\zeta\in\R:|\zeta|\leq 1\},\quad {\rm  supp}\ \!(\rho_m)\subseteq \{\xi\in\R^m:|\xi|\leq m^{-1}\},\quad \|\vartheta\|_{L^1(\R)}=\|\rho_m\|_{L^1(\R^m)}=1.$$
\item [iii)]  For any $n,\ell\in\N$,  we set  $\vartheta_\ell(\zeta)=\ell\vartheta(\ell\zeta)$ for any $\zeta\in\R$, and
\begin{align}
&\phi_n(x)=\int_{\R^n}\rho_n(\xi-P_n x)\phi\Big(\sum_{i=1}^n\xi_ie_i\Big)\! d\xi,\label{phi-approx}\\
&\psi_{\ell}(t,x,y,z):=\int_{\R^\ell}\int_{\R^\ell}\int_{\R}\rho_\ell\left(\xi-P_\ell x\right)\rho_\ell(\eta-Q_\ell z)\vartheta_\ell(y-\zeta)\psi\Big(t,\sum_{i=1}^\ell\xi_ie_i,\zeta,\sum_{j=1}^\ell \eta_jh_j\Big)d\zeta \ \!d\eta\ \! d\xi.\label{psi-approx}
\end{align}
\end{itemize} 
\end{definition}
It is not hard to prove the following lemma.
\begin{lemma}\label{L:3.4} Le $\phi$ and $\psi$ satisfy Hypothesis \ref{ip-psiphi}. Then the following hold.
\begin{itemize}
\item [(i)] For any $n\in\N$, the function $\phi_n$ in \eqref{phi-approx} is Fr\'echet differentiable,  satisfies estimate \eqref{stime_appr_1},  and
\begin{align*}
\lim_{n\rightarrow+\infty}\phi_n(x)=\phi(x), \quad x\in E.
\end{align*}
\item [(ii)] For any $\ell\in\N$, the function $\psi_{\ell}$ in \eqref{psi-approx} is Fr\'echet differentiable with respect to $x,y,z$, satisfies estimates \eqref{stime_appr_1}-\eqref{stime_appr_2}, and 
\begin{align*}
\lim_{\ell\rightarrow+\infty}\psi_\ell(t,x,y,z)=\psi(t,x,y,z), \quad (t,x,y,z)\in [0,T]\times E\times \R\times H.
\end{align*}
\end{itemize}
\end{lemma}

We can now  give a stability result for the Markovian BSDE in \eqref{fbsde} related to a forward process  $X$ taking values in the Banach space $E$,  when the final datum and the generator are approximated respectively by the sequences  $(\phi_n)_{n\geq 1}$ and  $(\psi_{\ell})_{\ell\geq 1}$.
Notice that a similar result is proved in \cite{Mas-spa}, where the forward process $X$ takes its values in a Hilbert space $H$: there the final datum and the generator are approximated in the norm of the uniform convergence by means of their inf-sup convolutions. 
Clearly, 
the following result holds true if we approximate only $\psi$ or $\phi$. 
\begin{proposition}
\label{prop-convp-bsdequadr}
Assume that Hypotheses \ref{ip_forward} and \ref{ip-psiphi} hold
true. For any   $(t, x) \in [0, T ] \times E$, let $(X,Y, Z)$ be a solution to the FBSDE \eqref{fbsde}.
Let $(Y^{n,l},Z^{n,l})$ be the solution of the BSDE in the forward-backward system
\begin{equation}\label{fbsde-n}
    \left\{\begin{array}{ll}
    \dis dX_\tau =
AX_\tau d\tau+ F( X_\tau)d\tau+(-A)^{-\alpha}dW_\tau,& \tau\in
[t,T],
\\\dis
X_t=x,
\\\dis
 dY^{n,l}_\tau=-\psi_l(\tau,X_\tau,Y^{n,l}_\tau,Z^{n,l}_\tau)\;d\tau+Z^{n,l}_\tau\;dW_\tau, & \tau\in[t,T], 
  \\\dis
  Y^{n,l}_T=\phi_n(X_T),
\end{array}\right.
\end{equation}
that is, the FBSDE \eqref{fbsde}
with final datum equal to $\phi_n$  in \eqref{phi-approx} in  place of $\phi$,  and with generator $\psi_l$  in \eqref{psi-approx} in place of $\psi$.
Then, for all $p\geq 1$, the unique solution of the Markovian BSDE in \eqref{fbsde} is such that
\[
 \Vert Y-Y^{n,l}\Vert _{\cals^p([t,T])}+\Vert Z-Z^{n,l}\Vert_{\calm^p([t,T];H)}\rightarrow 0 \qquad \text{ as }n,l\rightarrow \infty.
\]
\end{proposition}
\proof
Thanks to \eqref{stime_appr_1}, \eqref{stime_appr_2} and to Proposition \ref{teo-stimap-bsdequadr}, the pair of processes $(Y^{n,l},Z^{n,l})$ is bounded in $\cals^p([t,T])\times\calm^p([t,T];H)$, uniformly with respect to $n,l$.
The BSDE satisfied by the pair of the difference processes $(Y^{n,l}-Y,Z^{n,l}-Z)$ is
\begin{equation*}
    \left\{\begin{array}{l}
 d(Y^{n,l}_\tau-Y_\tau)=\left(\psi(\tau,X_\tau,Y_\tau,Z_\tau)-\psi_l(\tau,X_\tau,Y^{n,l}_\tau,Z^{n,l}_\tau)\right)d\tau+(Z^{n,l}_\tau-Z_\tau)\;dW_\tau, \quad \tau\in[t,T], 
  \\
  Y^{n,l}_T-Y_T=\phi_n(X_T)-\phi(X_T).
\end{array}\right.
\end{equation*}
Writing the previous equation in the integral form, we get
\begin{align*}
&Y^{n,l}_\tau-Y_\tau \\
&= \phi_n(X_T)-\phi(X_T)  -\int_\tau^T(Z^{n,l}_s-Z_s)\;dW_s
+\int_\tau^T\left(\psi(s,X_s,Y_s,Z_s)-\psi_l(s,X_s,Y_s,Z_s)\right)\;ds \\
&+\int_\tau^T\left(\psi_l(s,X_s,Y_s,Z_s)-\psi_l(s,X_s,Y_s,Z^{n,l}_s)\right)\;ds\\
&+\int_\tau^T\left(\psi_l(s,X_s,Y_s,Z^{n,l}_s)-\psi_l(s,X_s,Y^{n,l}_s,Z^{n,l}_s)\right)\;ds \\
=& \phi_n(X_T)-\phi(X_T)+\int_\tau^T\left(\psi(s,X_s,Y_s,Z_s)-\psi_l(s,X_s,Y_s,Z_s)\right)\;ds\\
&+\int_\tau^T\frac{\psi_l(s,X_s,Y_s,Z_s)-\psi_l(s,X_s,Y_s,Z^{n,l}_s)}{Z_s-Z^{n,l}_s}\left(
Z_s-Z^{n,l}_s\right)\;ds  -\int_\tau^T(Z^{n,l}_s-Z_s)\;dW_s\\
& +\int_\tau^T\frac{\psi_l(s,X_s,Y_s,Z^{n,l}_s)-\psi_l(s,X_s,Y_s^{n,l},Z^{n,l}_s)}{Y_s-Y^{n,l}_s}\left(Y_s-Y^{n,l}_s\right)\;ds \\
=& \phi_n(X_T)-\phi(X_T) +\int_\tau^T\left(\psi(s,X_s,Y_s,Z_s)-\psi_l(s,X_\tau,Y_s,Z_s)\right)\;ds
-\int_\tau^T(Z^{n,l}_s-Z_s)\;dW^{n,l}_s \\
&
+\int_\tau^T\frac{\psi_l(s,X_s,Y_s,Z^{n,l}_s)-\psi_l(s,X_s,Y^{n,l}_s,Z^{n,l}_s)}{Y_s-Y^{n,l}_s}\left(Y_s-Y^{n,l}_s\right)\;ds,
\end{align*}
where in the last passage we have used that 
\begin{align*}
W^{n,l}_\tau=W_\tau-\int_t^\tau\frac{\psi_l(s,X_s,Y_s,Z_s)-\psi_l(s,X_s,Y_s,Z^{n,l}_s)}{Z_s-Z^{n,l}_s}ds, \quad \tau\geq t,
\end{align*}
which, by the Girsanov Theorem (see, e.g., \cite[Theorem~10.14]{dpz92}), is a cylindrical Wiener process  under an equivalent probability measure $Q^{n,l}$.
Taking the $Q^{n,l}$-conditional expectation  $\E_{Q^{n,l}}^{\calf_\tau}[\cdot]:= \E_{Q^{n,l}}[\cdot| \calf_\tau]$,  we get 
\begin{align*}
 Y^{n,l}_\tau-Y_\tau 
=& \E_{Q^{n,l}}^{\calf_\tau}[
\phi_n(X_T)-\phi(X_T)]
+\E_{Q^{n,l}}^{\calf_\tau}\Big[\int_\tau^T\left(\psi(s,X_s,Y_s,Z_s)-\psi_l(s,X_s,Y_s,Z_s)\right)\;ds\Big]\\
&+\E_{Q^{n,l}}^{\calf_\tau}\Big[\int_\tau^T\frac{\psi_l(s,X_s,Y_s,Z^{n,l}_s)-\psi_l(s,X_s,Y^{n,l}_s,Z^{n,l}_s)}{Y_s-Y^{n,l}_s}\left(Y_s-Y^{n,l}_s\right)\;ds\Big].
\end{align*}
By taking the absolute value, the expectation and by applying the Gronwall lemma, we deduce that, for all $p \geq 1$, $Y^{n,l}\rightarrow Y$ in $\cals^p([t,T])$ as $n,\ell\rightarrow\infty$, with respect to the probability measure $Q^{n,l}$ and also with respect to the original probability measure. 

For what concerns the estimate of $Z-Z^{n,l}$, by applying  the It\^o formula to $|Y^{n,l}-Y|^2$ we get
\begin{align*}
&\E\Big[|Y^{n,l}_t-Y_t|^2\Big]+\E\Big[\int_t^T\vert Z^{n,l}_\tau-Z_\tau\vert^2_H\,d\tau\Big]\\
&=\E\Big[|\phi_n(X_T)-\phi(X_T)|^2\Big] -2\E\Big[\int_t^T\left(Y^{n,l}_\tau-Y_\tau\right)\left(\psi(\tau,X_\tau,Y_\tau,Z_\tau)-\psi_l(\tau,X_\tau,Y^{n,l}_\tau,Z^{n,l}_\tau)\right)d\tau\Big]\\&
\leq \E\Big[|\phi_n(X_T)-\phi(X_T)|^2\Big]+2\E\Big[\int_t^T\vert Y^{n,l}_\tau-Y_\tau\vert\vert\psi(\tau,X_\tau,Y_\tau,Z_\tau)-\psi_l(\tau,X_\tau,Y_\tau,Z_\tau)\vert\,d\tau\Big]\\ 
&+2\E\Big[\int_t^T\vert Y^{n,l}_\tau-Y_\tau\vert\vert\psi_l(\tau,X_\tau,Y_\tau,Z_\tau)-\psi_l(\tau,X_\tau,Y^{n,l}_\tau,Z^{n,l}_\tau)\vert\,d\tau\Big]\\
&\leq \E\Big[|\phi_n(X_T)-\phi(X_T)|^2\Big]+2\E\Big[\sup_{\tau\in[t,T]}\vert Y^{n,l}_\tau-Y_\tau\vert\int_t^T\vert\psi(\tau,X_\tau,Y_\tau,Z_\tau)-\psi_l(\tau,X_\tau,Y_\tau,Z_\tau)\vert\,d\tau\Big]\\ &
+C\E\Big[\sup_{\tau\in[t,T]}\vert Y^{n,l}_\tau-Y_\tau\vert\int_t^T\left(1+\vert Y_\tau-Y^{n,l}_\tau\vert+\vert Z_\tau-Z^{n,l}_\tau \vert_{H}\left(1+\vert Z_\tau\vert_{H}+\vert Z^{n,l}_\tau \vert_{H}
\right)\right)d\tau\Big].
\end{align*}
Let us consider the right-hand side of the  above inequality.
Thanks to   estimates (\ref{stime_appr_1}), (\ref{stime_appr_2}) and  to the boundedness of $Y$ and $Y^{n,l}$ in $\cals^p([t,T];E)$,  the first two terms converge to $0$ as $n,l\rightarrow \infty$ by the dominated convergence theorem. For what concerns the third term, by applying H\^older's inequality with $p,q$ conjugate exponents, we get
\begin{align*}
&\E\Big[\sup_{\tau\in[t,T]}\vert Y^{n,l}_\tau-Y_\tau\vert\int_t^T\left(1+\vert Y_\tau-Y^{n,l}_\tau\vert+\vert Z_\tau-Z^{n,l}_\tau \vert_{H}\left(1+\vert Z_\tau\vert_{H}+\vert Z^{n,l}_\tau \vert_{H}\right)\right) d\tau\Big]\\
&\leq C\Big(\E\Big[\sup_{\tau\in[t,T]}\vert Y^{n,l}_\tau-Y_\tau\vert^p\Big]\Big)^{\frac{1}{p}}\Big( \E\Big[\Big(\int_t^T(1+\vert Z_\tau\vert_{H}^2+\vert Z^{n,l}_\tau \vert_{H}^2)
 d\tau\Big)^q\Big]\Big)^\frac{1}{q}\rightarrow 0
\end{align*}
as $n,l\rightarrow \infty$.
 The stability result for $p=2$ follows, and   we can pass  to the case of general $p$ in a usual way.
\endproof

We now state a result on differentiability  for the solution of a Markovian BSDE with generator with quadratic growth, with respect to the initial datum $x$.
\begin{proposition}\label{P:BC_Diff_BSDE}
Assume that Hypotheses \ref{ip_forward} and \ref{ip-psiphi} hold
true, and for any   $(t, x) \in [0, T ] \times E$, let $(X^{t,x},Y ^{t,x}, Z^{ t,x})$ be a solution to the FBSDE \eqref{fbsde}.
Assume moreover  that $\phi$ is G\^ateaux differentiable with  bounded derivative, and that $\psi$
is G\^ateaux differentiable with respect to $x$, $y$ and $z$. Then the triple of processes $(X^{t,x},Y^{t,x},Z^{t,x})$ is  G\^ateaux differentiable as a map from $E$ with values in $
\cals^2((t,T];E)\times \cals^2([t,T])\times \calm^2([t,T];H)$ and, for any $h\in {E}$,  
\begin{align}\label{eq_gradYgradZ}
\left\{
\begin{array}{rcl}
-d\nabla_x Y_\tau^{t,x} h&=& \nabla_x \psi(\tau,X_\tau^{t,x},Y_\tau^{t,x},Z^{t,x}_\tau)\nabla_x X^{t,x}_\tau h \, d\tau+\nabla_y \psi(\tau,X_\tau^{t,x},Y_\tau^{t,x},Z^{t,x}_\tau)\nabla_x Y^{t,x}_\tau h \,d\tau\vspace{2mm} \\
& & +\nabla_z \psi(\tau,X_\tau^{t,x},Y_\tau^{t,x},Z^{t,x}_\tau)\nabla_x Z^{t,x}_\tau h \, d\tau-\nabla_x Z^{t,x}_\tau h \,dW_\tau, \quad \tau\in[t,T], \vspace{2mm} \\
\nabla_x Y_T^{t,x} h & =& \nabla_x \phi(X_T^{t,x})\nabla_x X_T^{t,x}h, \vspace{2mm} \\
d\nabla_x X_\tau^{t,x}h & = & A \nabla_x X_\tau^{t,x}h d\tau+\nabla F(X_\tau^{t,x})\nabla_x X_\tau^{t,x}h \, d\tau, \quad \tau\in[t,T], \vspace{1mm} \\
\nabla_x X_t^{t,x}h & = & h.
\end{array}
\right.
\end{align}
Moreover, there exists
a constant $C$, only dependent on $T,\,A,\,F,\,K_\psi,\,L_\psi,\,K_\phi$, such that 
\begin{align}
&\E\Big[\sup_{\tau\in [t,\,T]}|\nabla Y_\tau^{t,x}h|^2 + \int_{t}^T |\nabla_x Z_\tau^{t,x}h|_{H}^2\, d\tau \Big] \leq C |h|_{H}^2.\label{new_est_gradY_gradZ} 
\end{align}
\end{proposition}
\proof
In  the case of a Markovian BSDE with generator $\psi$ quadratic with respect to $Z$ and related to a forward process taking values in a Hilbert space, the result is given  in Theorem 4.5 of \cite{BriFu}.  Since in Proposition \ref{P:unificata} we have proved the differentiability of $X^{t,x}$ with respect to $x\in E$,  the same  conclusions hold when the forward process takes  values in the Banach space $E$, 
namely 
\begin{align*}
&\E\Big[\sup_{\tau\in [t,\,T]}|\nabla_x Y_s^{t,x}h|^2 + \int_{t}^T |\nabla_x Z_\tau^{t,x}h|_{H}^2\, d\tau \Big] \leq C |h|_{E}^2.
\end{align*} 
The stronger estimate  \eqref{new_est_gradY_gradZ} comes from    Proposition \ref{P:unificata}, estimate (\ref{estimate_der_H}).
\endproof

\subsection{Identification of \texorpdfstring{$Z$}{Z} and a priori estimates on \texorpdfstring{$(Y,Z)$}{(Y,Z)}}
We now prove an a priori estimate on $Z^{t,x}$ depending only on the $L^{\infty}$-norm of the final datum. The novelty towards \cite{Mas-spa} is that we work in a Banach space and
the pseudo-inverse of the diffusion operator
is the unbounded operator $(-A)^{\alpha}$. In order to get this estimate and also for the subsequent results   of the paper, it will be crucial to prove the identification
\[
  Z_t^{t,x} =\nabla_x Y_t^{t,x}(-A)^{-\alpha},
\]
which is new in the Banach space framework and in the case of quadratic generator  with respect to $z$. 
We have to make the following assumption:
\begin{hypothesis}\label{ipG-agg}
There exists a Banach space $E_0\subset D((-A)^\alpha)$ dense in $H$ such that $(-A)^{-\alpha}E_0 \subset E$ and $(-A)^{-\alpha}:E_0\rightarrow E$ is continuous. 
\end{hypothesis}
\begin{remark}
Notice that if $\mathcal D\subset\R^2$ is a bounded open domain with smooth boundary, $H=L^2(\mathcal D)$ and $A$ is the Laplace operator in dimension $2$ with Dirichlet boundary conditions, then we can take $E_0=D((-A)^{\frac{1}{2}})$ and all the requirements of Hypothesis \ref{ipG-agg} are verified.
\end{remark}

\begin{theorem}
\label{prop-identific-Z}  
Assume that Hypotheses \ref{ip_forward}, \ref{ip-psiphi} and \ref{ipG-agg} hold true, that $\phi$ is G\^ateaux differentiable with  bounded derivative, and that $\psi$
is G\^ateaux differentiable with respect to $x$, $y$ and $z$. For any   $(t, x) \in [0, T ] \times E$, let $(X^{t,x},Y ^{t,x} , Z^{ t,x} )$ be the solution to the FBSDE \eqref{fbsde}. 
Then the triple of processes $(X^{t,x},Y^{t,x},Z^{t,x})$ is G\^ateaux differentiable as a map from $E$ with values in $
\cals^2((t,T];E)\times \cals^2([t,T])\times \calm^2([t,T];H)$. Moreover, setting 
$v(t, x) = Y_t^{t,x}$, then,  $\P$-a.s.,
\begin{align}
Y_s^{t,x} &= v(s, X_s^{t,x} ),\quad s\in[t,T], \label{Id_Y}\\
Z_s^{t,x} h &=\nabla_x v(s, X_s^{t,x} )\nabla_x  X_s^{t,x}(-A)^{-\alpha}h,  \quad  \textup{a.e.} \ s\in[t,T], \ h \in E_0.\label{Id_Z}
\end{align}
\end{theorem}
\proof
The differentiability properties of $(X^{t,x},Y^{t,x},Z^{t,x})$  and the  identification formula  \eqref{Id_Y} directly   follow respectively  from Proposition \ref{P:BC_Diff_BSDE} and  formula \eqref{Id_Y_Markov_prop} in Theorem  \ref{teo-ex-bsdequadr}.

Let us now prove  identification formula (\ref{Id_Z}) for $Z$. Since we are in a  Banach space framework, we will follow the lines of the proof of Theorem 3.17 in \cite{MasBan}. However, a substantial difference with respect to \cite{MasBan} is that  
here we deal with  a generator $\psi$ with quadratic growth with respect to $z$, instead of   Lipschitz continuous. 
Fix $t \in [0,\,T]$. By  the definition of the function $v$,  we can write 
\begin{equation}\label{repr_v-forward}
v\left(\tau,X^{t,x}_\tau\right) +\int_{t}^{\tau}Z_{\sigma}^{t,x}dW_{\sigma}=v(t,x)+\int_{t}^{\tau}\psi_\sigma d\sigma, \quad 0\leq t \leq \tau \leq T, 
\end{equation}
where we have used the notation $\psi_\sigma:=\psi(\sigma, X_\sigma^{t,x}, Y_\sigma^{t,x},Z_\sigma^{t,x})$.
Notice that towards \cite{MasBan} we do not have $\psi\in \calm^2([0,T])$, but we only know that $\psi\in L^p(\Omega ,L^1(0,T;\R))$ for any $p\geq 2$.
As in \cite{MasBan},  we define a
family $\mathcal{S}$ of predictable processes with real values in the following way:
{\small
\begin{align*}
\mathcal{S} = \Big\{&\textup{predictable processes } \eta: 
\textup{for any } \,\,k=0,\ldots,2^{n}-1,\\
&\eta_t 1_{\left[\frac{kT}{2^{n}},\frac{\left( k+1\right)
T}{2^{n}}\right)}(t) = \eta^{k}\left( W_{t_{1}},\ldots,W_{t_{l_{k}}}\right)\,\,\textup{for}\,\, 0\leq t_{1}\leq\cdots\leq
t_{l_{k}}\leq\dfrac{kT}{2^{n}},\\
& \eta^{k}\,\, \textup{bounded functions in}\,\, C^{\infty}(\mathbb{R}
^{l_{k}},\mathbb{R}) \textup{ with bounded derivatives of all orders}\Big\}.
\end{align*}}
We will briefly write $\eta_{t}=\eta_{t}(W_{\cdot}) $, where
by $W_{\cdot}$\ we mean the trajectory of $W$ up to time $t$. 

 Let us set  $\xi_{t}:=\eta_{t}\varsigma$ for 
$\varsigma\in E_0$. 
From now on we fix $s>t$, and  $\delta>0,$ small enough such that $s-\delta>t$.
We also identify $H$ with its dual $H^{\ast}$, and we
write $\xi$ for $\xi^{\ast}$. Multiplying both sides of
(\ref{repr_v-forward}), with $\tau$ replaced by $s$, by $ {\int_{s-\delta}^{s}}
\xi_{\sigma}dW_{\sigma}$  and  taking the expectation, we get%
\begin{equation}\label{first}
\mathbb{E}\Big[ v\left( s,X^{t,x}_{s}\right)
\int_{s-\delta}^{s}\xi_{\sigma }dW_{\sigma}\Big] 
=\mathbb{E}\Big[
\int_{t}^{s}\psi_{\sigma}d\sigma\int_{s-\delta
}^{s}\xi_{\sigma}dW_{\sigma}\Big] 
+\mathbb{E}\Big[
\int_{t}^{s}Z^{t,x}_{\sigma}dW_{\sigma}\int_{s-\delta}^{s}\xi_{\sigma}dW_{\sigma}\Big].
\end{equation}
It is immediate that
\[
\mathbb{E}\Big[{\int_{t}^{s-\delta}} \psi_{\sigma}d\sigma
{\int_{s-\delta}^{s}} \xi_{\sigma}dW_{\sigma}\Big] =0,\quad
\mathbb{E}\Big[ {\int_{t}^{s}} Z^{t,x}_{\sigma}dW_{\sigma}
{\int_{s-\delta}^{s}} \xi_{\sigma}dW_{\sigma}\Big]
=\mathbb{E}\Big[ {\int_{s-\delta}^{s}}
Z^{t,x}_{\sigma}\xi_{\sigma}d\sigma\Big], 
\]
so \eqref{first} simplifies in   
\begin{equation*}
\mathbb{E}\Big[ v\left( s,X^{t,x}_{s}\right)
\int_{s-\delta}^{s}\xi_{\sigma }dW_{\sigma}\Big] 
=\mathbb{E}\Big[
\int_{s-\delta}^{s}\psi_{\sigma}d\sigma\int_{s-\delta
}^{s}\xi_{\sigma}dW_{\sigma}\Big] 
+\mathbb{E}\Big[
\int_{s-\delta}^{s}Z^{t,x}_{\sigma}\, \xi_{\sigma}d \sigma\Big].
\end{equation*}
By dividing both sides of the previous equality by $\delta$ and   letting $\delta\rightarrow0$,
we get 
\begin{equation}\label{ultima}
\mathbb{E}\Big[ Z^{t,x}_{s}\xi_{s}\Big] =\lim_{\delta\rightarrow0}\frac
{1}{\delta}\mathbb{E}\Big[ v\left( s,X^{t,x}_{s}\right) \int_{s-\delta}^{s}%
\xi_{\sigma}dW_{\sigma}\Big]- \lim_{\delta\rightarrow0}\frac
{1}{\delta}\mathbb{E}\Big[
\int_{s-\delta}^{s}\psi_{\sigma}d\sigma\int_{s-\delta
}^{s}\xi_{\sigma}dW_{\sigma}\Big].
\end{equation}
We will prove that
\begin{align}
&\lim_{\delta\rightarrow0}\frac{1}{\delta}\mathbb{E}\Big[ v\left(
s,X^{t,x}_{s}\right) \int_{s-\delta}^{s}\xi_{\sigma}dW_{\sigma}\Big]
=\mathbb{E}\Big[ \nabla_{x} \big(v( s,X^{t,x}_{s})\big) (-A)^{-\alpha}\xi_{s}\Big],\label{toprove_1}\\
&\lim_{\delta\rightarrow0}\frac
{1}{\delta}\mathbb{E}\Big[
\int_{s-\delta}^{s}\psi_{\sigma}d\sigma\int_{s-\delta
}^{s}\xi_{\sigma}dW_{\sigma}\Big]=0. \label{toprove_2}
\end{align}
If \eqref{toprove_1} and \eqref{toprove_2} hold, then, by \eqref{ultima},   for every $\eta\in\mathcal{S}$, $\mathbb{E}[Z^{t,x}_{\sigma}\varsigma\eta
_{s}] =\mathbb{E}[\nabla v(\sigma,X^{t,x}_{\sigma}) (-A)^{-\alpha}\varsigma \eta_{s}]$ for
almost every $\sigma\in[t,T]$. By the arbitrariness of $\eta$,  we would have, for almost every $\sigma\in[t,T]$, $Z^{t,x}_{\sigma}\varsigma=\nabla_{x}
\big(v(\sigma,X^{t,x}_{\sigma})\big) (-A)^{-\alpha}\varsigma,$ $\mathbb{P}$-a.s. for
all $\varsigma\in E_0$, and the formula (\ref{Id_Z}) would follow. 

 Let us thus show that \eqref{toprove_1} and \eqref{toprove_2} hold true. 
We start by proving \eqref{toprove_1}. 
One  proceeds as in \cite{B}, following also \cite{MasBan}. In particular,  for $0\leq t\leq\sigma\leq T$, we define
%
\begin{equation}\label{W^epsilon}
W_{\sigma}^{\varepsilon}=W_{\sigma}-\varepsilon
{\displaystyle\int_{t}^{\sigma}} \xi_{r}\left(
W_{\cdot}^{\varepsilon}\right) {\it dr},
\end{equation}
where $\xi_{r}(W_{\cdot}^{\varepsilon})$ depends on the
trajectories of $W_{\cdot}^{\varepsilon}$ up to time $r$, and the
dependence is given by the definition  of
$\eta$. The process $(W_{\sigma}^{\varepsilon})_{\sigma}$ is
defined as the solution of (\ref{W^epsilon}), which is not
considered as a stochastic differential equation, as specified in
\cite[ p.~476]{B}. Equation~(\ref{W^epsilon}) can be solved step
by step in each interval
\[
\left[\frac{\it kT}{2^{n}}, \frac{\left( k+1\right)
T}{2^{n}}\right),\quad k=0,\ldots,2^{n}-1.
\]
$(W_{\sigma}^{\varepsilon}) _{\sigma}$ is
well defined for every $0\leq\sigma\leq T$, see \cite{MasBan} for more details. Moreover,
$W_{\sigma}^{\varepsilon }$ is a function of the trajectories of
$W$ up to time $\sigma$, that is,
$W_{\sigma}^{\varepsilon}=W_{\sigma}^{\varepsilon}(W_{\cdot})$, and 
we can write
\[
W_{\sigma}^{\varepsilon}=W_{\sigma}-\varepsilon%
{\int_{t}^{\sigma}} \xi_{r}\left( W_{\cdot}^{\varepsilon}\left(
W_{\cdot}\right) \right) {\it dr},\qquad 0\leq t\leq\sigma\leq T.
\]
Now we define a probability measure $Q_{\varepsilon}$ such that%
\[
\dfrac{dQ_{\varepsilon}}{d\mathbb{P}}=\exp\Big(\varepsilon%
{\displaystyle\int_{t}^{T}}
\xi_{\sigma}\left( W_{\cdot}^{\varepsilon}\left( W_{\cdot}\right) \right)
dW_{\sigma}-\dfrac{\varepsilon^{2}}{2}%
{\displaystyle\int_{t}^{T}}
\left\vert \xi_{\sigma}\left( W_{\cdot}^{\varepsilon}\left( W_{\cdot
}\right) \right) \right\vert ^{2}d\sigma\Big).
\]
By the Girsanov Theorem,
under $Q_{\varepsilon}$, $W_{\sigma}^{\varepsilon}=W_{\sigma}-\varepsilon%
{\int_{t}^{\sigma}} \xi_{r}(W_{\cdot}^{\varepsilon}(W_{\cdot}))
{\it dr}$ is 
a cylindrical
Wiener process in $H$. By this construction of
$(W_{\sigma}^{\varepsilon})_{\sigma}$, it is also clear that for
every $0\leq\sigma\leq T$, $W_{\sigma}^{\varepsilon}$ is pathwise
differentiable with respect to\ $\varepsilon$ and
$\frac{d}{d\varepsilon}_{\mid\varepsilon=0}W_{\sigma}^{\varepsilon}=-%
\displaystyle{\int_{t}^{\sigma}} \xi_{r}(W_{\cdot}){\it dr}$, see also
\cite[p.~476]{B}.

\noindent By (\ref{repr_v-forward}), the random varaible $v(s,X^{t,x}_{s})$ is square integrable and 
\begin{align*}
\mathbb{E}[v^2( s,X^{t,x}_{s})]  
\leq &  c\Big\{
1+\mathbb{E}\Big[\Big(\int_{t}^{s}\xi_{\sigma}dW_{\sigma}\Big)
^{2}\Big]+\mathbb{E}\Big[\Big| \int
_{t}^{s}\psi_{\sigma}d\sigma\Big| ^{2}\Big]\Big\} \\
\leq&  c\Big\{ 1+\mathbb{E}\Big[\int_{t}^{s}|
\xi_{\sigma}|
_{H}^{2}d\sigma\Big] + \mathbb{E}\Big[\Big| \int
_{t}^{s}\psi_{\sigma}d\sigma\Big| ^{2}\Big]\Big\} 
< \infty.
\end{align*}

Therefore, by the Cauchy--Schwarz inequality the expectation of
$v(s,X^{t,x}_{s}) \displaystyle{\int_{s-\delta}^{s}}
\xi_{\sigma}dW_{\sigma}$ is well defined. We claim that 
\begin{eqnarray}\label{second}
&\mathbb{E}\Big[ v\left(s,X^{t,x}_{s}\right)
\displaystyle\int_{s-\delta}^{s}\xi_{\sigma
}dW_{\sigma}\Big] =\displaystyle \frac{d}{d\varepsilon}_{\mid\varepsilon=0}%
\mathbb{E}_{Q_{\varepsilon}}\left[ v\left( s,X^{t,x}_{s}\right) \right].
\end{eqnarray}
As a matter of fact, 
\begin{align*}
&\frac{d}{d\varepsilon}_{\mid\varepsilon=0}%
\mathbb{E}_{Q_{\varepsilon}}\left[ v\left( s,X^{t,x}_{s}\right) \right] =
\frac{d}{d\varepsilon}_{\mid\varepsilon=0}\mathbb{E}\Big[
v\left( s,X^{t,x}_{s}\right) \exp\Big(
\varepsilon\int_{s-\delta}^{s}\xi_{\sigma }{\it
dW}_{\sigma}-\frac{\varepsilon^{2}}{2}\int_{s-\delta}^{s}\Vert
\xi_{\sigma}\Vert _{H}^{2}d\sigma\Big) \Big] \\
& =\lim_{\varepsilon\rightarrow0}\mathbb{E}\Big[ v\left(
s,X^{t,x}_{s}\right)
\frac{1}{\varepsilon}  \Big\{\exp\Big( \varepsilon%
{\int_{s-\delta}^{s}}
\xi_{\sigma}dW_{\sigma}-\frac{\varepsilon^{2}}{2}%
{\int_{s-\delta}^{s}}
\left| \xi_{\sigma}\right| _{H}^{2}d\sigma\Big) -1\Big\}\Big]  
=\mathbb{E}\Big[ v\left( s,X^{t,x}_{s}\right)
\int_{s-\delta}^{s}\xi_{\sigma }{\it dW}_{\sigma}\Big],
\end{align*}
where in the last passage 
we have used the dominated
convergence theorem being $\xi$  bounded.
\newline Now notice that,  in $(\Omega,\mathcal{F},Q_{\varepsilon})$, $X^{t,x}$ is a mild
solution to the equation
\[
dX^{t,x}_{\tau}=AX^{t,x}_{\tau}d\tau+F\left(X^{t,x}_{\tau}\right) d\tau+(-A)^{-\alpha}\varepsilon
\xi_{\tau}\left( W_{\cdot}^{\varepsilon}\right) d\tau+(-A)^{-\alpha}dW_{\tau
}^{\varepsilon},\text{ \ \ \ }\tau\in\left[ s-\delta,T\right].
\]
On the other hand,  in
$(\Omega,\mathcal{F},\mathbb{P})$, we consider the process
$X^{\varepsilon}$ which is a mild solution to the equation
\[
\left\{
\begin{array}[c]{@{}l}%
dX_{\tau}^{\varepsilon}=AX_{\tau}^{\varepsilon}d\tau+F\left( X_{\tau
}^{\varepsilon}\right) d\tau+(-A)^{-\alpha}\varepsilon\xi_{\tau}\left( W_{\cdot}\right)
d\tau+(-A)^{-\alpha}dW_{\tau},\text{ \ \ \ }\tau\in\left[ s-\delta,T\right], \\
 X_{s-\delta}^{\varepsilon}=X^{t,x}_{s-\delta}.
\end{array}
\right.
\]
Then the process $X^{t,x}$ under $Q_{\varepsilon}$ and the process
$X^{\varepsilon}$
under $\mathbb{P}$ have the same law, so \eqref{second} yields 
\begin{eqnarray}\label{third}
&\mathbb{E}\Big[ v\left( s,X^{t,x}_{s}\right)
\displaystyle\int_{s-\delta}^{s}\xi_{\sigma
}dW_{\sigma}\Big] = \displaystyle\frac{d}{d\varepsilon}_{\mid\varepsilon =0}\mathbb{E}\left[
v\left( s,X_{s}^{\varepsilon}\right) \right].
\end{eqnarray}
Let us set $\overset{\cdot}{X}_{\tau}:=\frac{d}{d\varepsilon}
_{\mid\varepsilon=0}X_{\tau}^{\varepsilon}$ and $\Delta^{\varepsilon}X_{\tau}=\frac{X_{\tau}^{\varepsilon}-X_{\tau}%
}{\varepsilon}$, $\mathbb{P}$-a.s. for any $\tau\in [s-\delta,T]$.
Arguing as in \cite{MasBan}, one can prove that
\begin{align}
&\lim_{\varepsilon \rightarrow 0}|
\Delta^{\varepsilon}X_{\tau}-\overset{\cdot}{X}_{\tau }|
_{E}=0, \quad
\overset{\cdot}{X}_{\tau}=\int_{s-\delta}^{\tau}\nabla_{x} X_\tau^{\sigma,X_{\sigma}^{t,x}} (-A)^{-\alpha}\xi_{\sigma}d\sigma, \quad \tau\in [s-\delta,T], \quad \P\textup{-a.s.}
 \label{Xdotid}
\end{align}
Formula \eqref{Xdotid} in turn allows to show that 
\begin{align*}
&\frac{d}{d\varepsilon}_{\mid\varepsilon=0}\mathbb{E}\Big[ v\left(
s,X_{s}^{\varepsilon}\right) \Big] =\mathbb{E}\Big[ \nabla_{x}
v\left( s,X^{t,x}_{s}\right) \overset{\cdot}{X}_{s}\Big]= \mathbb{E}\Big[ \nabla_{x}
v\left(s,X^{t,x}_{s}\right) \int_{s-\delta}^{s}\nabla_{x} X_s^{\sigma,X_{\sigma}^{t,x}} (-A)^{-\alpha}\xi_{\sigma}d\sigma\Big], \end{align*}
so that formula \eqref{third} gives 
\begin{eqnarray}\label{fourth}
&\mathbb{E}\Big[ v\left( s,X^{t,x}_{s}\right)
\displaystyle\int_{s-\delta}^{s}\xi_{\sigma
}dW_{\sigma}\Big] = \mathbb{E}\Big[ \nabla_{x}
v\left( s,X^{t,x}_{s}\right) \displaystyle\int_{s-\delta}^{s}\nabla_{x} X_s^{\sigma,X_{\sigma}^{t,x}} (-A)^{-\alpha}\xi_{\sigma}d\sigma\Big].
\end{eqnarray}
By \eqref{fourth} we have 
\begin{align*}
&\lim_{\delta\rightarrow0}\frac {1}{\delta}\mathbb{E}\left[ v\left( s,X^{t,x}_{s}\right) \int_{s-\delta}^{s}\xi_{\sigma}dW_{\sigma}\right] 
 =\lim_{\delta\rightarrow 0}\frac{1}{\delta}\mathbb{E}\left[
\nabla_{x} v\left( s,X^{t,x}_{s}\right) \int_{s-\delta}^{s}\nabla_{x} X_s^{\sigma,X_{\sigma}^{t,x}} (-A)^{-\alpha}\xi_{\sigma}d\sigma\right] \\
& =\mathbb{E}\Big[ \nabla_{x} v\left(s,X^{t,x}_{s}\right) \nabla_{x}
X_s^{s,X_s^{t,x}}  (-A)^{-\alpha}\xi_{s}\Big]  =\mathbb{E}\Big[ \nabla_{x}\big(v(s,X^{t,x}_{s})\big)
 (-A)^{-\alpha}\xi_{s}\Big]
\end{align*}
so \eqref{toprove_1} is proved.

It remains to prove \eqref{toprove_2}. 
Recalling 
 identifications \eqref{Id_Y_Markov_prop}-\eqref{Id_Z_Markov_prop}, 
we have 
\begin{align}\label{thirdBIS}
&\frac{1}{\delta}\mathbb{E}\Big[\int_{s-\delta}^{s}\psi(\sigma, X_\sigma^{t,x}, Y_\sigma^{t,x},Z_\sigma^{t,x})d\sigma
\int_{s-\delta}^{s}\xi_{\sigma
}dW_{\sigma}\Big]\notag\\
&= \frac{1}{\delta}\mathbb{E}\Big[\int_{s-\delta}^{s}\psi(\sigma, X_\sigma, v(\sigma, X_\sigma), u(\sigma, X_\sigma))d\sigma
\int_{s-\delta}^{s}\xi_{\sigma
}dW_{\sigma}\Big] \notag\\
&= \frac{1}{\delta}\frac{d}{d\varepsilon}_{\mid\varepsilon =0}\mathbb{E}\Big[
\int_{s-\delta}^{s}\psi(\sigma, X^\varepsilon_\sigma, v(\sigma, X^\varepsilon_\sigma), u(\sigma, X^\varepsilon_\sigma))d\sigma\Big], 
\end{align} 
which is the analogous of formula \eqref{third} with $f_\delta:= \displaystyle\int_{s-\delta}^{\cdot}\psi_{\sigma}d\sigma$ in place of $v$.  
Now we notice that
\begin{align}\label{coll_1}
	\frac{1}{\delta}\frac{d}{d\varepsilon}_{\mid\varepsilon =0}\psi(\sigma, X^\varepsilon_\sigma, v(\sigma, X^\varepsilon_\sigma), u(\sigma, X^\varepsilon_\sigma)) = \nabla_x \psi_\sigma \, \frac{\overset{\cdot}{X}_{\sigma}}{\delta} + \nabla_y \psi_\sigma \,\frac{\overset{\cdot}{Y}_\sigma}{\delta} + \nabla_z \psi_\sigma \,\frac{\overset{\cdot}{Z}_\sigma}{\delta}, \quad \sigma \in [s-\delta, T], 
\end{align}
where we have used the notation 
$$
 (\overset{\cdot}{Y}, \overset{\cdot}{Z}) :=\Big(\frac{d}{d\varepsilon}_{\mid\varepsilon=0} Y, \frac{d}{d\varepsilon}_{\mid\varepsilon=0} Z\Big) =(\nabla_x v(\sigma,X^\varepsilon_{\sigma})\overset{\cdot}{X}, \nabla_x u(\sigma,X^\varepsilon_{\sigma})\overset{\cdot}{X}) . 
$$
 By \eqref{Xdotid} and \eqref{estimate_der}, 
we have  
\begin{align}\label{coll_2}
	\frac{|\overset{\cdot}{X}_{\tau}|_{E}}{\delta} \leq \frac{1}{\delta}\Big |\int_{s-\delta}^{\tau}\nabla_{x} X_\tau^{\sigma,X_\sigma^{t,x}}  (-A)^{-\alpha}\xi_{\sigma}d\sigma\Big|_{E}\leq C, \quad \tau \in [s-\delta, T].
\end{align}
On the other hand, the pair of processes 
$(\overset{\cdot}{Y}, \overset{\cdot}{Z})$
is solution  to the FBSDE 
\begin{align}\label{eq_dotYdotZ}
\left\{
\begin{array}{rcl}
-d\overset{\cdot}{Y}_\tau&=& \nabla_x \psi(\tau,X_\tau^{t,x},Y_\tau^{t,x},Z^{t,x}_\tau)\overset{\cdot}{X}_\tau  d\tau+\nabla_y \psi(\tau,X_\tau^{t,x},Y_\tau^{t,x},Z^{t,x}_\tau)\overset{\cdot}{Y}_\tau  d\tau\vspace{2mm} \\
& & +\nabla_z \psi(\tau,X_\tau^{t,x},Y_\tau^{t,x},Z^{t,x}_\tau)\overset{\cdot}{Z_\tau} d\tau-\overset{\cdot}{Z}_\tau  dW_\tau, \quad \tau\in[s-\delta,T], \vspace{2mm} \\
\overset{\cdot}{Y}_T  & =& \nabla_x \Phi(X_T^{t,x})\overset{\cdot}{X_T}, \vspace{2mm} \\
d \overset{\cdot}{X}_\tau & = & A \overset{\cdot}{X}_\tau d\tau+\nabla F(X_\tau^{t,x})\overset{\cdot}{X}_\tau d\tau, \quad \tau\in[s-\delta,T],\vspace{1mm} \\
\overset{\cdot}{X}_{s-\delta} & = & 0.
\end{array}
\right.
\end{align}
Moreover, taking into account \eqref{coll_2} and the linearity of the BSDE \eqref{eq_dotYdotZ}, we get that the pair $(\overset{\cdot}{Y}, \overset{\cdot}{Z})$ satisfies the estimates 
\begin{align}
&\sup_{\tau\in [s-\delta,\,T]}\frac{|\overset{\cdot}{Y}_\tau|}{\delta}	 \leq C, \quad \P\textup{-a.s.},\label{est_gradYNEW}\\
&\frac{1}{\delta}\E\Big[\int_{s-\delta}^T 
 |\overset{\cdot}{Z}_\tau|_H^2
  d\tau \Big] \leq C.\label{est_gradY_gradZNEW}
\end{align}
By Hypothesis \ref{ip-psiphi}, 
\begin{equation}\label{coll_3}
|\nabla_x \psi_\sigma|_{E^*} \leq C, \quad |\nabla_y \psi_\sigma| \leq C, \quad  |\nabla_z \psi_\sigma|_H \leq C(1 + |z|_H), \quad \sigma\in[s-\delta,T].
\end{equation}
Therefore, collecting \eqref{coll_1}-\eqref{coll_2},  \eqref{est_gradYNEW}-\eqref{est_gradY_gradZNEW} and \eqref{coll_3}, \eqref{thirdBIS} gives 
\begin{align*}
&\frac{1}{\delta}\mathbb{E}\Big[\int_{s-\delta}^{s}\psi(\sigma, X_\sigma^{t,x}, Y_\sigma^{t,x},Z_\sigma^{t,x})d\sigma
\int_{s-\delta}^{s}\xi_{\sigma
}dW_{\sigma}\Big]\\
&= \frac{1}{\delta}\mathbb{E}\Big[
\int_{s-\delta}^{s}\frac{d}{d\varepsilon}\psi(\sigma, X^\varepsilon_\sigma, v(\sigma, X^\varepsilon_\sigma), u(\sigma, X^\varepsilon_\sigma) )d\sigma\Big]\\
&=\mathbb{E}\Big[
\int_{s-\delta}^{s}\Big(\nabla_x \psi_\sigma \, \frac{\overset{\cdot}{X}_{\sigma}}{\delta} + \nabla_y \psi_\sigma \,\frac{\overset{\cdot}{Y}_\sigma}{\delta}  + \nabla_z \psi_\sigma \,\frac{\overset{\cdot}{Z}_\sigma}{\delta} \Big)d\sigma\Big]\\
&\leq C \delta + C  \mathbb{E}\Big[
\int_{s-\delta}^{s} |\nabla_y \psi_\sigma| \,\frac{|\overset{\cdot}{Y}_\sigma|}{\delta}d\sigma\Big]+ C \mathbb{E}\Big[
\int_{s-\delta}^{s} |\nabla_z \psi_\sigma|_H \,\frac{|\overset{\cdot}Z_\sigma|_H}{\delta} d\sigma\Big]\leq C \delta +\E\Big[\int_{s-\delta}^s |Z_\sigma|^2_Hd \sigma\Big]
\end{align*}
which goes to zero as $\delta$ goes to zero. This shows that \eqref{toprove_2} holds true and concludes the proof. 
\endproof

\begin{corollary}\label{rem:identif-Z}
Under the assumptions of Theorem \ref{prop-identific-Z} we have  
$$
Z_{s}^{t,x}h=\nabla v(s,X_{s}^{t,x})(-A)^{-\alpha}h, \quad h \in H, \,\, {\rm for \,\,a.e.} \ s\in[t,T],  \,\, \mathbb{P}\textup{-a.s.},
$$
where  $\nabla_x v(s,x)(-A)^{-\alpha}$ denotes an extension of the operator  $\nabla_x v(s,x)(-A)^{-\alpha}: E_0 \rightarrow \R$ to the whole space $H$.
Moreover,    there exists a constant $C$, that may depend also on $\nabla_x\phi$,
$\nabla_x\psi$ and $L_\psi$, such that
\begin{equation}\label{stimaZ-diffle}
\vert Z_s^{t,x}\vert_H \leq C, \quad {\rm for \,\,a.e.} \ s\in[t,T], \ \P\textup{-a.s.}
\end{equation} 

\end{corollary}
\proof 
Since 
$E_0$ is dense in $H$,  
by \eqref{Id_Z} in  Theorem \ref{prop-identific-Z},  for almost every $s \in [0,\,T]$ and almost surely with respect to the law of $X$, 
the operator $\nabla_x v(s,x)(-A)^{-\alpha}: E_0 \rightarrow \R$ extends to an operator defined on the whole $H$, which we still denote  $\nabla v(s,x) (-A)^{-\alpha}$.

\noindent Moreover, from \eqref{Id_Z} and by the Markov property, we get 
\[
Z_\sigma^{t,x}=Z_\sigma^{\sigma, X_\sigma^{t,x}}=\nabla_x Y_\sigma^{\sigma,k}\mid_{k=X_\sigma^{t,x}} (-A)^{-\alpha}, 
\quad {\rm for \,\,a.e.} \ \sigma\in[0,T], \ \P\textup{-a.s.}
\]
The conclusion \eqref{stimaZ-diffle} follows from the fact that     $\sup_{ \sigma} |\nabla_x Y_\sigma^{\sigma,k}| \leq C$ by \eqref{new_est_gradY_gradZ}, where $C$ is a constant that does not depend on $k$. 
\endproof

Now we use the previous result to give a priori estimates on $Z^{t,x}$.
\begin{proposition}\label{prop-aprioriZ}
Assume that Hypotheses \ref{ip_forward} and \ref{ip-psiphi} hold true, and for any   $(t, x) \in [0, T ] \times E$, let $(X^{t,x},Y ^{t,x} , Z^{ t,x} )$ be the solution to the FBSDE \eqref{fbsde}. 
Then  there exists a positive constant $C_T$ only depending on $T,\;A,\; F\;K_\phi,\;L_\psi$, $K_\psi$ such that
\begin{align}
|Z_t^{t,x}h| &\leq C_T (T-t)^{-1/2}|h|_H, \quad\P\textup{-a.s.},\ h\in H, \label{stimabismut}\\
|\nabla_x Y_t^{t,x}h|&\leq C_T(T-t)^{-1/2-\alpha}|h|_{H}, \quad \P\textup{-a.s.}, \ h\in H.\label{stimabismut2}
\end{align}
\end{proposition}
\proof
In the following $C_T$ will denote a positive constant which may depend on $T,L_\psi, K_\psi$, $K_\phi$ but  not on $\nabla_x\phi$, and that may vary from line to line. We fix $(t,x) \in [0\,T] \times E$.

We start by proving estimate \eqref{stimabismut}. 
We first  take $\phi$ and $\psi$ differentiable. 
By Proposition \ref{P:BC_Diff_BSDE},  the triple of processes $(X^{t,x}, Y^{t,x},Z^{t,x})$ is G\^ateaux differentiable as a map from $E$ with values in
$\cals^2((t,T];E)\times \cals^2([t,T])\times \calm^2([t,T];H)$, and  for any  $h\in E_0$, the triple of processes $(\nabla_x X^{t,x}, \nabla_x Y^{t,x} h, \nabla_x Z^{t,x} h)$ is solution to \eqref{eq_gradYgradZ}, and satisfies estimate 
\eqref{new_est_gradY_gradZ}.

Let us now introduce the process
\begin{align*}
W^{\Q}_\tau:=W_\tau-\int_t^\tau \nabla_z \psi(s,X_s^{t,x},Y_s^{t,x},Z_s^{t,x})ds, \quad \tau\in [t,T],
\end{align*}
where $\Q$ is the probability measure such that $W^{\Q}$ is a Brownian motion in $(\Omega,\mathcal F,(\mathcal F_t)_{t\geq0}, \Q)$.

\noindent Let us fix $h \in E_0$. 
 Arguing as in \cite[Proposition 3.6]{Mas-spa} it follows that 
\begin{align*}
F_\tau^{t,x}h
:= & e^{\int_t^\tau\nabla_y\psi(s,X^{t,x}_s,Y^{t,x}_s,Z^{t,x}_s)ds}\nabla_xY^{t,x}_\tau h \\
& +\int_t^\tau e^{\int_t^s \nabla_y\psi(r,X^{t,x}_r,Y^{t,x}_r,Z^{t,x}_r)dr}\nabla_x\psi(s,X^{t,x}_s,Y^{t,x}_s,Z^{t,x}_s)\nabla_x X_s^{t,x} hds, \quad \tau \in [t,\,T].
\end{align*}
Therefore, $(|F_\tau^{t,x}h|^2)_{\tau \in [t,T]}$ is a $\Q$-submartingale, which
implies, thanks to  identification formula \eqref{Id_Z},   that
\begin{align}
\label{fico}
\mathbb E^{\Q}\Big[\int_t^\tau |F_s^{t,x}h|^2ds\Big]
\geq (\tau-t)|F_t^{t,x}h|^2=(\tau-t)|Z_t^{t,x}(-A)^{\alpha}h|, \quad\tau \in [t,\,T].
\end{align}
Further, since $\psi$ is differentiable and Lipschitz continuous with respect to $x$  and $y$, and $\nabla_xX^{t,x}$ is bounded (see \eqref{estimate_der_H}), we deduce that 
\begin{align}
\label{martina}
|F_\tau^{t,x}h|^2\leq C_T \left(|\nabla_x Y_{\tau}^{t,x}h|^2+|h|_H^2\right), \quad \tau \in [t,\,T], \ \P\textup{-a.s.}
\end{align}
It remains to estimate $|\nabla_xY_\tau^{t,x}h|$. To this aim, we recall the well-known estimate
\begin{align}
\label{tiziano}
\mathbb E^{\Q}\Big[\Big(\int_t^\tau |Z_s^{t,x}|_{H}^2ds \Big)^{p/2}\Big]\leq C\|\Phi\|^p_\infty,\quad \tau \in [t,\,T], 
\end{align}
for some $C>0$ and any $p<+\infty$. Formulas \eqref{tiziano}, \eqref{Id_Z} and \eqref{sper_est_der2} give
\begin{align*}
E^{\Q}\Big[\int_t^\tau|\nabla Y_s^{t,x}h|^2 ds\Big]
\leq \mathbb E^{\Q}\Big[\int_t^\tau|Z^{t,x}_s|^2_H|(-A)^\alpha \nabla X_s^{t,x}h|_{H}^2 ds\Big]
\leq {c\|\Phi\|_\infty^2| (-A)^\alpha h|_{H}^2},\quad \tau \in [t,\,T], 
\end{align*}
which, together with \eqref{fico} and \eqref{martina}, allows us to conclude that
\begin{align*}
|Z_t^{t,x} (-A)^\alpha h|^2
\leq \frac{C_T}{T-t}|(-A)^\alpha h|_{H}^2, \quad h \in E_0.
\end{align*}
Let now fix $h \in H$. We notice that in this case  we can write $h=(-A)^\alpha(-A)^{-\alpha} h$. Therefore,
\begin{align*}
|Z_t^{t,x} h|^2
= &|Z_t^{t,x}(-A)^\alpha (-A)^{-\alpha}h|^2
\leq \frac{C_T}{T-t}|(-A)^\alpha (-A)^{-\alpha} h|_{H}^2
\leq \frac{C_T}{T-t}|h|_H^2, \quad h \in H, 
\end{align*}
which provides \eqref{stimabismut} in the case of $\psi$ and $\phi$ differentiable. 

\noindent Finally, the case $\psi$ and $\phi$ non differentiable can be obtained by approximating $\psi$ and $\phi$ with  $\psi_n$ and $\phi_n$ in \eqref{psi-approx} and \eqref{phi-approx}, respectively. For the proof we refer to \cite[Proposition 3.6]{Mas-spa}.

 Let us now prove estimate \eqref{stimabismut2}. Again, at first we prove the result when $\psi$ and $\phi$ are differentiable and then we generalize it by approximation. Let us fix $h\in E_0$. For any $t<\eta<\tau\leq T$, the submartingale property of $(|F_s^{t,x}|_H^2)_{s\in[t,T]}$ gives
\begin{align}\label{intgradYge}
\mathbb E^{\Q}\Big[\int_\eta^\tau |F_s^{t,x}h|^2ds\Big]
=\int_\eta^\tau \mathbb E^{\Q}[|F_s^{t,x}h|^2]ds
\geq \int_\eta^\tau \mathbb E^{\Q}[|F_\eta^{t,x}h|^2]ds
= (\tau-\eta)\mathbb E^{\Q}[|F_\eta^{t,x}h|^2].
\end{align}
Moreover, for any $\tau\in(t,T]$ we split
\begin{align}\label{Yrewr}
\mathbb E^{\Q}\Big[\int_t^\tau |F_s^{t,x}h|^2ds\Big]
= & \mathbb E^{\Q}\Big[\int_t^{(t+\tau)/2} |F_s^{t,x}h|^2ds\Big]+\mathbb E^{\Q}\Big[\int_{(t+\tau)/2}^\tau |F_s^{t,x}h|^2ds\Big]=:I_1+I_2.
\end{align}
Let  us evaluate separately $I_1$ and $I_2$. Concerning $I_1$, 
identification formula \eqref{Id_Z}, \eqref{tiziano} and \eqref{sper_est_der} give
\begin{align}\label{intgradYle}
&\mathbb E^{\Q}\Big[\int_\eta^\tau |\nabla Y_s^{t,x}h|^2ds\Big]\notag\\
&\leq \mathbb E^{\Q}\Big[\int_\eta^\tau |Z^{t,x}_s|_H^2|(-A)^\alpha\nabla X_s^{t,x}h|_{H}^2ds\Big]\notag\\
&+ C_T(\tau-t)^{-2\alpha}\Big(|x|_E^{2m+1}+\E^{\Q}\Big[\sup_{\tau\in[t,T]}|w^A(\tau)|_E^{2m+1}\Big]\Big)^2|h|_H^2\mathbb E^{\Q}\Big[\int_\eta^\tau |Z^{t,x}_s|_H^2 ds\Big] \notag\\
&\leq  C_T \|\Phi\|^{2}_\infty |h|_H^2 (\eta-t)^{-2\alpha}.
\end{align} 
Hence, from  \eqref{martina}, \eqref{intgradYge} and \eqref{intgradYle}  it follows that
\begin{align}
\label{blink}
E^{\Q}[|F_\eta^{t,x}h|^2_H]\leq C_T \|\Phi\|^{2}_\infty|h|_H^2\Big(\frac{(\eta-t)^{-2\alpha}}{(\tau-\eta)}+ 1\Big), \quad t<\eta<\tau\leq T.
\end{align}
By applying Fubini's theorem and  \eqref{blink}, we infer that
\begin{align}\label{estintT1}
I_1
&\leq  C_T\|\Phi\|_\infty^2|h|_H^2 \Big[(\tau-t)^{-1}\int_t^{(t+\tau)/2}(s-t)^{-2\alpha}ds + (\tau-t)^{-1}\int_t^{(t+\tau)/2}ds\Big]\notag\\
&=C_T|h|_H^2 [(\tau-t)^{-2\alpha}+ (\tau-t)] \leq   C_T |h|_H^2(\tau-t)^{-2\alpha}, \quad  \quad t<\eta<\tau\leq T.
\end{align}
As far as $I_2$ is concerned, we take advantage from \eqref{martina} and  \eqref{intgradYle}. Then, for  $t<\eta<\tau\leq T$ we get 
\begin{align}\label{estintT2}
I_2
\leq C_T |h|_H^2 \left[(\tau-t)+ (\tau-t)^{-2\alpha}\right]=  C_T |h|_H^2[(\tau-t)^{-2\alpha}+(\tau-t)] \leq   C_T |h|_H^2(\tau-t)^{-2\alpha}.
\end{align}
Thus collecting \eqref{intgradYge}, \eqref{estintT1} and  \eqref{estintT2}, we have
\begin{align*}
(\tau-\eta)\mathbb E^{\Q}[|F_\eta^{t,x}h|^2] \leq \mathbb E^{\Q}\int_t^\tau |F_s^{t,x}h|^2ds
\leq C_T |h|_H^2 (\tau-t)^{-2\alpha},  \quad t<\eta<\tau\leq T,
\end{align*} 
so that 
\begin{align}
\label{tria}
\mathbb E^{\Q}[|F_\eta^{t,x}h|^2]\leq C_T\frac{1}{(\tau-\eta)(\tau-t)^{2\alpha}}|h|_H^2, \quad h \in E_0, \,\,\,t\leq \eta <\tau\leq T.
\end{align}
Let us now fix  $h\in H$, and let us consider a sequence $(h_n)\subset E_0$ such that $h_n\rightarrow h$ as $n\rightarrow+\infty$ in $H$. Taking \eqref{tria} with $h$ replaced by $h_n$ and letting $n\rightarrow +\infty$, it follows that
\begin{align}
\label{tria2}
\mathbb E^{\Q}[|F_\eta^{t,x}h|^2]\leq \frac{C_T}{(\tau-\eta)(\tau-t)^{2\alpha}}|h|_H^2, \quad h \in H,\,\,\,t\leq \eta <\tau\leq T.
\end{align}
Inequality \eqref{stimabismut2} follows from \eqref{tria2} by taking $\tau=T$ and $\eta=t$.
\endproof

\section{The Bismut-Elworthy formula and the semilinear Kolmogorov equation: the Lipschitz case}
\label{sez-Bismut-lip}
Recall that  we deal with a process $X$ taking values in a Banach space and solution to equation (\ref{forward}), with  special diffusion operator $
(-A)^{-\alpha}$ with pseudo-inverse $(-A)^\alpha$ which is  not bounded.

\noindent In the present section we  adequate  to  our framework the results in \cite{futeBismut}.
More precisely, in  Subsection \ref{B_lips} we present a nonlinear version of the Bismut-Elworthy
formula in the case of  Lipschitz generator, which extends   the one provided in \cite{futeBismut} in the case of a  process $X$ taking  values in a Hilbert space, and  with a bounded  diffusion operator 
 with bounded inverse. 
Providing the Bismut-Elworthy  formula in the case of  Lipschitz generator is a   fundamental step in order to obtain the analogous  formula in the quadratic case. Moreover, it allows us to give an  existence and uniqueness result in the Banach framework for the semilinear Kolmogorov related to the process $X$, and with coefficients $\phi$ and $\psi$ not necessarily differentiable, see Subsection \ref{K_lip}.

For $0\leq t< s\leq T$ and $h\in H$ we define the real valued random variables
\begin{equation}\label{def-U}
 U^{h,t,x}_s:=\dfrac{1}{s-t}\int_t^s \<(-A)^{\alpha
 }\nabla_xX_r^{t,x}h,dW_r \>.
\end{equation}
Notice that,  for any $h \in H$, the process $ U^{h,t,x}$ is well defined  thanks to  formula (\ref{weight_estimate_derivative}) in  Proposition \ref{P1.9}. 
In what follows we prove some useful estimates on the process $ U^{h,t,x}$.

\begin{lemma}\label{lemma:boundU-q}
 Assume that Hypotheses \ref{ip_forward}
hold true. For any $(t,x) \in [0,\,T] \times E$,  let $X^{t,x}$ be the unique mild solution to \eqref{forward}.  Then, for any $h\in H$ and for any $q\geq 1$,
\begin{equation}\label{bound-Uq}
 \left(\E[\vert U^{h,t,x}_s \vert^q]\right)^{1/q}\leq
C\left(s-t \right)^{-(\frac{1}{2} + \alpha)}\vert h\vert_H, 
\end{equation}
and also
\begin{equation}
\Big(\E\Big[\sup_{s\in[\frac{t+T}{2},T]}\vert U^{h,t,x}_s \vert^q\Big]\Big)^{1/q}\leq C(T-t)^{-(\frac{1}{2} + \alpha)} |h|_H.\label{bound-Uq-sup}
\end{equation}
\end{lemma}
\proof  We compute
\begin{align*}
 \E[\vert U^{h,t,x}_s \vert^q] &=
\E\Big[\Big\vert \dfrac{1}{s-t}\int_t^s \<(-A)^\alpha\nabla_xX_r^{t,x}h,dW_r \>\Big\vert^q\Big] \leq  \dfrac{1}{(s-t)^q}\E \Big[\Big(\int_t^s \vert(-A)^\alpha\nabla_x X_r^{t,x}h\vert^2\,dr \Big)^{q/2}\Big]\\ \nonumber
& \leq \dfrac{1}{(s-t)^q}C ((s-t)^{1-2\alpha} |h|_H^2)^{q/2}=C (s-t)^{-q(\frac{1}{2} + \alpha)} |h|_H^q,
\end{align*}
where in the latter inequality we have used formula \eqref{weight_estimate_derivative} of Proposition \ref{P1.9} with $\varepsilon = \alpha$.
Analogously, we have 
\begin{align*}
 \E\Big[\sup_{s\in[\frac{t+T}{2},T]}\vert U^{h,t,x}_s \vert^q \Big]
\leq C  \dfrac{1}{(T-t)^q}\E \Big[\Big(\int_t^T \vert (-A)^\alpha\nabla_xX_r^{t,x}h\vert^2\,dr \Big)^{q/2} \Big]\leq C (T-t)^{-q(\frac{1}{2} + \alpha)} |h|_H^q.
\end{align*}
\endproof

\subsection{The Bismut formula}\label{B_lips}

We can now give a version of the Bismut-Elworthy formula in the case of Lipschitz generator and in the Banach space framework. 
We consider only the case of final datum $\phi$ and generator $\psi$ bounded with respect to $x$, since 
we aim to treat such a model  in the quadratic case.
 We start with  
 the case when the coefficients are also differentiable. An analogous result is proved in \cite{futeBismut} in the Hilbert space framework using the Malliavin calculus. Since here the process $X$ takes its values in a Banach space, we avoid  the use of the Malliavin calculus,   by exploiting instead techniques similar to the ones used in the proof of Theorem \ref{prop-identific-Z}.

In the rest of the section we will assume the following, that substitutes Hypothesis \ref{ip-psiphi}.
\begin{hypothesis}\label{ip-psiphi-lip}
 The functions  $\mathbb{\phi}:E\rightarrow\mathbb{R}$  and $\psi:[0,T]\times E\times \R\times H
\rightarrow\mathbb{R}$  in \eqref{fbsde} satisfy the following.
\begin{itemize}
\item[(i)] $\mathbb{\phi}$ is continuous,  and there exist a nonnegative constant $K_\phi$ such that 
$\vert \phi(x)\vert\leq K_\phi$
for every  $x\in E$.
\item[(ii)]	$\psi$ 
is measurable and,  for every fixed $t\in[0,T]$, the map $\psi(t,\cdot,\cdot,\cdot):E\times\R\times H\rightarrow\R$
is continuous.
 Moreover, there exist nonnegative constants $L_\psi$ and $K_\psi$
such that
\begin{align*}
&\vert \psi(t,x_1,y_1,z_1)- \psi(t,x_2,y_2,z_2)\vert\leq L_\psi\left(
\vert x_1-x_2\vert_E+
\vert y_1-y_2\vert+\vert z_1-z_2\vert_{H}\right),\\
&\vert \psi(t,x,0,0)\vert\leq K_\psi,
\end{align*}
for every $t\in[0,T]$, 
$x_1, x_2 \in E$, 
$y_1,y_2\in\R$ and $z_1,z_2\in H$.
\end{itemize}
\end{hypothesis}
\begin{theorem}\label{teoBismutLip-diffle}
Let Hypotheses \ref{ip_forward} and \ref{ip-psiphi-lip} hold
true, and for any   $(t,x) \in  [0,\,T] \times E$, let $(X^{t,x},Y^{t,x},Z^{t,x})$ be a solution of the forward-backward system  \eqref{fbsde}, and let $U^{h,t,x}$
be the process defined in \eqref{def-U}. 
Assume moreover  that $\phi$ is G\^ateaux differentiable with  bounded derivative,  and that $\psi$ is G\^ateaux differentiable with respect to $x$, $y$ and $z$. 
Then for $ t\leq s\leq T$, $x\in E,\,h\in H$,
\begin{equation}\label{Bismut-lip}
\E\left[ \nabla_x\,Y^{t,x}_sh \right]=
\E\Big[\int_s^T\psi\left(r,X_r^{t,x},Y_r^{t,x},Z_r^{t,x}\right)U^{h,t,x}_r\,dr\Big]
+\E\left[  \phi(X_T^{t,x})U^{h,t,x}_T\right].
\end{equation}
\end{theorem}
\proof 
Let $\xi $ be a given square integrable $E_0$-valued predictable process,  and 
$X^{\varepsilon,t,x}$ be  a mild solution to the equation
\begin{equation}\label{X_epsilon_Bism}
\left\{
\begin{array}[c]{@{}l}%
dX_{\tau}^{\varepsilon,t,x}=AX_{\tau}^{\varepsilon,t,x}d\tau+F\left(X_{\tau
}^{\varepsilon,t,x}\right) d\tau+(-A)^{-\alpha}\varepsilon\xi_{\tau}
d\tau+(-A)^{-\alpha}dW_{\tau},\text{ \ \ \ }\tau\in\left[ t,T\right], \\
 X_{t}^{\varepsilon,t,x}=x.
\end{array}
\right.
\end{equation}
 We also consider the pair of processes $(Y^{\varepsilon,t,x}, Z^{\varepsilon,t,x})$ solution to the Markovian BSDE
\begin{equation}\label{YZ_epsilon_Bism}
\left\{
\begin{array}{ll}
-dY^{\varepsilon,t,x}_\tau=\psi(\tau,X_\tau^{\varepsilon,t,x},Y_\tau^{\varepsilon,t,x},Z^{\varepsilon,t,x}_\tau)d\tau-Z^{\varepsilon,t,x}_\tau  dW_\tau, & \tau\in[t,T], \vspace{2mm} \\
Y^{\varepsilon,t,x}_T =  \phi(X_T^{\varepsilon,t,x}). & 
\end{array}
\right.
\end{equation}
Arguing similarly to the proof of Theorem \ref{prop-identific-Z}, we define
\begin{equation}\label{XYZ_dot}
\overset{\cdot}{X}_{\tau}:=\frac{d}{d\varepsilon}
_{\mid\varepsilon=0}X_{\tau}^{\varepsilon,t,x},\; \overset{\cdot}{Y}_{\tau}:=\frac{d}{d\varepsilon}
_{\mid\varepsilon=0}Y_{\tau}^{\varepsilon,t,x},\;\overset{\cdot}{Z}_{\tau}:=\frac{d}{d\varepsilon}
_{\mid\varepsilon=0}Z_{\tau}^{\varepsilon,t,x}, \quad \tau\in[t,T], 
\end{equation}
which are solution to the forward-backward system 
(\ref{eq_dotYdotZ}) with $s-\delta =t$.
We  already 
know (see 
 formula (\ref{Xdotid}) with $s-\delta =t$) 
  that
\begin{equation}\label{Xdotid_Bism}
\overset{\cdot}{X}_{\tau}=\int_{t}^{\tau}\nabla_{x} X_\tau^{\sigma,X_\sigma^{t,x} }(-A)^{-\alpha}\xi_{\sigma}d\sigma, \quad \tau\in[t,T], \ \P\textup{-a.s.}
\end{equation}
Now we want to prove a similar identification for the pair $(\overset{\cdot}{Y},\overset{\cdot}{Z})$.
To this aim, for any $\sigma \in [t,T]$,  we consider the Markovian BSDE in (\ref{fbsde}) on the time interval $[\sigma,T]$, and with initial condition $y$ given at time $\sigma$; from Proposition \ref{P:BC_Diff_BSDE} we know that the derivative with respect to $y\in E$ in the direction $h\in E$ satisfies the following BSDE, that we write in integral form: for any  $\tau\in[t,T]$, $\P$-a.s., 
\begin{align}\label{eq_gradYgradZ_sigma-y}
\left\{
\begin{array}{rcl}
\nabla_x Y_\tau^{\sigma,y} h&= \nabla_x \phi(X_T^{\sigma,y})\nabla_x X_T^{\sigma,y}h-\displaystyle\int_\tau^T\nabla_x Z^{\sigma,x}_r h \,dW_r+\displaystyle\int_\tau^T \big(\nabla_x \psi(r,X_r^{\sigma,y},Y_r^{\sigma,y},Z^{\sigma,y}_r)\nabla_x X^{\sigma,y}_r h \vspace{2mm} \\
& +\nabla_y \psi(r,X_r^{\sigma,y},Y_r^{\sigma,y},Z^{\sigma,y}_r)\nabla_x Y^{\sigma,y}_r h \,dr+\nabla_z \psi(r,X_r^{\sigma,y},Y_r^{\sigma,y},Z^{\sigma,y}_r)\nabla_x Z^{\sigma,y}_r h \big) dr. 
\end{array}
\right.
\end{align}
 Let us take  $y=X_\sigma^{t,x}$ and  $h=(-A)^{-\alpha}\xi_\sigma$ in \eqref{eq_gradYgradZ_sigma-y}, and let us integrate  both sides with respect to $\sigma\in [t,\tau]$. By inverting the order of integration where necessary, and using the Markov property, it is immediate to get
\begin{align*}
&\int_t^\tau \nabla_x Y_\tau^{\sigma,X_\sigma^{t,x}}(-A)^{-\alpha}\xi_\sigma \,d\sigma \notag\\
&= \displaystyle\int_t^\tau\nabla_x \phi(X_T^{\sigma,X_\sigma^{t,x}})\nabla_x X_T^{\sigma,X_\sigma^{t,x}}(-A)^{-\alpha}\xi_\sigma \,d\sigma
-\displaystyle\int_\tau^T\Big (\displaystyle\int_t^\tau\nabla _xZ^{\sigma,X_\sigma^{t,x}}_r (-A)^{-\alpha}\xi_\sigma \,d\sigma \Big)dW_r\notag\\ 
&+\int_\tau^T \Big(\int_t^\tau \nabla_x\psi\left(r,X_r^{\sigma,X_\sigma^{t,x}},Y_r^{\sigma,X_\sigma^{t,x}}, Z^{\sigma,X_\sigma^{t,x}}_r\right)\nabla _xX^{\sigma,X_\sigma^{t,x}}_r (-A)^{-\alpha}\xi_\sigma\,d\sigma\Big) dr \notag \\
 & +\int_t^\tau \Big(\int_t^\tau\nabla_y \psi(r,X_r^{\sigma,X_\sigma^{t,x}},Y_r^{\sigma,X_\sigma^{ t,x}},Z^{\sigma,X_\sigma^{t,x}}_r)\nabla _xY^{\sigma,X_\sigma^{t,x}}_r (-A)^{-\alpha}\xi_\sigma \,d\sigma\Big) dr \notag\\ 
&+\int_t^\tau \Big(\int_t^\tau\nabla_z \psi(r,X_r^{\sigma,X_\sigma^{t,x}},Y_r^{\sigma,X_\sigma^{t,x}},Z^{\sigma,X_\sigma^{t,x}}_r)\nabla _xZ^{\sigma,X_\sigma^{t,x}}_r (-A)^{-\alpha}\xi_\sigma\,d\sigma \Big) dr.
\end{align*}
By \eqref{YZ_epsilon_Bism} and \eqref{XYZ_dot}, together with (\ref{Xdotid_Bism}), we can conclude  that
\begin{equation}\label{YZdotid_Bism}
\overset{\cdot}{Y}_{\tau}=\int_{t}^{\tau}\nabla _xY_\tau^{\sigma,X_\sigma^{t,x}}(-A)^{-\alpha}\xi_{\sigma}d\sigma,\quad \overset{\cdot}{Z}_{\tau}=\int_{t}^{\tau}\nabla _xZ_\tau^{\sigma,X_\sigma^{t,x} }(-A)^{-\alpha}\xi_{\sigma}d\sigma, \quad \tau \in (t,T], \ \P\textup{-a.s.},
\end{equation}
since  these two pairs of processes satisfies the same BSDE. By density, arguing as in Corollary \ref{rem:identif-Z}, we infer that formulas \eqref{Xdotid_Bism} and \eqref{YZdotid_Bism} hold true for any square integrable $H$-valued predictable process $\xi$. Now, let $\eta\in E$, and let us  take
\begin{equation}\label{xi}
\xi_\tau:=(-A)^\alpha \nabla_x X_\tau^{t,x}\eta, \quad \tau \in (t,T].
\end{equation}
Notice that, since  $(-A)^\alpha \nabla _xX_\tau^{t,x}\eta\in D((-A)^{1/2-\alpha})$ $\P$-a.s., thanks to Proposition \ref{P1.9},  $(-A)^\alpha \nabla _xX_\tau^{t,x}\eta\in H$ for any $\tau \in (t,T]$, $\P$-a.s., 
and so 
\[
\int_t^\tau e^{(\tau-\sigma)A}(-A)^{-\alpha }(-A)^\alpha \nabla_x X_\sigma^{t,x}\eta\,d\sigma=\int_t^\tau e^{(\tau-\sigma)A}\nabla _xX_\sigma^{t,x}\eta\, d\sigma,\quad \tau \in (t,T], \ \P\textup{-a.s.},
\]
which belongs to $E$. Therefore,  for all $\tau \in (t,T]$ we have $\overset{\cdot}{X}_{\tau}\in E$  $\P$-a.s., where $\overset{\cdot}{X}$ denotes the mild solution to the forward equation in  (\ref{eq_dotYdotZ}) with $s-\delta =t$ with $\xi$ given by \eqref{xi}. With this choice of $\xi$ equalities (\ref{Xdotid_Bism}) and (\ref{YZdotid_Bism}) can be rewritten as
\begin{align}\label{XYZdotid_Bism}
&\overset{\cdot}{X}_{\tau}=\int_{t}^{\tau}\nabla_x X_\tau^{\sigma,X_\sigma^{t,x} }\nabla_xX_\sigma^{t,x}\eta d\sigma =(\tau-t)\nabla_xX_\tau^{t,x}\eta,\quad  
\overset{\cdot}{Y}_{\tau}=\int_{t}^{\tau}\nabla_xY_\tau^{\sigma,X_\sigma^{t,x} }\nabla_xX_\sigma^{t,x}\eta d\sigma  =(\tau-t)\nabla_xY_\tau^{t,x}\eta, \\
&\overset{\cdot}{Z}_{\tau}=\int_{t}^{\tau}\nabla_xZ_\tau^{\sigma,X_\sigma^{t,x} }\nabla_xX_\sigma^{t,x}\eta d\sigma=(\tau-t)\nabla_xZ_\tau^{t,x}\eta, \quad \tau\in[t,T], \ \P\textup{-a.s.}\nonumber
\end{align}
Let us now set 
\begin{align*}
\overset{\cdot}{\psi}(\tau,t,x) 
:=& \nabla_x \psi(\tau,X_\tau^{t,x},Y_\tau^{t,x},Z^{t,x}_\tau)\overset{\cdot}{X^{t,x}_\tau}  d\tau +\nabla_y \psi(\tau,X_\tau^{t,x},Y_\tau^{t,x},Z^{t,x}_\tau)\overset{\cdot}{Y^{t,x}_\tau}  d\tau \notag\\
 &+\nabla_z \psi(\tau,X_\tau^{t,x},Y_\tau^{t,x},Z^{t,x}_\tau)\overset{\cdot}{Z^{t,x}_\tau} d\tau,
 \\
\overset{\cdot}{\phi}(X_T^{t,x}) 
:=& \nabla_x \phi(X_T^{t,x})\overset{\cdot}{X^{t,x}_T}.  
\end{align*}
By (\ref{XYZdotid_Bism}), $\overset{\cdot}{\psi}$ and $\overset{\cdot}{\phi}$ 
can be rewritten as
\begin{align}\label{dot_psi-bis}
\overset{\cdot}{\psi}(\tau,t,x) &=(\tau-t)\left(\nabla_x \psi(\tau,X_\tau^{t,x},Y_\tau^{t,x},Z^{t,x}_\tau)\nabla_xX^{t,x}_\tau \eta d\tau+\nabla_y \psi(\tau,X_\tau^{t,x},Y_\tau^{t,x},Z^{t,x}_\tau)\nabla_xY^{t,x}_\tau\eta d\tau\right.\vspace{2mm}\notag\\
&\left. +\nabla_z \psi(\tau,X_\tau^{t,x},Y_\tau^{t,x},Z^{t,x}_\tau)\nabla_xZ^{t,x}_\tau\eta d\tau \right), 
\\
\overset{\cdot}{\phi}(X_T^{t,x}) &=(T-t)\nabla_x \phi(X_T^{t,x})\nabla_xX^{t,x}_T\eta. \label{dot_phi-bis}
\end{align}
Notice that the right-hand sides in (\ref{dot_psi-bis}) and in (\ref{dot_phi-bis}) are nothing else (modulo a renormalization) than the terms appearing in the right-hand sides of the first two equations in  (\ref{eq_gradYgradZ}).
Now we aim at finding an expression for $\overset{\cdot}{\psi}$ and $\overset{\cdot}{\phi}$ that does not involve the derivative of $\psi, \,\phi,\,X,\,Y$ and $Z$: this in turn will furnish an expression of $\nabla_xY$ that does not involve the derivatives of $\psi, \,\phi,\,X,\,Y$ and $Z$, as in formula (\ref{Bismut-lip}).
To this end, let us consider the process
\begin{equation}\label{W^e}
W_{\sigma}^{\varepsilon}=W_{\sigma}-\varepsilon%
{\int_{t}^{\sigma}}(-A)^\alpha \nabla_xX_{r}^{t,x}\eta{\it dr},\qquad 0\leq t\leq\sigma\leq T,
\end{equation}
and let us define a probability measure  $Q_{\varepsilon}$ such that 
\[
\dfrac{dQ_{\varepsilon}}{d\mathbb{P}}=\exp\Big( \varepsilon%
{\displaystyle\int_{t}^{T}}
\langle (-A)^\alpha \nabla_xX_{\sigma}^{t,x}\eta,\,
dW_{\sigma}\rangle-\dfrac{\varepsilon^{2}}{2}%
{\displaystyle\int_{t}^{T}}
\left\vert (-A)^\alpha \nabla_xX_{\sigma}^{t,x}\eta\right\vert_H ^{2}d\sigma\Big) .
\]
By the Girsanov theorem,  
under $Q_{\varepsilon}$ $(W_{\sigma}^{\varepsilon})_{\sigma\in[t,T]}$ is a 
cylindrical Wiener process in $H$.
Arguing as in the proof of Theorem \ref{prop-identific-Z},  we also notice that the process $X$ under $Q_{\varepsilon}$ and the process
$X^{\varepsilon}$
under $\mathbb{P}$ have the same law.
Therefore,  
\begin{align*}
\E[\nabla_x Y_\tau^{t,x} \eta] &=\E[\nabla_x \Phi(X_T^{t,x})\nabla_x X_T^{t,x}\eta] +\E\Big[ \int_\tau^T\big( \nabla_x \psi(\sigma,X_\sigma^{t,x},Y_\sigma^{t,x},Z^{t,x}_\sigma)\nabla_x X^{t,x}_\sigma \eta \vspace{2mm} \\
& +\nabla_y \psi(\sigma,X_\sigma^{t,x},Y_\sigma^{t,x},Z^{t,x}_\sigma)\nabla_x Y^{t,x}_\sigma \eta +\nabla_z \psi(\sigma,X_\sigma^{t,x},Y_\sigma^{t,x},Z^{t,x}_\sigma)\nabla_x Z^{t,x}_\sigma \eta\big) \, d\sigma\Big]
\vspace{2mm} \\
&=\frac{1}{T-t}\E[\overset{\cdot}{\phi}(X_T^{t,x})]+\E
 \Big[\int_\tau^T\frac{1}{\sigma-t}\overset{\cdot}{\psi}(\sigma,t,x)\,d\sigma\Big].
\end{align*}
By differentiating inside the expectation with respect to $\varepsilon$ and changing the order of integration, 
we get 
\begin{align*}
\E [\overset{\cdot}{\psi}(\sigma,t,x)]&
=\E\Big[\frac{d}{d\varepsilon}
_{\mid\varepsilon=0}\psi(\sigma, X^{\varepsilon,t,x}_\sigma,Y^{\varepsilon,t,x}_\sigma,Z^{\varepsilon,t,x}_\sigma)\Big]
=\frac{d}{d\varepsilon}
_{\mid\varepsilon=0}\E^{\Q^\varepsilon}[\psi(\sigma, X^{t,x}_\sigma,Y^{t,x}_\sigma,Z^{t,x}_\sigma)]\\
&=\E  \Big[\psi(\sigma, X^{t,x}_\sigma,Y^{t,x}_\sigma,Z^{t,x}_\sigma)\int_t^\sigma  \<(-A)^{\alpha}\nabla_x X^{t,x}_r \eta,dW_r\>\Big],
\end{align*}
and so, recalling   \eqref{def-U}, 
\begin{align*}
\E \Big[\int_\tau^T\frac{1}{\sigma-t}\overset{\cdot}{\psi}(\sigma,t,x)\,d\sigma\Big]&=\E \Big[\int_\tau^T\frac{1}{\sigma-t} \Big(\int_t^\sigma  \<(-A)^{\alpha}\nabla_x X^{t,x}_r \eta ,dW_r\> \Big) \psi(\sigma, X^{t,x}_\sigma, Y^{t,x}_\sigma, Z^{t,x}_\sigma)d\sigma\Big]\\ \vspace{2mm}
&=\E \Big[\int_\tau^TU^{\eta,t,x}_\sigma\psi(\sigma, X^{t,x}_\sigma,Y^{t,x}_\sigma,Z^{t,x}_\sigma)d\sigma\Big].
\end{align*}
Similarly,
$\E[\overset{\cdot}{\phi}(X_T^{t,x})]=\E[\phi(X_T^{t,x})U^{\eta,t,x}_T]$, and this proves \eqref{Bismut-lip} when $\eta\in E$. The general case with $\eta \in H$ follows by density, thanks to estimates \eqref{stimabismut2} and \eqref{bound-Uq}.
\endproof

In the next result we  remove the differentiability assumption on $\psi$ and $\phi$ in Theorem \ref{teoBismutLip-diffle}. 

\begin{theorem}\label{teoBismutLip}
Let Hypotheses \ref{ip_forward} and \ref{ip-psiphi-lip} hold true, and for any  $(t,x) \in [0,\,T]\times E$, let $(X^{t,x},Y^{t,x},Z^{t,x})$ be a solution of the forward-backward system  \eqref{fbsde},
 and let $U^{h,t,x}$
be the process defined in \eqref{def-U}.
Then,  for $ t\leq s\leq T$, $x\in E,\,h\in H$, the Bismut formula given in \eqref{Bismut-lip} holds true.
\end{theorem}
\proof The proof follows the same lines of the one of Theorem 3.10 in \cite{futeBismut}. The main ingredients are  formula (\ref{Id_Z}) in Theorem \ref{prop-identific-Z} and   Proposition \ref{prop-convp-bsdequadr}, which provide respectively  the identification of $Z$ in the Banach space case and with the diffusion operator $(-A)^{-\alpha}$, and the   stability result for the  BSDE in (\ref{fbsde}) when the generator and the final datum are approximated  by (\ref{psi-approx})-(\ref{phi-approx}). We underline that approximations (\ref{phi-approx})-(\ref{psi-approx}) preserve the boundedness and the growth, and are only of pointwise type. Notice that in \cite{futeBismut}, the final datum and the generator are approximated by means of their inf-sup convolutions, and so the approximation is uniform. However, thanks to  
the aforementioned stability properties for the BSDE,  our pointwise approximations (\ref{psi-approx})-(\ref{phi-approx}) are sufficient to obtain the desired result.
\endproof

\subsection{The semilinear Kolmogorov equation}\label{K_lip}

By means of Theorem  \ref{teoBismutLip}, we can give an existence and uniqueness result in the Banach framework for the semilinear Kolmogorov related to the the process $X$, and with coefficients $\phi$ and $\psi$ not necessarily differentiable, as it is assumed in \cite{MasBan}.  

\noindent Let $P_{t,\tau},\,t\leq\tau\leq T$, be the transition semigroup related to the process $X^{t,x}$ solution of the forward equation (\ref{forward}),
namely, for every bounded and measurable function $\varphi:E\rightarrow\R$,
$P_{t,\tau}[\varphi](x):=\E\,\varphi(X_\tau^{t,x})$.
We consider the following semilinear Kolmogorov equation
\begin{equation}
\left\{
\begin{array}
[c]{l}%
\frac{\partial v}{\partial t}(t,x)=-\call v\left(  t,x\right)
+\psi\left( t,x,v(t,x),\nabla^{(-A)^{-\alpha}}v(t,x)   \right)  ,\text{ \ \ \ \ }t\in\left[  0,T\right]
,\text{ }x\in E, \\
v(T,x)=\phi\left(  x\right),
\end{array}
\right.  \label{Kolmo}%
\end{equation}where
$\mathcal L $ is the generator of the transition
semigroup $(P_{t,s})_{0\leq t\leq s\leq T}$, that is, at least formally, 
$$
(\call f)(x)=\frac{1}{2}(\operatorname{Tr} ((-A)^{-\alpha}(-A^*)^{-\alpha} \nabla^2 f)(x)+\<Ax,\nabla f(x)\>+\<F(x),\nabla f(x)\>, \quad x\in E.
$$
We introduce the notion of mild solution of the nonlinear Kolmogorov
equation (\ref{Kolmo}), see e.g. \cite{fute}.
\begin{definition}
\label{defsolmildkolmo}A function $v:\left[  0,T\right]  \times E\rightarrow\mathbb{R}$ is a mild
solution of the semilinear Kolmogorov equation \eqref{Kolmo} if 
$v\in\calg^{0,1}\left(  \left[  0,T\right]  \times E\right)$, and 
\begin{equation}
v(t,x)=P_{t,T}\left[  \phi\right]  \left(  x\right)  +\int_{t}^{T}%
P_{t,s}\Big[  \psi(s,\cdot, v(s,\cdot), \nabla^{(-A)^{-\alpha}} v\left(  s,\cdot\right)) 
\Big]  \left(  x\right)  ds, 
\text{\ \ }t\in\left[  0,T\right]  ,\text{ }x\in E. \label{solmildkolmo}%
\end{equation}
\end{definition}

\begin{theorem}\label{teo_fey_kac_ip} Let Hypotheses \ref{ip_forward} 
 and \ref{ip-psiphi-lip}
hold true. 
Then the semilinear Kolmogorov equation \eqref{Kolmo} has a unique mild solution $v$ given by the formula
$$
v(t, x) = Y_t^{t,x} ,\qquad     (t, x) \in [0, T ] \times E,
$$ 
where, for any $(t,x)\in [0,\,T] \times E$,  $(X^{t,x},Y ^{t,x} , Z^{ t,x} )$ denotes the solution to the FBSDE \eqref{fbsde}. 
In addition, we
have, $\P$-a.s.,
\[
Y_s^{t,x} = v(s, X_s^{t,x} ),\quad   s \in [t,\,T], \quad      Z_s^{t,x} =\nabla_x v(s, X_s^{t,x} )\nabla_x  X_s^{t,x}(-A)^{-\alpha}, \quad \textup{a.e. }s \in [t,\,T].
\]
\end{theorem}
\proof If the data $\phi$ and $\psi$ are also differentiable,  the result can be proved  as in \cite{MasBan}, Theorem 6.2. When the data are not differentiable, the Bismut formula (\ref{Bismut-lip}) is still true, see Theorem \ref{teoBismutLip}, and the result can be proved arguing as in \cite{futeBismut}, Theorem 4.2.
\endproof

\section{The Bismut-Elworthy formula and the semilinear Kolmogorov equation: the quadratic case}
\label{sez-Bismut-quad}

We are ready to state and prove the main result of the paper, which is a nonlinear
Bismut-Elworthy formula as the one in Theorem \ref{teoBismutLip}, but in the case
of quadratic generator.   This in particular will  give an existence and uniqueness result for the Kolmogorov equation
\eqref{Kolmo}  in the quadratic case and in the Banach framework, see Theorem \ref{teo_fey_kac}. 
\begin{theorem}\label{teoBismut}
Let Hypotheses \ref{ip_forward} and 
\ref{ip-psiphi} hold true. 
For any $(t,x) \in [0,\,T] \times E$, let $(X^{t,x},Y^{t,x},Z^{t,x})$ be
the solution of the forward-backward system \eqref{fbsde} and let $U^{h,t,x}$
be the process defined in \eqref{def-U}.
Then,  for $ t\leq s\leq T$, $x\in E$ and $h\in H$, 
\begin{equation}\label{Bismut}
 \E\left[ \nabla_x\,Y^{t,x}_sh \right]=
\E\Big[\int_s^T\psi\left(r,X_r^{t,x},Y_r^{t,x},Z_r^{t,x}\right)U^{h,t,x}_r\,dr\Big]
+\E\left[  \phi(X_T^{t,x})U^{h,t,x}_T\right].
\end{equation}
\end{theorem}
 \proof 
 We split the proof into two steps: we first  prove the statement when $\psi$ is differentiable with respect to $x,y$ and $z$, and then  we remove this additional assumption.
 
 \vspace{2mm}
 \noindent {\bf STEP $1$}.
 We start by considering $\psi$ differentiable with
respect to $x$, $y$ and $z$. For all $n\geq 1$,  let us  denote by
  $(X^{t,x}, Y^{n,t,x},Z^{n,t,x})$ the solution of the Markovian BSDE in
 (\ref{fbsde}) with final datum equal to $\phi_n$ in (\ref{phi-approx}) in the place of $\phi$:
 \begin{equation}\label{bsde-n}
     \left\{\begin{array}{l}\dis
  dY^{n,t,x}_\tau=-\psi(\tau,X^{t,x}_\tau,Y^{n,t,x}_\tau,Z^{n,t,x}_\tau)\;d\tau+Z^{n,t,x}_\tau\;dW_\tau,\quad \tau \in [t,\,T], 
   \\\dis
   Y^{n,t,x}_T=\phi_n(X^{t,x}_T).
 \end{array}\right.
 \end{equation}
 By estimate (\ref{stimaZ-diffle}) in Corollary  \ref{rem:identif-Z},  for any $n\geq 1$,  there exists a constant $C(n)$, depending on $n$, which
 is bounded for every $n$ and blows up as $n\rightarrow\infty$, and  such that
 \begin{equation}\label{stimaZ-difflen}
 \vert Z_s^{n,t,x}\vert_H \leq C(n), \quad \P\textup{-a.s.}, \ \textup{a.e. }s\in[t,T]. 
 \end{equation}
 In particular, 
  $$
  \vert \psi(s,x,y,z_1)- \psi(s,x,y,z_2)\vert\leq C(n)\vert z_1-z_2\vert_H, \quad z_1,\,z_2\,\in H: \,\, \vert z_i\vert_H \leq C(n),\,i=1,2.
  $$
Therefore, the generator $\psi$ acts as a Lipschitz generator with respect to $z$ in the BSDE (\ref{bsde-n}),  so the Bismut-Elworthy formula stated in Theorem \ref{teoBismutLip} holds true for the BSDE (\ref{bsde-n}): for every $s \in [t,\,T]$, 
 \begin{equation}\label{Bismut-n}
  \E\left[ \nabla_x\,Y^{n,t,x}_sh \right]=
 \E\Big[\int_s^T\psi\left(r,X_r^{t,x},Y_r^{n,t,x},Z_r^{n,t,x}\right)U^{h,t,x}_r\,dr\Big]
 +\E\left[ \phi_n(X_T^{t,x})U^{h,t,x}_T\right].
 \end{equation}
At this point we aim at  taking the limit as $n\rightarrow\infty$ in \eqref{Bismut-n}.  
\newline We start by considering  the right-hand side of \eqref{Bismut-n}.
 By  the properties of the approximations $(\phi_n)_{n\geq1}$ together with (\ref{bound-Uq}), by the dominated convergence theorem and the pointwise convergence of $\phi_n$ to $\phi$ we have 
 \begin{align*}
   &\E\Big[\vert [ \phi_n(X_T^{t,x})-\phi(X_T^{t,x})]U^{h,t,x}_T\vert\Big]\leq 
(\E[\vert \phi_n(X_T^{t,x})-\phi(X_T^{t,x})\vert^2])^{1/2}
 (\E[\vert U^{h,t,x}_T\vert^2])^{1/2}\\
 &\leq C \left( T-t\right)^{-(1/2+\alpha)}(\E[\vert \phi_n(X_T^{t,x})-\phi(X_T^{t,x})\vert^2])^{1/2}\rightarrow 0 \,\, \textup{as } n\rightarrow\infty.
 \end{align*}
 Therefore, 
 \begin{equation*}
  \lim_{n\rightarrow\infty}\E\Big[  \phi_n(X_T^{t,x})U^{h,t,x}_T\Big]=
 \E\Big[  \phi(X_T^{t,x})U^{h,t,x}_T\Big].
 \end{equation*}
In order to compute the limit of the remaining term in the right-hand side of \eqref{Bismut-n}, 
 we will show that
 \[
  \lim_{n\rightarrow\infty}\E\Big[\int_t^T\vert\psi\left(r,X_r^{t,x},Y_r^{n,t,x},Z_r^{n,t,x}\right)U^{h,t,x}_r
 -\psi\left(r,X_r^{t,x},Y_r^{t,x},Z_r^{t,x}\right)U^{h,t,x}_r\vert\,dr\Big]=0.
 \]
 We notice that 
 \begin{align*}
  &\E\Big[\int_t^T\vert\psi\left(r,X_r^{t,x},Y_r^{n,t,x},Z_r^{n,t,x}\right)U^{h,t,x}_r
 -\psi\left(r,X_r^{t,x},Y_r^{t,x},Z_r^{t,x}\right)U^{h,t,x}_r\vert\,dr\Big]\\ \nonumber
 &=\E\Big[\int_t^{\frac{t+T}{2}}\vert\psi\left(r,X_r^{t,x},Y_r^{n,t,x},Z_r^{n,t,x}\right)U^{h,t,x}_r
 -\psi\left(r,X_r^{t,x},Y_r^{t,x},Z_r^{t,x}\right)U^{h,t,x}_r\vert\,dr\Big]\\ \nonumber
 &+\E\Big[\int_{\frac{t+T}{2}}^T\vert\psi\left(r,X_r^{t,x},Y_r^{n,t,x},Z_r^{n,t,x}\right)U^{h,t,x}_r
 -\psi\left(r,X_r^{t,x},Y_r^{t,x},Z_r^{t,x}\right)U^{h,t,x}_r\vert\,dr\Big]=:I + II.
 \end{align*}
 We start by estimating the term $I$. We have 
 \begin{align*}
  I= & \E\Big[\int_t^{\frac{t+T}{2}}\vert\psi\left(r,X_r^{t,x},Y_r^{n,t,x},Z_r^{n,t,x}\right)U^{h,t,x}_r
 -\psi\left(r,X_r^{t,x},Y_r^{t,x},Z_r^{t,x}\right)U^{h,t,x}_r\vert\,dr\Big]\\ \nonumber
  \leq & L_\psi \E\Big[\int_t^{\frac{t+T}{2}}\left(\vert Y_r^{n,t,x}-Y_r^{t,x}\vert \,\vert U^{h,t,x}_r\vert\right)\,dr \Big]\\
& +  L_\psi\E\Big[\int_t^{\frac{t+T}{2}}\left(\vert Z_r^{n,t,x}-Z_r^{t,x}\vert_H\left(1+\vert Z_r^{n,t,x}\vert_H
 +\vert Z_r^{t,x}\vert_H \right)\vert U^{h,t,x}_r\vert\right)\,dr\Big]=: 
 I_a+I_b.
 \end{align*} 
We recall that,  by  estimate (\ref{stimabismut}) in Proposition \ref{prop-aprioriZ}, and since
 $\Vert \phi_n\Vert_\infty\leq  K_\phi$, there exists a constant $C$, not depending on $n$, such that
 \begin{equation*}
 \vert Z^{n,t,x}_t\vert_H\leq C (T-t)^{-1/2}, \quad \P\textup{-a.s.}
 \end{equation*}
 So, since $Z^{n,t,x}_r=Z^{n,r,X_r^{t,x}}_r$ and $Z^{t,x}_r=Z^{r,X_r^{t,x}}_r$, for $r\in [t, \frac{t+T}{2}]$
 \begin{equation}\label{stimabismutn-bis}
 \vert Z^{n,t,x}_r\vert_H+\vert Z^{t,x}_r\vert_H\leq C \sup_{r\in [t, \frac{t+T}{2}]}(T-r)^{-1/2}\leq C(T-t)^{-1/2}, \quad \P\textup{-a.s.}
 \end{equation} 
We only show the convergence of $I_b$ since the convergence of $I_a$ follows in a simpler way by the boundedness of $Y^{t,x}$ and  of $Y^{n,t,x}$ (uniform in $n$), and by the convergence of $Y^{n,t,x}$
to $Y^{t,x}$ in $\cals^p([t,T]),\,p\geq2$. Using H\"older inequality with $p=\frac{2}{1-\beta}$ and $q=\frac{2}{1+\beta}$, for some $2\alpha<\beta<1$, together with estimate \eqref{bound-Uq-sup} in Lemma \ref{lemma:boundU-q}, we get
\begin{align*}
I_b
 &\leq C(T-t)^{-\frac{1}{2}} (T-t)^{-\frac{1}{2}\beta}\E\Big[\int_t^{\frac{t+T}{2}}\vert Z_r^{n,t,x}-Z_r^{t,x}\vert_H^{1-\beta} \vert U^{h,t,x}_r\vert\,dr\Big] \\ \nonumber
  &\leq C(T-t)^{-\frac{1}{2}} (T-t)^{-\frac{1}{2}\beta}
  \Big(\E\Big[\int_t^{\frac{t+T}{2}}\vert Z_r^{n,t,x}-Z_r^{t,x}\vert_H^2\,dr\Big]\Big)^{\frac{1-\beta}{2}}
  \Big(\int_t^{\frac{t+T}{2}}\E [\vert U^{h,t,x}_r\vert^{
  \frac{2}{1+\beta}}] \, dr \Big)^{\frac{1+\beta}{2}} \\ \nonumber
   &\leq C(T-t)^{-\frac{1}{2} (1+\beta)}  \Big(\E\Big[\int_t^{\frac{t+T}{2}}\vert Z_r^{n,t,x}-Z_r^{t,x}\vert_H^2\,dr\Big]\Big)^{\frac{1-\beta}{2}}\Big(\int_t^{\frac{t+T}{2}}
 \frac{1}{(r-t)^{\frac{1+2\alpha}{1+\beta}}}\,dr \Big)^{\frac{1+\beta}{2}}\\ \nonumber
  &\leq C(T-t)^{-\frac{1}{2}(1+\beta)} (T-t)^{(1-\frac{1+2\alpha}{1+\beta})\frac{1+\beta}{2}} 
  \Big(\E\Big[\int_t^{\frac{t+T}{2}}\vert Z_r^{n,t,x}-Z_r^{t,x}\vert_H^2\,dr\Big]\Big)^{\frac{1-\beta}{2}}
 \\ \nonumber
&\leq C(T-t)^{-\frac{1}{2}-\alpha} 
  \Big(\E\Big[\int_t^{\frac{t+T}{2}}\vert Z_r^{n,t,x}-Z_r^{t,x}\vert_H^2\,dr\Big]\Big)^{\frac{1-\beta}{2}} \rightarrow 0
 \end{align*}
  as $n\rightarrow \infty$, since $Z^{n,t,x}\rightarrow Z^{t,x}$ in $\calm^2([t,T];H)$. 
\newline Let us now  estimate the term  $II$. 
To this end, we recall that, by Theorem \ref{teo-stimap-bsdequadr},
$Y^{n,t,x}$, $Y^{t,x}$ are bounded in $\cals^p([t,T])$ and $Z^{n,t,x}$, $Z^{t,x}$ are bounded  in $\calm^{2p}([t,T];H)$,   by a constant independent
 on $n$. Moreover,  by Proposition \ref{prop-convp-bsdequadr},  $Y^{n,t,x}$ converges to $Y^{t,x}$ is $\cals^p([t,T])$ and $Z^{n,t,x}$ converges to $Z^{t,x}$ in $\calm^{2p}([t,T];H)$, for any $p \geq 1$.
 By  using again H\"older's inequality for some $p, q \geq 1$, $\frac{1}{p}+\frac{1}{q}=1$,  and estimate  (\ref{bound-Uq-sup}) in Lemma \ref{lemma:boundU-q},  we get
 \begin{align*}
 II&=\E\Big[\int_{\frac{t+T}{2}}^T\vert\psi\left(r,X_r^{t,x},Y_r^{n,t,x},Z_r^{n,t,x}\right)U^{h,t,x}_r
 -\psi\left(r,X_r^{t,x},Y_r^{t,x},Z_r^{t,x}\right)U^{h,t,x}_r\vert\,dr\Big]\\
 &\leq\E\Big[\sup_{s\in[\frac{t+T}{2},T]}\vert U^{h,t,x}_s\vert \int_{\frac{t+T}{2}}^T\vert\psi\left(r,X_r^{t,x},Y_r^{n,t,x},Z_r^{n,t,x}\right)
 -\psi\left(r,X_r^{t,x},Y_r^{t,x},Z_r^{t,x}\right)\vert\,dr\Big]\\
  &\leq\Big(\E\sup_{s\in[\frac{t+T}{2},T]}\vert U^{h,t,x}_s\vert^q\Big)^{\frac{1}{q}}
 \Big(\E\Big[\int_{\frac{t+T}{2}}^T \vert\psi\left(r,X_r^{t,x},Y_r^{n,t,x},Z_r^{n,t,x}\right)
  -\psi\left(r,X_r^{t,x},Y_r^{t,x},Z_r^{t,x}\right)\vert\,dr\Big]^p\Big)^{\frac{1}{p}}\\
 &\leq C\dfrac{1}{(T-t)^{\frac{1}{2}+\alpha}}\Big(\E\Big[\int_{\frac{t+T}{2}}^T(\vert Y_r^{n,t,x}-Y_r^{t,x}\vert+\vert Z_r^{n,t,x}-Z_r^{t,x}\vert_{H}(1+\vert Z_r^{n,t,x}\vert_{H}+\vert Z_r^{t,x}\vert_{H}))\,dr\Big]^p\Big)^{\frac{1}{p}}\\
  &\leq C\dfrac{1}{(T-t)^{\frac{1}{2}+\alpha}}\Big\{\dfrac{T-t}{2}
 \, \Big(\E\Big[\sup_{r\in[t,T]}\vert Y_r^{n,t,x}-Y_r^{t,x}\vert^{p}\Big]\Big)^{\frac{1}{p}}\\ 
 &+\Big(\E\Big[\Big(\int_{\frac{t+T}{2}}^T
  \vert Z_r^{n,t,x}-Z_r^{t,x}\vert_{H}^2\,dr\Big)^{\frac{p}{2}}
  \Big(\int_{\frac{t+T}{2}}^T\Big(1+\vert Z_r^{n,t,x}\vert_{H}+\vert Z_r^{t,x}\vert_{H} \Big)^2\,dr\Big)^{\frac{p}{2}}\Big]\Big)^{\frac{1}{p}}\Big\}\\
 &\leq C\dfrac{1}{(T-t)^{\frac{1}{2}+\alpha}}\Big\{\dfrac{T-t}{2}
\Big(\E\Big[\sup_{r\in[t,T]}\vert Y_r^{n,t,x}-Y_r^{t,x}\vert^{p}\Big]\Big)^{\frac{1}{p}} \\
 &+\Big(\E\Big[\int_{\frac{t+T}{2}}^T
  \vert Z_r^{n,t,x}-Z_r^{t,x}\vert_{H}^2\,dr\Big]^{p}\Big)^{\frac{1}{2p}}\Big(\E
  \Big[\int_{\frac{t+T}{2}}^T\Big(1+\vert Z_r^{n,t,x}\vert_{H}+\vert Z_r^{t,x}\vert_{H} \Big)^2\,dr\Big]^{p}\Big)^{\frac{1}{2p}}\Big\}\rightarrow 0,
 \end{align*}
 as $n\rightarrow \infty$. 
Collecting all the previous results, we deduce that, for every $s\in[t,T]$,
 \begin{equation}\label{eq}
 \lim_{n\rightarrow \infty} \E\left[ \nabla_x\,Y^{n,t,x}_sh \right]=
 \E\Big[\int_s^T\psi\left(r,X_r^{t,x},Y_r^{t,x},Z_r^{t,x}\right)U^{h,t,x}_r\,dr\Big]
 +\E\left[ \phi(X_T^{t,x})U^{h,t,x}_T\right].
 \end{equation}
 In particular, by taking $s=t$ in (\ref{eq}),
 \begin{equation*}
 \lim_{n\rightarrow \infty}  \nabla_x\,Y^{n,t,x}_th =
 \E\Big[\int_t^T\psi\left(r,X_r^{t,x},Y_r^{t,x},Z_r^{t,x}\right)U^{h,t,x}_r\,dr\Big]
 +\E\left[ \phi(X_T^{t,x})U^{h,t,x}_T\right], 
 \end{equation*}
which shows that  $\lim_{n\rightarrow \infty}  \nabla_x\,Y^{n,t,x}_t h$ exists. Moreover, 
 arguing as in the end of the proof of Theorem 4.1 in \cite{Mas-spa}, we deduce that
 $\lim_{n\rightarrow \infty}  \nabla_x\,Y^{n,t,x}_th = \nabla_x\,Y^{t,x}_t h$  for all $h \in H$.

\vspace{2mm}
\noindent {\bf STEP {$2$}}. Let us now  remove the differentiability assumptions on $\psi$.
For any $k\geq 1$, let $\psi_k$ be the function defined in \eqref{psi-approx}.
From Lemma \ref{L:3.4} 
we know that
 $\psi_k$ is differentiable and  it preserves the Lipschitz constant, so that
\[
 \vert \nabla_x \psi_k \vert_{E^*} \leq L_\psi,\quad\vert \nabla_y \psi_k \vert \leq L_\psi.
\]
Moreover, from Lemma \ref{L:3.4} we have $\psi_k(t,x,y,z)\rightarrow\psi(t,x,y,z)$ as $k\rightarrow+\infty$ for any $(t,x,y,z)\in [0,T]\times E\times \R\times H$, and for any $t\in[0,T],\, x\in E,\, y\in \R,\,z_1,z_2\in H,$
\begin{align}
\label{stimapsi_kpsi}
|\psi_k(t,x,y_1,z_1)-\psi_k(t,x,y_2,z_2)|\leq L_\psi(|y_1-y_2|+|z_1-z_2|_H(1+|z_1|_H+|z_2|_H)), 
\end{align}
for any $k\in\N$.
We consider the BSDE with generator equal to $\psi_k$ in the place of $\psi$:
 \begin{equation}\label{bsde-k}
     \left\{\begin{array}{l}\dis
  dY^{k,t,x}_\tau=-\psi_k(\tau,X^{t,x}_\tau,Y^{k,t,x}_\tau,Z^{k,t,x}_\tau)\;d\tau+Z^{k,t,x}_\tau\;dW_\tau, \quad \tau\in[t,T], 
   \\\dis
   Y^{k,t,x}_T=\phi(X^{t,x}_T).
 \end{array}\right.
 \end{equation}
By the first part of the proof, for any $k\geq 1$,  
\begin{equation}\label{Bismut-k}
 \E\left[ \nabla_x\,Y^{k,t,x}_sh \right]=
\E\Big[\int_s^T\psi_k\left(r,X_r^{t,x},Y_r^{k,t,x},Z_r^{k,t,x}\right)U^{h,t,x}_r\,dr\Big]
+\E\left[ \phi(X_T^{t,x})U^{h,t,x}_T\right].
\end{equation}
We aim at  taking the limit as $k\rightarrow\infty.$ 
We start by considering   the first term in the right-hand side of \eqref{Bismut-k}, and we will show that
\[
 \lim_{k\rightarrow\infty}\E\Big[\int_t^T\vert\psi_k\left(r,X_r^{t,x},Y_r^{k,t,x},Z_r^{k,t,x}\right)U^{h,t,x}_r
-\psi\left(r,X_r^{t,x},Y_r^{t,x},Z_r^{t,x}\right)U^{h,t,x}_r\vert\,dr\Big]=0.
\]
We start by splitting the integral above as follows:
\begin{align*}
 &\E\Big[\int_t^T\vert\psi_k\left(r,X_r^{t,x},Y_r^{k,t,x},Z_r^{k,t,x}\right)U^{h,t,x}_r
-\psi\left(r,X_r^{t,x},Y_r^{t,x},Z_r^{t,x}\right)U^{h,t,x}_r\vert\,dr\Big]\\ \nonumber
&=\E\Big[\int_t^{\frac{t+T}{2}}\vert\psi_k\left(r,X_r^{t,x},Y_r^{k,t,x},Z_r^{k,t,x}\right)U^{h,t,x}_r
-\psi\left(r,X_r^{t,x},Y_r^{t,x},Z_r^{t,x}\right)U^{h,t,x}_r\vert\,dr\Big]\\ \nonumber
&+\E\Big[\int_{\frac{t+T}{2}}^T\vert\psi_k\left(r,X_r^{t,x},Y_r^{k,t,x},Z_r^{k,t,x}\right)U^{h,t,x}_r
-\psi\left(r,X_r^{t,x},Y_r^{t,x},Z_r^{t,x}\right)U^{h,t,x}_r\vert\,dr\Big]=:I+II.
\end{align*}
In order to estimate the term  $I$, we notice that
\begin{align*}
I&\leq \E\Big[\int_t^{\frac{t+T}{2}}\vert\psi_k\left(r,X_r^{t,x},Y_r^{k,t,x},Z_r^{k,t,x}\right)U^{h,t,x}_r
-\psi_k\left(r,X_r^{t,x},Y_r^{t,x},Z_r^{t,x}\right)U^{h,t,x}_r\vert\,dr\Big]\\
&+\E\Big[\int_t^{\frac{t+T}{2}}\vert\psi_k\left(r,X_r^{t,x},Y_r^{t,x},Z_r^{t,x}\right)U^{h,t,x}_r
-\psi\left(r,X_r^{t,x},Y_r^{t,x},Z_r^{t,x}\right)U^{h,t,x}_r\vert\,dr\Big]=:I_a+I_b.
\end{align*}
Concerning $I_a$, by \eqref{appr_psil_phin}-\eqref{stime_appr_2} we can  argue as for $I$ in Step $1$,  and get   that 
$I_a\rightarrow 0 \,\,\textup{ as } \,\,k\rightarrow+\infty$.
\\
\noindent Let us now consider the term $I_b$. From Hypothesis \ref{ip-psiphi} and formulas \eqref{appr_psil_phin} and \eqref{stime_appr_2} it follows 
\begin{align}
\label{stima_psik_psi}
 |\psi_k\left(r,X_r^{t,x},Y_r^{t,x},Z_r^{t,x}\right)
-\psi\left(r,X_r^{t,x},Y_r^{t,x},Z_r^{t,x}\right)|
\leq & C(1+|Y_r^{t,x}|+|Z_r^{t,x}|^2_H),
\end{align}
where $C$ is a positive constant depending on $L_\psi$ and $K_\psi$. Arguing as for $I$ is Step $1$, it is possible to prove that
$$
r\mapsto(1+|Y_r^{t,x}|+|Z_r^{t,x}|^2_H)|U^{h,t,x}_r| \in L^1\Big(\Omega;L^1\Big(t,\frac{t+T}{2};\R\Big)\Big).
$$
On the other hand, recalling that $\psi_k\rightarrow \psi$ pointwise as $k\rightarrow+\infty$,  we get that 
$I_b\rightarrow0 \,\, \textup{as } \,\,k\rightarrow+\infty$
by the dominated convergence theorem.

\noindent Let us now  estimate  $II$. To this end,  we notice that 
\begin{align*}
II
\leq & \E\Big[\int_{\frac{t+T}{2}}^T\vert\psi_k\left(r,X_r^{t,x},Y_r^{k,t,x},Z_r^{k,t,x}\right)U^{h,t,x}_r
-\psi_k\left(r,X_r^{t,x},Y_r^{t,x},Z_r^{t,x}\right)U^{h,t,x}_r\vert\,dr\Big] \\
& + \E\Big[\int_{\frac{t+T}{2}}^T\vert\psi_k\left(r,X_r^{t,x},Y_r^{t,x},Z_r^{t,x}\right)U^{h,t,x}_r
-\psi\left(r,X_r^{t,x},Y_r^{t,x},Z_r^{t,x}\right)U^{h,t,x}_r\vert\,dr\Big]=:II_a+II_b.
\end{align*}
Arguing as for the term  $II$ in Step $1$,  we deduce that $II_a\rightarrow0$ as $k\rightarrow0$. 
As far as $II_b$ is considered, we get
\begin{align*}
II_b
\leq &\Big(\E\Big[\sup_{s\in[\frac{t+T}2,T]}|U_s^{h,t,x}|^q\Big]\Big)^{1/q}\Big(\E\Big[\int_{\frac{t+T}2}^T|\psi_k\left(r,X_r^{t,x},Y_r^{t,x},Z_r^{t,x}\right)
-\psi\left(r,X_r^{t,x},Y_r^{t,x},Z_r^{t,x}\right)|dr\Big]^p\Big)^{1/p}.
\end{align*}
Arguing as for $II$ in Step $1$ it follows that 
\begin{align*}
r\mapsto(1+|Y_r^{t,x}|+|Z_r^{t,x}|^2_H)|U^{h,t,x}_r| \in L^p\Big(\Omega;L^1\Big(\frac{t+T}{2},T;\R\Big)\Big).
\end{align*}
Since $\psi_k$ pointwise converges to $\psi$, we can again apply the dominated convergence theorem which gives $II_b\rightarrow0$ as $k\rightarrow+\infty$.
We can thus conclude  that, for every $s\in[t,T]$,
\begin{equation*}
\lim_{k\rightarrow \infty} \E\left[ \nabla_x\,Y^{k,t,x}_sh \right]=
\E\Big[\int_s^T\psi\left(r,X_r^{t,x},Y_r^{t,x},Z_r^{t,x}\right)U^{h,t,x}_r\,dr\Big]
+\E\left[ \phi(X_T^{t,x})U^{h,t,x}_T\right].
\end{equation*}
As in the end of Step $1$,  arguing as at the end of Theorem 4.1 in \cite{Mas-spa} we can show that,  for any $s\in[t,T]$, 
$
 \lim_{k\rightarrow \infty} \E\left[ \nabla_x\,Y^{k,t,x}_sh \right]
= \E\left[ \nabla_x\,Y^{t,x}_sh \right]$.
\endproof
We now state two corollaries: the former is about integral  estimates of $\nabla_x\,Y^{t,x}$, the latter is about the identification of $\nabla_x\,Y^{t,x}$ with $Z^{t,x}$  without differentiability assumptions.  
 Notice that, by means of the Bismut formula  \eqref{Bismut}, we can also recover estimate \eqref{stimabismut2} on $\nabla_x Y^{t,x}$. 

 \begin{corollary}\label{cor-stima-nablaY}
Let $(t,x) \in [0,\,T] \times E$.  Under the assumptions of Theorem \ref{teoBismut},
the process $\nabla_x\,Y^{t,x}$  belongs to 
$\calm^2([t,T])$, and there exists a constant $C$ depending only on $L_\psi,\,K_\psi,\,K_\phi$ such that
\begin{equation}\label{stima-nablaY-square}
\E\Big[\int_t^T\vert \nabla_x\,Y^{t,x}_s\vert^2\, ds\Big] \leq C(T-t)^{-2\alpha}.
\end{equation}
\end{corollary}
\proof
Integrating \eqref{Bismut} between $t$ and $T$ we get
 \begin{align*}
&\int_t^T\vert \E [\nabla_x\,Y^{t,x}_s]\vert^2\, ds 
=\int_t^T\Big\vert\E\Big[\int_s^T\psi\left(r,X_r^{t,x},Y_r^{t,x},Z_r^{t,x}\right)U^{h,t,x}_r\,dr\big]
+\E\left[  \phi(X_T^{t,x})U^{h,t,x}_T\right]\Big\vert^2\,ds \\ \nonumber
&\leq C\int_t^T\Big\vert\E\Big[\int_s^T\psi\left(r,X_r^{t,x},Y_r^{t,x},Z_r^{t,x}\right)U^{h,t,x}_r\,dr\Big]\Big\vert^2\,ds
+\int_t^T\Big\vert\E\left[\phi(X_T^{t,x})U^{h,t,x}_T\right]\Big\vert^2\,ds=:I+II.
 \end{align*}
We have
\begin{equation*}
 II\leq \int_t^T \Vert \phi\Vert_\infty\frac{1}{(T-t)^{1+2\alpha}}\,dr=C(T-t)^{-2\alpha}.
 \end{equation*}
For what concerns $I$, we split it as 
\begin{align*}
 I
 &=  C\Big(\int_{\frac{t+T}{2}}^T\Big\vert\E\Big[\int_s^T\psi\left(r,X_r^{t,x},Y_r^{t,x},Z_r^{t,x}\right)U^{h,t,x}_r\,dr\Big]\Big\vert^2ds \\
 & +  \int_t^\frac{t+T}{2}\Big\vert\E\Big[\int_s^T\psi\left(r,X_r^{t,x},Y_r^{t,x},Z_r^{t,x}\right)U^{h,t,x}_r\,dr\Big]\Big\vert^2ds\Big)= :I_a+I_b.
\end{align*}
From \eqref{bound-Uq-sup} and Proposition \ref{teo-stimap-bsdequadr} we have
\begin{align*}
 I_a&\leq C\int_{\frac{t+T}2}^T\Big(\E\Big[\sup_{r\in[\frac{t+T}{2},T]}\vert U^{h,t,x}_r\vert
 \int_s^T\vert\psi\left(r,X_r^{t,x},Y_r^{t,x},Z_r^{t,x}\right)\vert\,dr\Big]\Big)^2 ds\\
 &\leq C\int_{\frac{t+T}2}^T\E\Big[\sup_{r\in[\frac{t+T}{2},T]}\vert U^{h,t,x}_r\vert^2\Big]
 \E\Big[\Big(\int_s^T\left(1+\vert Y_r^{t,x}\vert+\vert Z_r^{t,x}\vert_{H}^2\right)\,dr\Big)^2\Big] ds \\ 
 &\leq C(T-t)^{-1-2\alpha}\int_{\frac{t+T}2}^T
 \E\Big[\Big(\int_s^T\left(1+\vert Y_r^{t,x}\vert+\vert Z_r^{t,x}\vert_{H}^2\right)\,dr\Big)^2\Big]ds \leq C(T-t)^{-2\alpha}.
  \end{align*}
On the other hand, we consider the function under the integral sign in $I_b$ and we split it as follows:
\begin{align*}
&\E\Big[\int_s^T\psi\left(r,X_r^{t,x},Y_r^{t,x},Z_r^{t,x}\right)U^{h,t,x}_r\,dr\Big]\\
&=  \E\Big[\int_s^{\frac{t+T}{2}}\psi\left(r,X_r^{t,x},Y_r^{t,x},Z_r^{t,x}\right)U^{h,t,x}_r\,dr\Big]  +\E\Big[\int_{\frac{t+T}2}^T\psi\left(r,X_r^{t,x},Y_r^{t,x},Z_r^{t,x}\right)U^{h,t,x}_r\,dr\Big] \\
&= : I_b'+I_b''.
\end{align*}
We argue as in the proof of Theorem \ref{teoBismut},  Step $1$.
In particular,  
  arguing as for the estimate of $I$ 
   we infer that $|I_b'|\leq C(T-t)^{-1/2-\alpha}$ for some positive constant $C$. 
On the other hand, as far as $I_b''$ is considered, arguing as in the estimate of $II$, 
we get that  $|I_b''|\leq C(T-t)^{-1/2-\alpha}$ for some positive constant $C$. 
Hence,
\begin{align*}
I_b
\leq C\int_t^{\frac{t+T}2}(T-t)^{-1-2\alpha}ds=C(T-t)^{-2\alpha},
\end{align*}
  and this concludes the proof.
 \endproof
In the following we prove that the identification of $Z$ with the directional derivative of $Y$ remains true  also when $\phi$ and $\psi$ are not differentiable.
\begin{corollary}\label{cor-identif-Z}
  Under the assumptions of Theorem \ref{teoBismut}, for every $(t, x)\in[0,T]\times  E$, 
\begin{equation}\label{identif-Z}
  Z^{t,x}_t=\nabla_x\,Y^{t,x}_t (-A)^{-\alpha}.
\end{equation}
\end{corollary}
 \proof Let $\phi$ and $\psi$ be respectively approximated by  $\phi_n$ and $\psi_n$  in \eqref{phi-approx} and \eqref{psi-approx}, and let $(Y^{n,t,x},Z^{n,t,x})$ be the solution of the BSDE  with final datum $\phi_n$ and generator $\psi_n$.
 By Theorem \ref{teo_fey_kac_ip} we already know that 
 $Z_t^{n,t,x} =\nabla_x Y^{n,t,x}_t (-A)^{-\alpha}$.
 On the other hand, we have shown in Theorem \ref{teoBismut}
 that $x\mapsto Y_\tau^{t,x} = v(\tau, X_\tau^{t,x})$ is differentiable and that $\nabla_xY^{n,t,x}_\tau\rightarrow \nabla_x Y^{t,x}_\tau$,
  $dt\otimes d\P$ a.e. and a.s., as $n\rightarrow \infty$. Moreover, by computing the joint quadratic variation between the process
 $v^n(\tau,X_\tau^{t,x}):=Y^{n,t,x}_\tau,\,t\leq\tau\leq T$, and $\int_t^\cdot \xi_s\,dW_s,\,\xi\in 
 \calm^2([t,T];H)$,
 it turns out that
 \[
  \int_t^\tau \nabla v^n(s,X_s^{t,x})(-A)^{-\alpha}\xi_s\,ds
 =\int_t^\tau Z^{n,t,x}_s\xi_s\,ds, \;\P\text{-a.s.}, \ \textup{a.e.}\ \tau\in[t,T].
 \]
 By taking a subsequence (that for simplicity we call again $n$) and
 letting $n\rightarrow\infty $ in both sides, from Proposition \ref{prop-convp-bsdequadr} we get
 \[
  \int_t^\tau \nabla v(s,X_s^{t,x})(-A)^{-\alpha}\xi_s\,ds
 =\int_t^\tau Z^{t,x}_s\xi_s\,ds, \;  \textup{a.e.}\ \tau\in[t,T],\ \P\text{-a.s.},
 \]
 which gives
 formula \eqref{identif-Z}. 
 \endproof
Using Theorem  \ref{prop-identific-Z},  we can give an existence and uniqueness result for the Kolmogorov equation
\eqref{Kolmo} and we can provide  a Feynman-Kac formula in the quadratic case and in the Banach framework.
\begin{theorem}\label{teo_fey_kac} Let Hypotheses \ref{ip_forward} 
 and \ref{ip-psiphi}
hold true.
Then there  exists a unique mild solution $v(t,x)$ of the semilinear Kolmogorov equation \eqref{Kolmo}
given by the formula
\[
 v(t,x)=Y^{t,x}_t,
\]
where $(X^{t,x},Y ^{t,x} , Z^{ t,x} )$ is the solution to the FBSDE \eqref{fbsde}, 
and 
$\P$-a.s.,
\[
Y_s^{t,x} = v(s, X_s^{t,x} ), \quad       Z_s^{t,x} =\nabla_x v(s, X_s^{t,x} )\nabla_x  X_s^{t,x}(-A)^{-\alpha}, \ \textup{a.e. }s\in[t,T].
\]
In particular,
\begin{equation*}
\vert v(t,x)\vert\leq C,\qquad \vert \nabla_x v(t,x)\vert\leq C\left(T-t\right)^{-(\frac{1}{2}+\alpha)}.
\end{equation*}
If   in addition $\phi$ is G\^ateaux differentiable with bounded derivative, and $\psi$
is G\^ateaux differentiable with respect to $x$, $y$ and $z$,  then  
\begin{equation*}
\vert Z_s^{t,x}\vert_H \leq C.
\end{equation*}
\end{theorem}
\proof 
For the first part without differentiability assumptions on $\phi$ and $\psi$, it is enough to apply Theorem \ref{teoBismut} and Corollary \ref{cor-identif-Z} to get existence of the solution, as well as the estimate for $v$. The uniqueness follows from the uniqueness of the solution of the related BSDE. The estimate for $\nabla_x v(t,x)$ is a direct consequence of Proposition \ref{prop-aprioriZ}.
The second part of the result can be proved  in a standard way by means of  Proposition \ref{P:BC_Diff_BSDE} and the identification of $Z$ proved in Theorem \ref{prop-identific-Z}, see e.g. the proof of Theorem 6.2 in \cite{fute}.
\endproof

\section{A quadratic optimal control problem}\label{sez-appl-contr}
In this section we deal with the controlled state equation
\begin{align}
\label{2-sdecontrolforte}
\left\{
\begin{array}{ll}
dX_\tau^{u}=AX_\tau^{u}d\tau+F(X_\tau^{u})d\tau+Q u_\tau d\tau+(-A)^{-\alpha}dW_\tau, & \tau\in[t,T], \vspace{1mm} \\
X_t^{u}=x\in E,
\end{array}
\right.
\end{align}
where 
$Q=(-A)^{-\alpha}$ or $Q=I$, and $u$ is the control process  belonging to a suitable space $\mathcal U$ of $H$-valued functions.
 We will study the  optimal control problem associated to equation (\ref{2-sdecontrolforte}) with cost functional $J: [0,\,T] \times E \times \mathcal U\rightarrow \R$ defined by 
\begin{align}
\label{costoso}
J(t,x,u):=\mathbb E\Big[\int_t^T\ell(s,X_s^u,u_s)\, ds\Big]+\mathbb E[\Phi(X_T^u)],
\end{align}
that we are going to minimize over all admissible controls. We define  the value function of the optimal control problem as 
\begin{align}\label{V}
V(t,x):=\inf_{u\in \mathcal U}J(t,x,u), \quad x\in H, \ t\in[0,T].
\end{align}
For any $p\geq 1$, we introduce the spaces of admissible control processes
\begin{align*}
\mathcal U_p& :=\left\{u\in L^2(\Omega;L^p(0,T;H)):u \textrm{ is adapted}\right\}, \\
\mathcal U_p^\alpha & :=\left\{u\in L^2(\Omega;L^p(0,T;D((-A)^\alpha))):u\textrm{ is adapted}\right\},
\end{align*}
where $D((-A)^\alpha))$ is endowed with the norm
\[
|x|_\alpha:=|x|_H+|(-A)^\alpha x|_H.
\] 
We first prove some results about well posedness of the controlled equation (\ref{2-sdecontrolforte}). The main novelty towards Section \ref{Sec:2} and the known results in the literature is that the controls $u$ are not necessarily bounded, together with the fact that $X$ evolves in a Banach space E.

\noindent Beside Hypothesis \ref{ip_forward} we assume the following. 
\begin{hypothesis}
\label{ip_contr_E_H}
There exists $\beta>0$ such that $D((-A)^\beta)\subset E$ with continuous embedding.
\end{hypothesis}

\begin{remark}
\label{rmk_contr_E_H}
Let $A$ be an operator satisfying Hypothesis \ref{ip_forward}-(i). If Hypothesis \ref{ip_contr_E_H} holds true, then we have the following.
\begin{itemize}
\item[(i)] 
For any $t>0$ and $h\in H$,  $e^{tA}h\in E$ and there exists a positive constant $c$ such that
\begin{align}
\label{stima_E_H}
|e^{tA} h|_E\leq c t^{-\beta}|h|_H.
\end{align}
\item[(ii)] For any $t>0$ and $h\in H$, there exists a positive constant $c$ such that
\begin{align}
\label{stima_E_H2}
|e^{tA}(-A)^{-\alpha}h|_E\leq c t^{(-\beta+\alpha)\wedge 0}|h|_H.
\end{align}
\item[(iii)] For any $t>0$ and $h\in D((-A)^{\alpha})$, there exists a positive constant $c$ such that
\begin{align*}
|e^{tA}h|_E
= |e^{tA}(-A)^{-\alpha}(-A)^{\alpha}h|_E
\leq ct^{(-\beta+\alpha)\wedge0}|(-A)^\alpha h|_H
\leq ct^{(-\beta+\alpha)\wedge 0}|h|_\alpha.
\end{align*}
\end{itemize}
\end{remark}

\begin{remark}
Hypothesis \ref{ip_contr_E_H} may be replaced by the weaker condition in Remark \ref{rmk_contr_E_H}-(i).
However, this condition would not imply Remark \ref{rmk_contr_E_H}-(ii)-(iii).
\end{remark}

\begin{example}
Let $\mathcal D\subset \R^2$ be a bounded domain with smooth boundary. Set  $H=L^2(\mathcal D)$,  $E=C(\overline{\mathcal D})$, and  let $A$ be the Laplace operator with Dirichlet boundary conditions. Then, Hypothesis \ref{ip_contr_E_H} is satisfied with $\beta>1/2$.
\end{example}
We will deal with mild solutions to \eqref{2-sdecontrolforte}, namely  adapted processes $X^{t,x,u}:[t,T]\times \Omega\rightarrow E$ such that\begin{align}
\label{mild_sol_contr_syst}
X_\tau^{t,x,u}=e^{\tau A}x+\int_t^\tau e^{(\tau-s)A}F(X_s^{t,x,u})ds+\int_t^\tau e^{(\tau-s)A} Q u_sds
+\int_t^{\tau}e^{(\tau-s)A} (-A)^{-\alpha} dW_s,
\end{align}
for any $\tau\in[t,T]$, $\P$-a.s.
For any  $t\in[0,T]$, $u \in \mathcal U$, we set 
\begin{equation}\label{conv}
I^u(t,\tau):=\displaystyle\int_t^\tau e^{(\tau-s)A}Q u_s \, ds,\quad \tau \in [t,\,T].
\end{equation}
\begin{lemma}
\label{stima_E_contr}
Let $A$ be an operator satisfying Hypothesis \ref{ip_forward}-(i),  and assume that Hypothesis  \ref{ip_contr_E_H} holds true  for some positive   constant $\beta$. 
Let  $p\geq1$, and  set $p'$ be the conjugate exponent of $p$, i.e., $p^{-1}+(p')^{-1}=1$.
Then the following hold. 
\begin{itemize}
\item[(i)]
Case $Q=(-A)^{-\alpha}$ and  $p'[(\beta-\alpha)\vee0]<1$.

For any $u\in\mathcal U_p$, $I^u(t,\tau)\in E$ for any $\tau\in[t,T]$, $\P$-a.s., and there exists a positive constant $c_{\alpha,\beta,p,T}$ such that
\begin{align}\label{est_I}
|I^u(t,\tau)|_E\leq c_{\alpha,\beta,p,T}\|u\|_{L^p(0,T;H)}, \quad \tau\in[t,T], \ \P\textup{-a.s.}
\end{align}
\item[(ii)]
Case $Q=I$ and $p'[(\beta-\alpha)\vee0]<1$.

For any $u\in\mathcal U_p^\alpha$,  $I^u(t,\tau)\in E$ for any $\tau\in[t,T]$, $\P$-a.s., and there exists a positive constant $c_{\alpha,\beta,p,T}$ such that
\begin{align*}
|I^u(t,\tau)|_E\leq c_{\alpha,\beta,p,T}\|u\|_{L^p(0,T;D((-A)^\alpha))}, \quad \tau\in[t,T], \ \P\textup{-a.s.}
\end{align*}
\item[(iii)]
Case $Q=I$ and  $p'\beta<1$.

For any $u\in\mathcal U_p$,  $I^u(t,\tau)\in E$ for any $\tau\in[t,T]$, $\P$-a.s., and $I^u(t,\tau)$ satisfies estimate \eqref{est_I} for some  positive constant $c_{\alpha,\beta,p,T}$.
\end{itemize}
\end{lemma}
\proof
Let us prove item $(i)$, items $(ii)$ and $(iii)$ follow from similar arguments. From Hypothesis \ref{ip_contr_E_H},  we have
\begin{align*}
|e^{(\tau-s)A}(-A)^{-\alpha}u_s|_E\leq c(\tau-s)^{(-\beta+\alpha)\wedge 0}|u_s|_H \quad \textup{a.e.}\,\, s\in(t,\tau), \P\textup{-a.s.}
\end{align*}
Therefore,
\begin{align*}
\Big|\int_t^\tau e^{(\tau-s)A} (-A)^{-\alpha}u_sds\Big|_E
\leq & \int_t^\tau \Big|e^{(\tau-s)A}(-A)^{-\alpha} u_s\Big|_E ds
\leq c\int_t^\tau(\tau-s)^{(-\beta+\alpha)\wedge 0}|u_s|_H ds \\
\leq &  c\Big(\int_t^\tau(\tau-s)^{[(-\beta+\alpha)\wedge 0] p'}ds\Big)^{1/p'}\|u\|_{L^p(0,T;H)} \\
\leq &c(\tau-t)^{(-\beta+\alpha)\wedge0+1/p'}\|u\|_{L^p(0,T;H)}, \quad \P\textup{-a.s.}
\end{align*}
\endproof

Thanks to Lemma \ref{stima_E_contr}, arguing as in \cite[Theorem 7.11]{dpz92} we deduce the following result, which is the counterpart of Proposition \ref{P:unificata}-(i) for the controlled equation.
\begin{proposition}
\label{prop:mild_sol_contr_prob_E}
Let  Hypothesis \ref{ip_forward} holds true,  and assume that Hypothesis  \ref{ip_contr_E_H} holds true  for some positive   constant $\beta$. 
Let  $t\in[0,T]$, $p\geq1$, and  set $p'$ be the conjugate exponent of $p$.
Then the following hold.
\begin{itemize}
\item[(i)] Case $Q=(-A)^{-\alpha}$,  $p'[(\beta-\alpha)\vee0]<1$.

For any $x\in E$ and  $u\in\mathcal U_p$,  there exists a unique mild solution $X_\tau^{t,x,u}$ to  \eqref{2-sdecontrolforte}  belonging to 
$\cals^2((t,T];E)$. Moreover, there exists a positive constant $c$ such that, for any $\tau\in[t,T]$,
\begin{align}
\label{stima_sol_contr_gen}
|X_\tau^{t,x,u}|_E\leq c \Big(|x|_E+\|u\|_{L^p(t,T;H)}^{2m+1}+\sup_{\tau\in[t,T]}|w^A(t,\tau)|_E^{2m+1}\Big), \quad \P\textup{-a.s.}
\end{align}
\item[(ii)]Case $Q=I$, $p'[(\beta-\alpha)\vee0]<1$. 
 
 For any $x\in E$ and  $u\in\mathcal U_p^\alpha$, there exists a unique mild solution $X_\tau^{t,x,u}$ to  \eqref{2-sdecontrolforte}  belonging to 
 $\cals^2((t,T];E)$. Moreover, 
there exists a positive constant $c$ such that, for any $\tau\in[t,T]$,
\begin{align}
\label{stima_sol_contr_gen2}
|X_\tau^{t,x,u}|_E\leq c \Big(|x|_E+\|u\|_{L^p(t,T;D((-A)^\alpha))}^{2m+1}+\sup_{\tau\in[t,T]}|w^A(t,\tau)|_E^{2m+1}\Big), \quad \P\textup{-a.s.}
\end{align}
\item[(iii)] Case $Q=I$,  $p'\beta<1$. 

For any $x\in E$ and $u\in\mathcal U_p$, there exists a unique mild solution $X_\tau^{t,x,u}$ to  \eqref{2-sdecontrolforte}  belonging to
$\cals^2((t,T];E)$. Moreover, there exists a positive constant $c$ such that,  for any $\tau\in[t,T]$, 
\begin{align}
\label{stima_sol_contr_gen3}
|X_\tau^{t,x,u}|_E\leq c \Big(|x|_E+\|u\|_{L^p(t,T;H)}^{2m+1}+\sup_{\tau\in[t,T]}|w^A(t,\tau)|_E^{2m+1}\Big), \quad \P\textup{-a.s.}
\end{align}
\end{itemize}
\end{proposition}
\proof
We show item $(i)$, the proof of items $(ii)$ and $(iii)$ being analogous.  Since by Lemma \ref{stima_E_contr} the convolution defined in 
\eqref{conv} is a well defined $E$-valued process for any $u\in \mathcal U_p$, it is possible to argue as in \cite[Theorem 7.11]{dpz92}. Therefore, by applying the fixed point theorem we infer that for any $t\in[0,T]$,  $x\in E$ and  $u\in \mathcal U_p$, there exists a unique mild solution $X^{\alpha,t,x,u}$ to \eqref{2-sdecontrolforte} with $F$ replaced by its Yosida approximations $F_\alpha$, $\alpha>0$, such that $X^{\alpha,t,x,u}$ satisfies \eqref{stima_sol_contr_gen}. Further, the sequence $\{X^{\alpha,t,x,u}\}_{\alpha>0}$ converges as $\alpha\rightarrow0$ to the mild solution $X^{t,x,u}$ to \eqref{2-sdecontrolforte}. In particular,  estimate \eqref{stima_sol_contr_gen} holds true also for $X^{t,x,u}$.
\endproof

\subsection{The structure condition: the case \texorpdfstring{$Q=(-A)^{-\alpha}$}{Q=}}
\label{subsec_str_cond}
In this section we deal with control processes  $u \in \mathcal U_2$, and with  
the controlled equation 
\begin{align}
\label{2-sdecontrolforte_str_cond}
\left\{
\begin{array}{ll}
dX_\tau^{u}=AX_\tau^{u}d\tau+F(X_\tau^{u})d\tau+(-A)^{-\alpha}u_\tau d\tau+(-A)^{-\alpha}dW_\tau, & \tau\in[t,T], \vspace{1mm} \\
X_t^{u}=x\in E,
\end{array}
\right.
\end{align}
 satisfying the so called structure condition: the control affects the system only through the noise. 

We make the following assumptions on the cost functional \eqref{costoso}. 
\begin{hypothesis}\label{H: control1}
Let $\phi:E\rightarrow \R$ and $\ell:[0,T]\times E\times H\rightarrow \R$  be two measurable functions satisfying the following properties.
\begin{description}
\item[(i)] $\phi$ is continuous and bounded.
\item[(ii)] 
 For all $t\in[0,T]$, $u\in H$, the function $x\mapsto \ell(t,x,u)$ is bounded and continuous from $E$ onto $\R$.  
 For all $t\in[0,T]$, $x\in E$, the function $u\mapsto \ell(t,x,u)$ is continuous from $H$ onto $\R$. Further, there exist  $c, C,R$ positive constants such that, for all $t\in[0,T]$, $x\in E$, $u\in H$,
\begin{align}
&0\leq \ell(t,x,u)\leq c(1+|u|_H)^2, \label{2belgio}\\
&\ell(t,x,u)\geq C|u|^2_{H}, \quad   |u|_{H}\geq R. \label{2francia}
\end{align}
\item[(iii)] 
 There exists a positive constant $L>0$ such that, for all  $t\in[0,T]$, $x_1,x_2\in E$, $u\in H$,
$$
|\ell(t,x_1,u)-\ell(t,x_2,u)|\leq L|x_1-x_2|_E, 
$$
\end{description}
\label{2inghilterra}
\end{hypothesis}
\begin{remark}\label{R: estl}
Under Hypothesis \ref{H: control1}-(ii), it is easy to see that there exist  $c, R$  positive constants such that  
$$
 \,\,\ell(t,x,u)\geq c(|u|_H^2-R^2), \quad t\in[0,T], \ x\in E, \ u\in H.
$$	
\end{remark}
We introduce the Hamiltonian function
\begin{align}
\label{imballaggio1}
\psi(t,x,z):=\inf_{u\in H}\left\{\ell(t,x,u)+\langle z,u\rangle_H\right\}, \quad t\in[0,T], \ x\in E, \ z\in H. 
\end{align} 
Arguing as in \cite[Lemma 3.1]{FT06} we deduce an analogous result.
\begin{lemma}\label{viaggione}
Let Hypotheses \ref{2inghilterra} be satisfied. Then, the function $\psi$ in \eqref{imballaggio1}  is Borel measurable, and there exists a positive constant $C$ such that
\begin{align}
\label{2russia}
-C(1+|z|^2_H)\leq \psi(t,x,z)\leq \ell(t,x,u)+|z|_H| u|_H, \quad t\in[0,T], \ x\in E, \ z,u \in H.
\end{align}
Further, if the minimum in \eqref{imballaggio1} is attained, it is attained in a ball of radius $C(1+|z|_H)$, i.e., 
\begin{align}
\psi(t,x,z)&=\inf_{u\in H, |u|_H\leq C(1+|z|_H)}\left\{\ell(t,x,u)+\langle z, u\rangle_H\right\},
\quad t\in[0,T], \,\, x\in E,\  z\in H, \label{stima_min_locale}\\
\psi(t,x,z)&\leq\ell(t,x,u)+\langle z,u\rangle_H, \quad |u|_H\geq C(1+|z|_H).\notag
\end{align}
Finally, there exists a positive constant $C$ such that, for any $x_1,x_2\in E$, $z_1,z_2\in H$, 
\begin{align}
\label{2podio}
|\psi(t,x_1,z_1)-\psi(t,x_2,z_2)|\leq C(|x_1-x_2|_E+|z_1-z_2|_H(1+|z_1|_H+|z_2|_H)), \quad t\in[0,T].
\end{align}
\end{lemma}
The HJB equation associated to the control problem (\ref{V}), related to the controlled state equation \eqref{2-sdecontrolforte_str_cond},
  is given by
\begin{equation}
\left\{
\begin{array}
[c]{l}%
\frac{\partial v}{\partial t}(t,x)=-\call v\left(  t,x\right)
+\psi\left( t,x,v(t,x),\nabla^{(-A)^{-\alpha}}v(t,x)   \right)  ,\text{ \ \ \ \ }t\in\left[  0,T\right]
,\text{ }x\in E, \\
v(T,x)=\phi\left(  x\right),
\end{array}
\right.  \label{Hjb-structure}%
\end{equation}
where $\psi$ is defined in (\ref{imballaggio1}).
The HJB equation (\ref{Hjb-structure}) turns out to be a semilinear Kolmogorv equation as (\ref{Kolmo}), with $\psi$ and $\phi$ satisfying Hypotehsis \ref{ip-psiphi}. So by Theorem \ref{teo_fey_kac} its    mild solution can be represented in terms of the solution $(X^{t,x}, Y^{t,x}, Z^{t,x})$ of the forward-backward system (\ref{fbsde}).

In the following Theorem we state and prove the fundamental relation, and we characterize the optimal control with a feedback law.
\begin{theorem}
\label{thm_rel_fond_opt_cont_str}
Let  Hypotheses \ref{ip_forward},  \ref{H: control1} hold true,  and assume that Hypothesis  \ref{ip_contr_E_H} holds true  with a  constant $\beta$ 
	  such that 
$\beta-\alpha<1/2$. Let  $X^{t,x,u}$ be the mild solution of \eqref{2-sdecontrolforte_str_cond},$V(t,x)$ be the value function of the control problem \eqref{V}, and $v $ be the mild solution of the HJB equation \eqref{Hjb-structure}.
Then, for any $(t, x) \in[0,T]\times  E$ and
  $u\in \mathcal U_2$, the so called fundamental relation holds true:
\begin{align*}
v(t,x)
= J(t,x,u)+\E\Big[\int_\tau^T\left(\psi(s,X_s^{t,x,u},Z_s^{t,x})-\ell(s, X_s^{t,x,u},Z_s^{t,x})-Z_s^{t,x}u_s\right)ds\Big].
\end{align*}
In particular,
$v(t,x)\leq V(t,x)$, for all $(t,x) \in[0,T]\times E$.
Moreover, if there exists a measurable function $\gamma:[0,T]\times E\times H\rightarrow H$ satisfying 
\begin{align*}
\psi(t,x,z)=\ell(t,x,\gamma(t,x,z))+\langle z,\gamma(t,x,z)\rangle_H, \quad  t\in[0,T], \,x\in E,\  z\in H,
\end{align*}
then 
$$
v(t,x)=V(t,x)
$$ and, thanks to \eqref{stima_min_locale}, the process $\bar u$ defined by
\begin{align*}
\bar u _s:= \gamma(s,X_s^{x,\bar u},\nabla_x v(s, X_s^{x, \bar u} )\nabla_x  X_s^{x, \bar u}(-A)^{-\alpha}) \quad  \textup{for} \,\textup{-a.e.} \, \, s\in(0,T), \quad \P\textup{-a.s.},
\end{align*} 
belongs to $\mathcal U_2$ and it is optimal.
\end{theorem}
\proof
The proof is standard and follows the same  lines of \cite[Proposition 4.1]{FT06}. We notice that, by Proposition \ref{prop:mild_sol_contr_prob_E}-$(i)$, problem \eqref{2-sdecontrolforte_str_cond} admits a unique mild solution $X^{t,x,u}$ for any $u\in\mathcal U_2$.  Further, for any $u\in\mathcal U_2$, we introduce the family of stopping times $\tau_n$ defined by
\begin{align*}
\tau_n:=\inf\Big\{\tau\in[t,T]:\int_t^\tau|u_s|_H^2ds>n\Big\}, \quad n\in\N.
\end{align*}
Then we  
proceed as in \cite[Proposition 4.1]{FT06}, by applying the Girsanov Theorem and using the fact that $\psi$ satisfies Hypothesis \ref{ip-psiphi}-(ii), and that the pair of processes $(Y^{t,x},Z^{t,x})$, solution to the Markovian BSDE in (\ref{fbsde}), are identified respectively with the solution $v$ of the HJB equation (\ref{Hjb-structure}) and with its directional derivative $\nabla^{(-A)^{-\alpha}}v$. Namely, by Theorems \ref{teoBismut} and \ref{teo_fey_kac}, $Y_s^{t,x}=  v(s,X_s^{t,x}) $ and  $Z_s^{t,x}= \nabla_x v(s, X_s^{x, \bar u} )\nabla_x  X_s^{x, \bar u}(-A)^{-\alpha}$.
\endproof

\subsection{The case \texorpdfstring{$Q=I$}{Q=I} with a special running cost}\label{sez-cost-alfa}
In the present section we deal with control processes  $u \in \mathcal U_2$, and with the controlled equation 
\begin{align}
\label{Q=I}
\left\{
\begin{array}{ll}
dX_\tau^{u}=AX_\tau^{u}d\tau+F(X_\tau^{u})d\tau+u_\tau d\tau+(-A)^{-\alpha}dW_\tau, & \tau\in[t,T], \vspace{1mm} \\
X_t^{u}=x\in E.
\end{array}
\right.
\end{align}
The controlled equation \eqref{Q=I} has  a different structure towards \eqref{2-sdecontrolforte_str_cond} considered in Subsection \ref{subsec_str_cond}, so the problem is different, and  we need different assumptions on the cost functional \eqref{costoso}.
\begin{hypothesis}\label{H_control2}
Let $\phi:E\rightarrow \R$ and $\ell:[0,T]\times E\times H\rightarrow \R\cup\{+\infty\}$  be two measurable functions satisfying the following properties.
\begin{description}
\item[(i)] $\phi$ is continuous and bounded.
\item[(ii)] 
 For all $t\in[0,T]$, $u\in D((-A)^\alpha)$, the function $x\mapsto \ell(t,x,u)$ is bounded and continuous from $E$ onto $\R$.  For all $t\in[0,T]$, $x\in E$, the function $u\mapsto \ell(t,x,u)$ is continuous from $D((-A)^\alpha)$ onto $\R$. Further, there exists $c, C, R$ positive constants  such that, for all $t\in[0,T]$, $x\in E$ and $u\in D((-A)^\alpha)$, \begin{align}
&0\leq \ell(t,x,u)\leq c(1+|u|_\alpha)^2,\label{3belgio}\\
&\ell(t,x,u)\geq C|u|^2_{\alpha}, \quad |u|_{\alpha}\geq R.\label{3francia}
\end{align}
\item[(iii)] 
There exists a positive constant $L>0$ such that, for any $t\in[0,T]$,  $u\in D((-A)^\alpha)$,  $x_1,x_2\in E$, 
\begin{align*}
|\ell(t,x_1,u)-\ell(t,x_2,u)|\leq L|x_1-x_2|_E.
\end{align*}
\end{description}
\label{ip_contr_pen}
\end{hypothesis}
\begin{remark}
Condition  \eqref{3francia} in Hypothesis \ref{H_control2} implies that, if $u$ does not take values in $D((-A)^\alpha)$, then $J(t,x,u)=+\infty$. In particular, 
$\inf_{u \in \mathcal U^\alpha_2} J(t,x,u)=\inf_{u \in \mathcal U_2} J(t,x,u)$, 
so we can limit ourselves to consider here the space  of admissible controls $\mathcal U^\alpha_2$.
\end{remark}
\begin{remark}
Under Hypothesis \ref{H_control2}-(ii), there exist positive constants $c, R$ such that, for any $t\in [0,T]$, $x\in E$, $u\in D((-A)^\alpha)$,  we have
$\ell(t,x,u)\geq c(|u|_\alpha^2-R^2)$. 
\end{remark}	

We introduce the Hamiltonian function
\begin{align}
\label{ham_pen_cost}
\psi^\alpha(t,x,z):=\inf_{u\in D((-A)^\alpha)}\left\{\ell(t,x,u)+\langle z,(-A)^\alpha u\rangle_H\right\}, \quad t\in[0,T], \ x\in E, \ z\in H. 
\end{align} 
Arguing again as in \cite[Lemma 3.1]{FT06}, we infer the following properties of  $\psi^\alpha$.
\begin{lemma}\label{lem_stime_ham}
Let Hypotheses \ref{ip_contr_pen} be satisfied. Then, the  function $\psi^\alpha$ in \eqref{ham_pen_cost} is Borel measurable and there exists a positive constant $C$ such that
\begin{align}
\label{3russia}
-C(1+|z|^2_H)\leq \psi^\alpha(t,x,z)\leq \ell(t,x,u)+|z|_H| u|_\alpha, \quad t\in[0,T], \ x\in E, \ z\in H,\ u \in D((-A)^\alpha). 
\end{align}
Further, if the minimum in \eqref{ham_pen_cost} is attained, it is attained in a ball of radius $C(1+|z|_H)$, i.e., 
\begin{align}
\label{stima_min_loc_pen_cost}
\psi^\alpha(t,x,z)=\inf_{u\in D((-A)^\alpha), |u|_\alpha\leq C(1+|z|_H)}\left\{\ell(t,x,u)+z(-A)^\alpha u\right\},
\quad t\in[0,T], \,\, x\in E,\  z\in H.
\end{align}
Finally, for any $x_1,x_2\in E$, $z_1,z_2\in H$, $\psi^\alpha$, there exists  a  positive constant $C$ such that
\begin{align}
\label{lipz_ham_pen_cost}
|\psi^\alpha(t,x_1,z_1)-\psi^\alpha(t,x_2,z_2)|\leq C(|x_1-x_2|_E+ |z_1-z_2|_H (1+|z_1|_H+|z_2|_H)), \quad t\in[0,T].
\end{align}
\end{lemma}
The HJB equation associated to the control problem (\ref{V}), related to the controlled state equation \eqref{Q=I},
  is given by
\begin{equation}
\left\{
\begin{array}
[c]{l}%
\frac{\partial v}{\partial t}(t,x)=-\call v\left(  t,x\right)
+\psi^\alpha\left( t,x,v(t,x),\nabla^{(-A)^{-\alpha}}v(t,x)   \right)  ,\text{ \ \ \ \ }t\in\left[  0,T\right]
,\text{ }x\in E, \\
v(T,x)=\phi\left(  x\right),
\end{array}
\right.  \label{Hjb-alfa}%
\end{equation}
where $\psi^\alpha$ is defined in (\ref{ham_pen_cost}).
Again, the HJB equation (\ref{Hjb-alfa}) turns out to be a semilinear Kolmogorv equation as (\ref{Kolmo}), with $\psi^\alpha$ and $\phi$ satisfying Hypotehsis \ref{ip-psiphi}. So by Theorem \ref{teo_fey_kac} its mild solution can be represented in terms of the solution $(X^{t,x}, Y^{t,x}, Z^{t,x})$ of the forward-backward system
\begin{equation}\label{fbsde-alfa}
    \left\{\begin{array}{ll}
\dis dX_\tau =
AX_\tau d\tau+ F(X_\tau)d\tau+(-A)^{-\alpha}dW_\tau,& \tau\in
[t,T],
\\\dis
X_t=x,
\\\dis
 dY_\tau=-\psi^\alpha(\tau,X_\tau,Y_\tau,Z_\tau)\;d\tau+Z_\tau\;dW_\tau, & \tau\in[t,T],
  \\\dis
  Y_T=\phi(X_T),
\end{array}\right.
\end{equation} which is nothing else than the forward-backward system (\ref{fbsde}) with $\psi$ instead of $\psi^\alpha$.

As in Subsection \ref{subsec_str_cond}, in the following Theorem we state and prove the fundamental relation, and we characterize the optimal control with a feedback law.
\begin{theorem}
\label{thm_rel_fond_opt_cont_spec_cost}
Let  Hypotheses \ref{ip_forward},  \ref{H_control2}  hold true, and assume that Hypothesis \ref{ip_contr_E_H} holds true with a  constant $\beta$ 
	  such that 
$\beta-\alpha<1/2$. Let  $X^{t,x,u}$ be the mild solution of \eqref{Q=I}, $V(t,x)$ be the value function of the control problem \eqref{V}, and $v $ be the mild solution of the HJB equation \eqref{Hjb-alfa}.
Then, for any  $(t,x)\in[0,T]\times  E$ and 
 $u\in \mathcal U_2^\alpha$, \begin{align*}
v(t,x)
= J(t,x,u)+\E\int_\tau^T\left(\psi^\alpha(s,X_s^{t,x,u},Z_s^{t,x})-\ell(s, X_s^{t,x,u},Z_s^{t,x})-Z_s^{t,x}(-A)^\alpha u_s\right)ds.
\end{align*}
In particular, 
$v(t,x)\leq V(t,x)$, for all $t\in[0,T]$, $x\in E$.
Moreover, if there exists a measurable function $\gamma^\alpha:[0,T]\times E\times H\rightarrow D((-A)^\alpha)$ satisfying 
\begin{align*}
\psi^\alpha(t,x,z)=\ell(t,x,\gamma^\alpha(t,x,z))+\langle z,(-A)^{\alpha}\gamma^\alpha(t,x,z)\rangle_H, \quad t\in[0,T], \,x\in E,\  z\in H,
\end{align*}
then 
$$
v(t,x)=V(t,x)
$$ 
and, thanks to \eqref{stima_min_loc_pen_cost}, the process 
\begin{equation}\label{un_eps_ott}
\bar u _s^\alpha:= \gamma^\alpha(s,X_s^{x,\bar u^\alpha},\nabla_x v(s, X_s^{x, \bar u^\alpha} )\nabla_x  X_s^{x, \bar u^\alpha}(-A)^{-\alpha}) \quad  \textup{for} \,\textup{-a.e.} \, \, s\in(0,T), \quad \P\textup{-a.s.},
\end{equation} 
belongs to $\mathcal U_2^\alpha$ and it  is optimal.
\end{theorem}

\proof
Notice that by Proposition \ref{prop:mild_sol_contr_prob_E}-$(ii)$, for any $u\in\mathcal U_2^\alpha$ there exists a unique mild solution $X^{t,x,u}$ to \eqref{Q=I} which satisfies \eqref{stima_sol_contr_gen2}. The proof is similar to the one of Theorem \ref{thm_rel_fond_opt_cont_str}. 
The main  difference consists in the fact that, for any given $u\in \mathcal U_2^\alpha$, we introduce a family of stopping times depending on the norm  $|\cdot|_\alpha$:
$$
\tau_n:=\inf\Big\{\tau\in[t,T]:\int_t^\tau|u_s|_\alpha^2ds>n\Big\}, \quad n\in\N.
$$
Then,  we set $u^n_\tau:=u_\tau\mathds1_{\tau\leq \tau_n}+u_0\mathds 1_{\tau >\tau_n}$, $u_0\in D((-A)^\alpha)$, and we introduce the process
\begin{align*}
W^n_\tau:=W_\tau+\int_t^\tau (-A)^\alpha u_s^nds.
\end{align*}
Afterwards, we  apply the Girsanov Theorem: writing
$u_s=(-A)^{-\alpha}(-A)^\alpha u_s$ 
in \eqref{Q=I}, we get that $X^{t,x,u^n}$ is mild solution to
\begin{align*}
\left\{
\begin{array}{ll}
dX_\tau=AX_\tau d\tau+F(X_\tau)d\tau+(-A)^{-\alpha}d W^n_\tau, & \tau\in[t,T], \vspace{1mm} \\
X_t=x\in E.
\end{array}
\right.
\end{align*} 
By \eqref{lipz_ham_pen_cost} in Lemma \ref{lem_stime_ham}, we see that Hypothesis \ref{ip-psiphi}-(ii) is verified by $\psi^\alpha$.
We conclude by arguing again as  in \cite[Proposition 4.1]{FT06} and in Theorem \ref{thm_rel_fond_opt_cont_str}.
\endproof


\subsection{The case \texorpdfstring{$Q=I$}{Q=I} with a  general running cost}

In this subsection we deal with the general controlled equation \eqref{Q=I} under Hypothesis \ref{2inghilterra} on the coefficients of the cost functional, and we consider control processes $u\in\mathcal U_2$. 
Unlike the two cases just treated, in this framework the HJB equation 
would not have 
 the structure of equation (\ref{Kolmo} ) since the Hamiltonian function 
would depend  on $\nabla v$, not only on the directional derivative $\nabla ^{(-A)^{-\alpha}}v$, see e.g. \cite{FabbriGozziSwiech},  formula (6.67) and the discussion related to formulas (4.278)-(4.279). Up to our knowledge, when $\phi$ in only continuous, the well posedness of 
 such an equation 
  is an open problem: in \cite{Ce} an  equation of this type is solved in mild sense with  Lipschitz type assumptions on the final datum $\phi$.
\newline For this reason,  we will not end up identifying the  value function \eqref{V} with the solution of the HJB equation, but instead we will approximate it. The following result will be used in the aforementioned approximation of the value function.
\begin{proposition}
\label{lips_est_gen_sol}
Assume that Hypothesis \ref{ip_forward} holds true. Let $t\in[0,T]$,  $x\in E$ and  $u,\tilde u\in \mathcal U_2$.
Then, 
\begin{align}
\label{lips_est_contr_gen_sol}
\sup_{\tau\in[t,T]}|X_\tau^{t,x,u}-X_\tau^{t,x,\tilde u}|_H^2
\leq\int_t^T|u_s-\tilde u_s|^2_H \,ds, \quad  \P\textup{-a.s.},
\end{align}
where $X^{t,x,u}$ and $X^{t,x,\tilde u}$ are  respectively the mild solutions to \eqref{Q=I} with control $u$ and $\tilde u$.
\end{proposition}
\proof
Let us set $L(\tau):=X_\tau^{t,x,u}-X_\tau^{t,x,\tilde u}$ and let us assume that $L(\tau)$ is a strict solution to
\begin{align*}
\left\{
\begin{array}{ll}
\frac d{d\tau}L(\tau)= AL(\tau)+F(X_\tau^{t,x,u})-F(X_\tau^{t,y,\tilde u})+u_\tau-\tilde u_\tau, & \tau\in[t,T], \vspace{1mm} \\
L(t)=0, \\
\end{array}
\right.
\end{align*}
otherwise we can use an approximation argument as in the proof of Proposition \ref{P:unificata}(ii). Then, the non-positivity of $A$, the dissipativity of $F$, the Cauchy-Schwartz inequality and the Young inequality give
\begin{align*}
\frac12\frac d{ds}|L(s)|^2_H\leq \langle u_s-\tilde u_s,L(s)\rangle_H
\leq \frac12 |u_s-\tilde u_s|_H^2+\frac12|L(s)|^2_H, \quad s\in[t,T], \quad \P\textup{-a.s.}\end{align*}
Integrating between $t$ and $\tau$ and  applying the Gronwall Lemma,  we get
\begin{align*}
|L(\tau)|^2_H\leq\int_t^\tau|u_s-\tilde u_s|_H^2ds, \quad \forall \tau\in [t,T],\quad  \P\textup{-a.s.},
\end{align*}
and we immediately deduce \eqref{lips_est_contr_gen_sol}.
\endproof

Thanks to Proposition \ref{lips_est_gen_sol} we deduce that, up to a subsequence, we can approximate $X^{t,x,u}$ in $H$ by means of mild solutions $X^{t,x,u^n}$ of problem \eqref{Q=I}, with $u$ replaced by $u^n$, where $(u_n)\subset \mathcal U_p$ satisfies $u^n\rightarrow u$ in $\mathcal U_p$. In the following Proposition we prove that a similar approximation holds true in $E$.

\begin{proposition}
\label{prop:conv_mild_sol_contr_E}
Let  Hypothesis \ref{ip_forward} holds true.   
Let  $t\in[0,T]$, $p\geq2$, and  set $p'$ be the conjugate exponent of $p$.
Assume that Hypothesis  \ref{ip_contr_E_H} holds true  for some positive   constant $\beta$ such that $p'\beta<1$. Let $u\in\mathcal U_p$ and  $(u^n)\subset \mathcal U_p$ be such that $u^n\rightarrow u$ in $\mathcal U_p$. Then, for any $x\in E$,
\begin{align}
\label{conv_mild_sol_contr_E}
\lim_{n\rightarrow+\infty}|X_\tau^{t,x,u^{k_n}}-X_\tau^{t,x,u}|_E=0, \quad \forall \tau\in[t,T], \quad \P\textup{-a.s.},
\end{align}
where $(u^{k_n})\subset (u^n)$ be such that $u^{k_n}_s\rightarrow u_s$ $\P$- a.s. for a.e. $s\in(t,T)$.
\end{proposition}
\proof
As usual, we limit ourselves to consider the case $t=0$. For any $n\in\N$, let us set $L^n:=X^{x,u^{k_n}}-X^{x,u}$, where $X^{x,u^{k_n}}$ and $X^{x,u}$ are mild solutions to \eqref{Q=I} with initial datum $x$ and control processes $u^{k_n}$ and $u$, respectively. Further, let us denote by $N$ the subset of $\Omega$ such that $\P(N)=0$ and $u^{k_n}_s\rightarrow u_s$ on $\Omega\setminus N$ for a.e. $s\in(t,T)$. Then, for any $t\in[0,T]$
\begin{align*}
L^n_t
= & \int_0^te^{(t-s)A}(F(X_s^{x,u^{k_n}})-F(X_s^{x,u}))ds
+ \int_0^te^{(t-s)A}(u^{k_n}_s-u_s)ds, \quad \P\textup{-a.s.}
\end{align*}
which gives
\begin{align*}
|L^n_t|_E
\leq \int_0^t|e^{(t-s)A}(F(X_s^{x,u^{k_n}})-F(X_s^{x,u}))|_Eds
+\int_0^t|e^{(t-s)A}(u^{k_n}_s-u_s)|_Eds =  :I^n_1(t)+I^n_2(t), \quad \P\textup{-a.s.}
\end{align*}
Let us estimate $I_1^n$ and $I_2^n$ separately. As far as $I^n_1$ is concerned, from the boundedness of $e^{tA}$ on $E$, Hypothesis \ref{ip_forward}-$4.$ and \eqref{stima_sol_contr_gen3}, it follows that
\begin{align*}
|e^{(t-s)A}(F(X_s^{x,u^{k_n}})-F(X_s^{x,u}))|_E<\infty, \quad \ s\in(0,T),
\end{align*}
on $\Omega\setminus N$. Further, from \eqref{stima_E_H} it follows that
\begin{align*}
|e^{(t-s)A}(F(X_s^{x,u^{k_n}})-F(X_s^{x,u}))|_E 
\leq c(t-s)^{-\beta}|F(X_s^{x,u^{k_n}})-F(X_s^{x,u})|_H, 
\end{align*}
on $\Omega\setminus N$, for any $s\in(0,T)$. Since $F$ is continuous on $H$, from \eqref{lips_est_contr_gen_sol} we infer that $|F(X_s^{x,u^n})-F(X_s^{x,u})|_H\rightarrow0$ on $\Omega\setminus N$ as $n\rightarrow+\infty$ for any $s\in(0,T)$. The dominated convergence theorem implies that $I^n_1\rightarrow0$ as $n\rightarrow+\infty$ on $\Omega\setminus N$.

Concerning $I_2^n$, from \eqref{stima_E_H} and  arguing as above we get
\begin{align*}
I_2^n(t)
\leq & c\int_0^t(t-s)^{-\beta}|u^{k_n}_s-u_s|_Hds
\leq c T^{-\beta+1/p'}\|u^{k_n}-u\|_{L^p(0,T;H)}\rightarrow0, \quad n\rightarrow+\infty,
\end{align*}
on $\Omega\setminus N$. This concludes the proof.
\endproof

\subsubsection{The approximate optimal control problem}
We will consider the Hamiltonian function $\psi^\alpha$ in \eqref{ham_pen_cost}
under  Hypothesis \ref{2inghilterra}. This prevents us to obtain directly estimates as those in Lemmas \ref{viaggione} and \ref{lem_stime_ham}, since we don't have the structure condition and the assumptions on $\ell$ are not sufficient to bound the term $(-A)^\alpha u$.
For this reason, for any $n\in\N$ we  introduce the function $\ell_n:[0,T]\times E\times H\longrightarrow \R$ defined by 
\begin{equation}\label{elln}
	\ell_n(s,x,u):= \ell(s,x,u)+\frac1n|(-A)^\alpha u|^2_H.  
\end{equation}
\begin{lemma}\label{L:est}
Let $A$ be an operator satisfying Hypothesis \ref{ip_forward}-(i). Then the function $\ell_n$  in \eqref{elln} satisfies the following conditions: for any $(t,x) \in[0,T]\times  H$, 
\begin{align*}
&0\leq \ell_n(t,x,u)\leq c_n(1+|u|_{\alpha}^2)\quad   u\in D((-A)^\alpha),\\
&\exists \,c_n, R>0: \,\,\ell_{n}(t,x,u)\geq c_n|u|_\alpha^2-cR^2, \quad  \ u\in D((-A)^\alpha).
\end{align*}	
\end{lemma}
\proof
The first inequality directly comes from \eqref{elln}. On the other hand, for any $t\in[0,T]$, $x\in E$ and $u\in D((-A)^\alpha)$, by Remark \ref{R: estl} we have
$$
\ell_n(t,x,u)\geq c\,(|u|_H^2-R^2)+\frac1n|(-A)^\alpha u|_H^2
\geq c_n|u|_\alpha^2-cR^2.
$$
\qed

For any $n\in \N$, we introduce the approximate Hamiltonian function
\begin{align}
\label{imballaggio}
\psi_n(t,x,z):=\inf_{u\in D((-A)^\alpha)}\left\{\ell_n(t,x,u)+\langle z,(-A)^\alpha u\rangle_H\right\}, \quad t\in[0,T], \ x\in E, \ z\in H. 
\end{align} 
Estimates in Lemma \ref{L:est} give the following result, which is analogous to Lemma \ref{lem_stime_ham}.
\begin{lemma}\label{lem_stime_ham_gen}
Let Hypothesis \ref{2inghilterra} be satisfied, and let $A$ be an operator satisfying Hypothesis \ref{ip_forward}-(i). Then, for any $n\in\N$,  the function $\psi_n$ in \eqref{imballaggio} is Borel measurable, and there exists a positive constant $C_n$ such that
\begin{align}
\label{4russia}
-C_n(1+|z|^2_H)\leq \psi_n(t,x,z)\leq \ell_n(t,x,u)+|z|_H| u|_\alpha, \quad t\in[0,T], \ x\in E, \ z\in H, \ u \in D((-A)^\alpha). 
\end{align}
Further, if the minimum in \eqref{imballaggio} is attained, it is attained in a ball of radius $C_n(1+|z|_H)$, i.e., 
\begin{align*}
\psi_n(t,x,z)=\inf_{u\in D((-A)^\alpha), |u|_\alpha\leq C_n(1+|z|_H)}\left\{\ell_n(t,x,u)+z(-A)^\alpha u\right\},
\quad t\in[0,T], \,\, x\in E,\  z\in H. 
\end{align*}
In particular,  
there exists a positive constant $C_n$ such that, for any $x_1,x_2\in E$,  $z_1,z_2\in H$, 
\begin{align}
\label{4podio}
|\psi_n(t,x_1,z_1)-\psi_n(t,x_2,z_2)|\leq C_n(|x_1-x_2|_E+|z_1-z_2|_H (1+|z_1|_H+|z_2|_H)), \quad t\in[0,T].
\end{align}
\end{lemma}
For any $n \in \N$, we introduce the approximate cost functional defined by
\begin{align}
\label{costoso2}
J_n(t,x,u)
& :=\mathbb E\Big[\int_t^T\ell_n(s,X_s^u,u_s)ds\Big]+\mathbb E[\Phi(X_T^u)],
\end{align}
and the associated approximated optimal control problem
\begin{align}\label{Vn}
V_n(t,x):=\inf_{u\in \mathcal U_2^\alpha}J_n(t,x,u), \qquad t\in[0,T], \,x\in E.
\end{align}

The HJB equation associated to the control problem (\ref{Vn}), related to the controlled state equation \eqref{Q=I},
  is given by
\begin{equation}
\left\{
\begin{array}
[c]{l}%
\frac{\partial v}{\partial t}(t,x)=-\call v\left(  t,x\right)
+\psi^n\left( t,x,v(t,x),\nabla^{(-A)^{-\alpha}}v(t,x)   \right)  ,\text{ \ \ \ \ }t\in\left[  0,T\right]
,\text{ }x\in E, \\
v(T,x)=\phi\left(  x\right),
\end{array}
\right.  \label{Hjb-n}%
\end{equation}
where $\psi^n$ is defined in (\ref{imballaggio}).
The HJB equation (\ref{Hjb-n}) is the analogous of (\ref{Hjb-alfa}) in Section \ref{sez-cost-alfa}. So again by Theorem \ref{teo_fey_kac}, its solution can be represented in terms of the solution $(X^{t,x}, Y^{n,t,x}, Z^{nt,x})$ of the forward-backward system
\begin{align}\label{fbsde-n-costo}
\left\{
\begin{array}{ll}
dY^{n,t,x}_\tau=-\psi_n(\tau,X^{t,x}_\tau,Z^{n,t,x}_\tau)d\tau+Z_\tau^{n,t,x} d W_\tau, & \tau\in[t,T],\vspace{1mm}\\
Y^{n,t,x}_T=\Phi(X^{t,x}_T), \vspace{1mm} \\
dX^{t,x}_\tau=AX^{t,x}_\tau d\tau+F(X^{t,x}_\tau)d\tau+(-A)^{-\alpha}d W_\tau, & \tau\in[t,T], \vspace{1mm} \\
X_t^{t,x}=x\in E,
\end{array}
\right.
\end{align}which is nothing else than the forward-backward system (\ref{fbsde-alfa}) with $\psi^n$ instead of $\psi^\alpha$.

We consider the following assumptions.

\begin{hypothesis}
\label{ip_impl_funct_n1}
For any $n\in\N$,  there exists a measurable function $\gamma_n:[0,T]\times E\times H\rightarrow D((-A)^\alpha)$ satisfying 
\begin{align}
\label{impl_funct_n_2}
\psi_n(t,x,z)=\ell_n(t,x,\gamma_n(t,x,z))+\langle z,(-A)^{\alpha}\gamma_n(t,x,z)\rangle_H, \quad t\in[0,T], \,x\in E,\  z\in H.
\end{align}
\end{hypothesis}

We state now  the analogous of Theorem \ref{thm_rel_fond_opt_cont_spec_cost} for the approximate optimal control problems \eqref{Vn}.
\begin{theorem}
\label{thm_rel_fond_opt_cont_gen_case}
Let Hypotheses \ref{ip_forward},  \ref{H: control1} hold true, and assume that Hypothesis \ref{ip_contr_E_H} holds true with $\beta<\frac12$. Let $X^{u,t,x}$ be the solution of equation \eqref{Q=I} and for any $n\in\N$, let $V_n$ be the function defined in \eqref{Vn}, and $v_n$ be the mild solution of the HJB equation \eqref{Hjb-n}.
%
Then, for any  $(t,x) \in[0,T]\times  E$ and $u\in \mathcal U_2^\alpha$, 
\begin{align*}
v_n(t,x)
= J_n(t,x,u)+\E\Big[\int_\tau^T\left(\psi_n(s,X_s^{t,x,u},Z_s^{n,t,x})-\ell_n(s, X_s^{t,x,u},Z_s^{n,t,x})-Z_s^{n,t,x}(-A)^\alpha u_s\right)ds\Big],
\end{align*}
where $(X^{t,x},Y^{n,t,x},Z^{n,t,x})$ is the solution to \eqref{fbsde-n-costo}.
In particular, $v_n(t,x)\leq V_n(t,x)$, for all $(t,x) \in[0,T]\times E$.
Finally, if Hypothesis \ref{ip_impl_funct_n1} holds true, then 
\begin{equation}\label{id_vnVn}
v_n(t,x)=V_n(t,x)
\end{equation} 
and, thanks to \eqref{stima_min_loc_pen_cost}, the process 
\begin{equation}\label{un_eps_ott2}
\bar u _s^n:= \gamma_n(s,X_s^{x,\bar u^n},\nabla_x v(s, X_s^{x, \bar u^n} )\nabla_x  X_s^{x, \bar u^n}(-A)^{-\alpha}) \quad  \textup{for} \,\textup{-a.e.} \, \, s\in(0,T), \quad \P\textup{-a.s.}
\end{equation} 
belongs to $\mathcal U_2^\alpha$ and it is optimal. 
\end{theorem}

\subsubsection{A characterization of the value function}

In the present section  we show that the value function $V$ of the optimal control problem  \eqref{V} can be approximated by the sequence $(v_n)$ of mild solutions to \eqref{Hjb-n}, that are identified with the approximated value functions $(V_n)$, see  formula \eqref{id_vnVn} in Theorem \ref{thm_rel_fond_opt_cont_gen_case}. As a byproduct, we deduce that the sequence $(\bar u^n)$ defined in \eqref{un_eps_ott2} is a minimizing sequence for \eqref{V},  and it is a  bounded sequence  in $\mathcal U_2$.

We start by introducing the Yosida approximations of $u\in\mathcal U_2$, namely a suitable sequence $(u_k)_{k\geq 1}\subset \mathcal U_2^\alpha$ which converges to $u$ in $\mathcal U_2$. Since $\bar u^n\in\mathcal U_2^\alpha$ for any $n\in\N$, this would allow  to approximate $V(t,x)$ in terms of $J(t,\,x,\bar u^n)$.

\begin{definition}\label{D:ueps_uespn}
For any $\varepsilon>0$, 
\begin{itemize}
	\item [(i)]  we denote by $u_\varepsilon\in\mathcal U_2$ any admissible control such that $J(t,x,u_\varepsilon)\leq V(t,x)+\varepsilon$.
\item [(ii)] we denote by  $u_{\varepsilon,k}$   the Yosida approximations of $u_\varepsilon$, i.e.,   
$$
u_{\varepsilon,k}(t,\omega): =
\begin{cases}
kR(k,A)\left(u_\varepsilon(t,\omega)\right), & \textrm{if}\,\,\,u_\varepsilon(t,\omega) \textrm{ is well defined}, \vspace{1mm} \\
0 & \textrm{otherwise}.
\end{cases}
$$
\end{itemize}
\end{definition}
\begin{lemma}\label{L:propuepsn}
Let $A$ be an operator satisfying Hypothesis \ref{ip_forward}-(i). Let $\varepsilon>0$ and $u_\varepsilon$, $u_{\varepsilon,k}$, with $k\in \N$, be the processes introduced in Definition \ref{D:ueps_uespn}. Then, for any $k\in\N$,  $\delta\in[0,1]$ and  $ (t,\omega)\in [0,T]\times \Omega$, 
	\begin{align}
&|u_{\varepsilon,k}(t,\omega)|_H\leq C|u_\varepsilon(t,\omega)|_H,\label{uepsnEst1}\\
&|(-A)^\delta u_{\varepsilon,k}(t,\omega)|_H\leq c_\delta k^\delta|u_\varepsilon(t,\omega)|_H,   \label{uepsnEst2}
\end{align}
for some positive constants $C$, $c_\delta$ not depending neither on $k$ nor on $u_\varepsilon$. In particular,  $u_{\varepsilon, k}\in\mathcal U_2^\delta$, $u_{\varepsilon,k}\rightarrow u_\varepsilon $ $\P$-a.s., a.e. in $[t,T]$ as $k\rightarrow+\infty$, and 
\begin{equation}
u_{\varepsilon,k}\rightarrow u_\varepsilon \,\,\textup{ in } \mathcal U_2\quad \textup{ as } n\rightarrow+\infty.\label{conv_uepsn}
\end{equation}
\end{lemma}

\proof
Estimate \eqref{uepsnEst1} directly follows from the properties of $R(k,A)$. Further, the fact that $u_{\varepsilon,k}\rightarrow u_\varepsilon $ $\P$-a.s., a.e. in $[t,T]$ as $k\rightarrow+\infty$ follows from the properties of Yosida approximations. Then, convergence \eqref{conv_uepsn} follows from the dominated convergence theorem, Finally, it easily follows that
\begin{align*}
|(-A)u_{\varepsilon,k}(t,\omega)|_H\leq kC|u_\varepsilon(t,\omega)|_H, 
\end{align*}
for any $(t,\omega)\in [0,T]\times \Omega$ and any $k\in\N$, where $C$ is the same positive constant as in \eqref{uepsnEst1}. Interpolation estimates give \eqref{uepsnEst2}. 
\endproof

\begin{proposition}\label{claim}
Let Hypotheses \ref{ip_forward},  \ref{H: control1}, \ref{ip_contr_E_H} hold true.
Let $\varepsilon>0$ and let $u_\varepsilon$,  $u_{\varepsilon,n}$, with $n\in\N$, be the processes introduced in Definition \ref{D:ueps_uespn}, and let $J$, $J_n$ be respectively  the cost functionals in  \eqref{costoso}, \eqref{costoso2}. Then for any $(t,x)\in[0,T]\times  E$ we have 
	$$
J_n(t,x,u_{\varepsilon,n})\rightarrow J(t,x,u_\varepsilon) , \quad  n\rightarrow+\infty. 
$$
\end{proposition} 
\proof
Since $u_{\varepsilon,n}$ pointwise converges to $u$, a.e. in $(0,T)$, $\P$-a.s., from \eqref{conv_uepsn} in Lemma \ref{L:propuepsn} and Proposition \ref{prop:conv_mild_sol_contr_E} it follows that $X_\tau^{t,x,u_{\varepsilon,n}}\rightarrow X_\tau^{t,x,u_\varepsilon}$ $\P$-a.s. in $E$ as $n\rightarrow+\infty$ for any $\tau\in [t,T]$. By dominated convergence theorem we deduce that
\begin{align*}
\mathbb E[\Phi(X_T^{u_{\varepsilon,n}})]\rightarrow \mathbb E[\Phi(X_T^{u_{\varepsilon}})], \quad  n\rightarrow+\infty. 
\end{align*}
To estimate the convergence of the approximate running cost $\ell_n$ in \eqref{elln},  we  consider separately the two terms in \eqref{elln}. We stress that
\begin{align*}
&\mathbb E \Big[  \int_t^T|\ell(s,X_s^{u_{\varepsilon,n}},(u_{\varepsilon,n})_s)-\ell(s,X_s^{u_{\varepsilon}},(u_{\varepsilon})_s)|ds\Big] \\
&\leq  \mathbb E \Big[\int_t^T|\ell(s,X_s^{u_{\varepsilon,n}},(u_{\varepsilon,n})_s)-\ell(s,X_s^{u_{\varepsilon}},(u_{\varepsilon,n})_s)|ds\Big] +\mathbb E \Big[\int_t^T|\ell(s,X_s^{u_{\varepsilon}},(u_{\varepsilon,n})_s)-\ell(s,X_s^{u_{\varepsilon}},(u_{\varepsilon})_s)|ds\Big].
\end{align*}
Arguing as above, from Hypothesis \ref{2inghilterra}-(iii) and \eqref{conv_mild_sol_contr_E} we get
\begin{align*}
|\ell(s,X_s^{u_{\varepsilon,n}},(u_{\varepsilon,n})_s)-\ell(s,X_s^{u_{\varepsilon}},(u_{\varepsilon,n})_s)|
\leq &L|X_s^{u_{\varepsilon,n}}-X_s^{u_{\varepsilon}}|_E \rightarrow0, \quad  n\rightarrow+\infty, \quad \P\textup{-a.s.}
\end{align*}
Further, from \eqref{2belgio} in Hypothesis \ref{2inghilterra}-(ii) and \eqref{uepsnEst1} we infer that 
\begin{align*}
|\ell(s,X_s^{u_{\varepsilon,n}},(u_{\varepsilon,n})_s)-\ell(s,X_s^{u_{\varepsilon,n}}, (u_{\varepsilon})_s)|
\leq c(1+|(u_{\varepsilon})_s|_H^2),
\end{align*}
for any $s\in(t,T)$, $\P$-a.s. By dominated convergence theorem we get
\begin{align*}
\mathbb E \Big[\int_t^T|\ell(s,X_s^{u_{\varepsilon,n}},(u_{\varepsilon,n})_s)-\ell(s,X_s^{u_{\varepsilon}},(u_{\varepsilon,n})_s)|ds\Big] \rightarrow 0, \quad n\rightarrow+\infty.
\end{align*}
Moreover, the continuity of $\ell$ with respect to $u$ and the dominated convergence theorem give
\begin{align*}
\int_t^T|\ell(s,X_s^{u_{\varepsilon}},(u_{\varepsilon,n})_s)-\ell(s,X_s^{u_{\varepsilon}},(u_{\varepsilon})_s)| \,ds\rightarrow0,\ n\rightarrow+\infty. 
\end{align*}
Finally, since $\alpha\in(0,1/2)$, from \eqref{uepsnEst2} with $\delta=\alpha$ we have
\begin{align*}
\frac1n\mathbb E\Big[\int_t^T|(-A)^\alpha (u_{\varepsilon,n})_s|_H^2 \, ds\Big]
\leq Cn^{2\alpha-1}\|u_\varepsilon\|_{\mathcal U_2}^2\rightarrow0, \ n\rightarrow+\infty,
\end{align*}
 and this concludes the proof. 
\endproof

The following  theorem constitutes the main result of the section.
\begin{theorem}\label{T:convJu_n}
Let Hypotheses \ref{ip_forward},  \ref{H: control1},  \ref{ip_impl_funct_n1} hold true, and  assume that Hypothesis  \ref{ip_contr_E_H} holds true with $\beta > \frac12$.  For any $n \in \N$, let $\bar u^n$ and $v_n$ denote respectively the process in \eqref{un_eps_ott2}  and the mild solution to \eqref{Hjb-n}. Let $V$, $J$  be respectively the functions in \eqref{V}, \eqref{costoso}. Then, for any $(t,x)\in[0,T]\times E$,
\begin{align}
V(t,x)&=\lim_{n\rightarrow+\infty}
v_n(t,x)= \lim_{n\rightarrow+\infty} J(t,x,\bar u^n).
\label{conv_gen}
\end{align}
Moreover,  $(\bar u^n)$ is bounded in $\mathcal U_2$.
\end{theorem}
\proof
Let $(t,x)  \in[0,T]$ and   $\varepsilon>0$. For any  $n\in\N$,  let $u_\varepsilon$, $u_{\varepsilon,n}$ be the processes introduced in Definition \ref{D:ueps_uespn}. By Proposition \ref{claim},  there exists $\bar n\in\N$ such that 
$|J_n(t,x,u_{\varepsilon,n})-J(t,x,u_\varepsilon)|\leq\varepsilon$, for  any $n\geq \bar n,
$
which in turn gives
\begin{equation}\label{First}
J_{n}(t,x,u_{\varepsilon,{n}})\leq V(t,x)+2\varepsilon, \quad n\geq \bar n.
\end{equation}
Notice that, from the definitions 
of $V$, $V_n$ in \eqref{V}, \eqref{Vn}, it  follows that 
\begin{equation}\label{VlessVn}
V(t,x)\leq V_n(t,x).
\end{equation}
Further, from \eqref{uepsnEst2} with $\delta=\alpha$ we get that $u_{\varepsilon,n}\in\mathcal U_2^\alpha$, and therefore the definition of $V_n$ implies
\begin{equation}\label{Second}
V_{n}(t,x)\leq J_{n}(t,x,u_{\varepsilon, n}), \quad n\geq \bar n.
\end{equation}
 Then, collecting \eqref{First}, \eqref{VlessVn} and \eqref{Second}, 
$
V(t,x)\leq V_{ n}(t,x)\leq J_{n}(t,x,u_{\varepsilon,n})\leq V(t,x) + 2 \varepsilon$ for any $n\geq \bar n$.
Hence, $V(t,x)\leq \limsup_{n\rightarrow+\infty}V_n(t,x)\leq V(t,x)+2\varepsilon$, and the arbitrariness of $\varepsilon$ gives
\begin{equation}\label{conv_Vn}
\lim_{n\rightarrow+\infty}V_n(t,x)=V(t,x). 
\end{equation}
Then the first equality in \eqref{conv_gen} follows from \eqref{conv_Vn},  recalling that,  by Theorem \ref{thm_rel_fond_opt_cont_gen_case},  $V_n(t,x)=v_n(t,x)$ for any $n\in\N$.

On the other hand,  since $V_n(t,x)=J_n(t,x,\bar u^n)$ for any $n\in\N$,  
\begin{align*}
V(t,x)\leq J(t,x, \bar u^n) \leq J_n(t,x, \bar u^n) = V_{ n}(t,x),
\end{align*}
so that,  taking into account \eqref{conv_Vn},  the second equality in \eqref{conv_gen} follows. 

Finally, let us prove  that  $(\bar u^n)$ is bounded in $\mathcal U_2$. Assume by contradiction that there exists a subsequence $(u_{k_n})\subset (\bar u^n)$ such that $\|u_{k_n}\|^2_{\mathcal U_2}\geq n$ for any $n\in\N$. Then,
\begin{align*}
n\leq & \|u_{k_n}\|^2_{\mathcal U_2}
= \int_t^T\Big(\int_{\{|u_{k_n}(s)|_H \leq  R\}}|u_{k_n}(s)|^2d\mathbb P\Big)ds+\int_t^T\Big(\int_{\{|u_{k_n}(s)|_H > R\}}|u_{k_n}(s)|^2d\mathbb P\Big)ds \\
\leq & TR^2+\int_t^T\Big(\int_{\{|u_{k_n}|_H > R\}}|u_{k_n}(s)|^2d\mathbb P\Big)ds.
\end{align*}
On the other hand, since $\ell$ is nonnegative and satisfies \eqref{3francia}, 
\begin{align*}
&\mathbb E \Big[ \int_t^T\ell(s,X_s^{k_n},u_{k_n}(s))ds\Big] \\
&=   \int_t^T\Big(\int_{\{|u_{k_n}(s)|_H \leq R\}}\ell(s,X_s^{k_n},u_{k_n}(s))d\mathbb P\Big)ds+ \int_t^T\Big(\int_{\{|u_{k_n}(s)|_H > R\}}\ell(s,X_s^{k_n},u_{k_n}(s))d\mathbb P\Big)ds \\
&\geq  \int_t^T\Big(\int_{\{|u_{k_n}(s)|_H > R\}}|u_{k_n}(s)|^2d\mathbb P\Big)ds.
\end{align*}
Therefore
$\mathbb E[  \int_t^T\ell(s,X_s^{k_n},u_{k_n}(s))ds]
\geq n-TR^2$, 
which contradicts \eqref{conv_gen}.
\endproof


\begin{thebibliography}{10}

\bibitem{B}
J.~M. Bismut.
\newblock Martingales, the {M}alliavin calculus and hypoellipticity under
  general {H}\"ormander's conditions.
\newblock {\em Z. Wahrsch. Verw. Gebiete}, 56:469--505, 1981.

\bibitem{BriFu}
P.~Briand and F.~Confortola.
\newblock {BSDE}s with stochastic lipschitz condition and quadratic {PDE} s in
  {H}ilbert spaces.
\newblock {\em Stochastic Process. Appl.}, 118(5):818--838, 2008.

\bibitem{Ce}
S.~Cerrai.
\newblock {\em Second order PDE's in finite and infinite dimension, A
  probabilistic approach}, volume 1762 of {\em Lecture Notes in Mathematics}.
\newblock Springer-Verlag, Berlin, 2001.

\bibitem{dpzErgodicity}
G.~Da~Prato and J.~Zabczyk.
\newblock {\em Ergodicity for infinite-dimensional systems}, volume 229 of {\em
  London Mathematical Society Lecture Note Series}.
\newblock Cambridge University Press, Cambridge, 1996.

\bibitem{DP3}
G.~Da~Prato and J.~Zabczyk.
\newblock {\em Second order partial differential equations in {H}ilbert
  spaces.}, volume 293 of {\em London Mathematical Society Note Series}.
\newblock Cambridge University Press, Cambridge, 2002.

\bibitem{dpz92}
G.~Da~Prato and J.~Zabczyk.
\newblock {\em Stochastic equations in infinite dimensions}, volume 152 of {\em
  Encyclopedia of Mathematics and its Applications}.
\newblock Cambridge University Press, Cambridge, second edition, 2014.

\bibitem{BaoDeHu}
F.~Delbaen, Y.~Hu, and X.~Bao.
\newblock Backward {SDE}s with superquadratic growth.
\newblock {\em Probab. Theory Related Fields}, 150(1-2):145--192, 2011.

\bibitem{ElKaPengQuenez}
S.~El~Karoui, N.~Peng and M.~C. Quenez.
\newblock Backward stochastic differential equations in finance.
\newblock {\em Mathematical finance}, 7(1):1–71, 1997.

\bibitem{FabbriGozziSwiech}
G.~Fabbri, F.~Gozzi, and A.~{\'S}wi{\c{e}}ch.
\newblock {\em Stochastic optimal control in infinite dimensions: {D}ynamic
  programming and {HJB} equations\textnormal{, with Chapter 6 by M. Fuhrman and
  G. Tessitore}}.
\newblock Springer, 2017.

\bibitem{futeBismut}
M.~Fuhrman and G.~Tessitore.
\newblock The {B}ismut-{E}lworthy formula for backward {SDE}s and applications
  to nonlinear {K}olmogorov equations and control in infinite dimensional
  spaces.
\newblock {\em Stoch. Stoch. Rep.}, 74(1-2):429--464, 2002.

\bibitem{fute}
M.~Fuhrman and G.~Tessitore.
\newblock Nonlinear {K}olmogorov equations in infinite dimensional spaces: the
  backward stochastic differential equations.
\newblock {\em Ann. Probab.}, 30(3):1397--1465, 2002.

\bibitem{FT05}
M.~Fuhrman and G.~Tessitore.
\newblock Generalized {D}irectional {G}radients, {B}ackward {S}tochastic
  {D}ifferential {E}quations and {M}ild {S}olutions of {S}emilinear {P}arabolic
  {E}quations.
\newblock {\em Appl Math Optim.}, 51(3):279–332, 2005.

\bibitem{FT06}
Y.~Fuhrman, M.~Hu and G.~Tessitore.
\newblock On a class of stochastic optimal control problems related to {BSDE}s
  with quadratic growth.
\newblock {\em SIAM Journal on Control and Optimization.}, 45(4):1279--1296,
  2006.

\bibitem{Go1}
F.~Gozzi.
\newblock Regularity of solutions of second order {H}amilton-{J}acobi equations
  in {H}ilbert spaces and applications to a control problem.
\newblock 20, 1995.

\bibitem{Go2}
F.~Gozzi.
\newblock Global regular solutions of second order {H}amilton-{J}acobi
  equations in {H}ilbert spaces with locally {L}ipschitz nonlinearities.
\newblock 198, 1996.

\bibitem{Kob}
M.~Kobylanski.
\newblock Backward stochastic differential equations and partial differential
  equations with quadratic growth.
\newblock {\em Ann. Probab.}, 28(2):558--602, 2000.

\bibitem{LU95}
A.~Lunardi.
\newblock {\em Analytic semigroups and optimal regularity in parabolic
  problems}.
\newblock Modern Birkh\"{a}user Classics. Birkh\"{a}user/Springer Basel AG,
  Basel, 1995.
\newblock [2013 reprint of the 1995 original] [MR1329547].

\bibitem{MasBanEJP}
F.~Masiero.
\newblock Regularizing properties for transition semigroups and semilinear
  parabolic equations in {B}anach spaces.
\newblock {\em Electron. J. Probab.}, 12:no. 13, 387--419, 2007.

\bibitem{MasBan}
F.~Masiero.
\newblock Stochastic optimal control problems and parabolic equations in
  {B}anach spaces.
\newblock {\em SIAM J. Control Optim.}, 47(1):251--300, 2008.

\bibitem{Mas1}
F.~Masiero.
\newblock Hamilton {J}acobi {B}ellman equations in infinite dimensions with
  quadratic and superquadratic {H}amiltonian.
\newblock {\em Discrete and Continuous Dynamical Systems}, 32(1A):223--263,
  2012.

\bibitem{Mas-spa}
F.~Masiero.
\newblock A {B}ismut-{E}lworthy formula for quadratic {BSDE}s.
\newblock {\em Stochastic Process. Appl.}, 125(5):1945--1979, 2015.

\bibitem{MR}
F.~Masiero and A.~Richou.
\newblock {HJB} equations in infinite dimensions with locally lipschitz
  hamiltonian and unbounded terminal condition.
\newblock {\em J. Differential Equations}, 257(6):1989--2034, 2014.

\bibitem{PaPe92}
\'{E}. Pardoux and S.~Peng.
\newblock Backward stochastic differential equations and quasilinear parabolic
  partial differential equations.
\newblock In {\em Stochastic partial differential equations and their
  applications ({C}harlotte, {NC}, 1991)}, volume 176 of {\em Lect. Notes
  Control Inf. Sci.}, pages 200--217. Springer, Berlin, 1992.

\bibitem{Par}
K.~R. Parthasarathy.
\newblock {\em Probability measures on metric spaces}, volume~3 of {\em
  Probability and Mathematical Statistics}.
\newblock Academic Press, Inc., New York-London, 1967.

\bibitem{PZ}
S.~Peszat and J.~Zabczyk.
\newblock Strong {F}eller property and irreducibility for diffusions on
  {H}ilbert spaces.
\newblock {\em Ann. Probab.}, 23(1):157--172, 1995.

\bibitem{Ri}
A.~Richou.
\newblock Markovian quadratic and superquadratic {BSDE}s with an unbounded
  terminal condition.
\newblock {\em Stochastic Process. Appl.}, 122(9):3173--3208, 2012.

\end{thebibliography}
\end{document}